\documentclass[reqno,11pt]{amsart}
\usepackage{psfig, amsmath, amsfonts, amssymb, amsthm, amscd}

\DeclareMathSymbol{\leqslant}{\mathalpha}{AMSa}{"36} 
\DeclareMathSymbol{\geqslant}{\mathalpha}{AMSa}{"3E} 
\DeclareMathSymbol{\eset}{\mathalpha}{AMSb}{"3F}     
\renewcommand{\leq}{\;\leqslant\;}                   
\renewcommand{\geq}{\;\geqslant\;}                   
\newcommand{\dd}{\,\text{\rm d}}             
\DeclareMathOperator*{\union}{\bigcup}       
\DeclareMathOperator*{\inter}{\bigcap}       
\newcommand{\sumtwo}[2]{\sum_{\substack{#1 \\ #2}}} 
\newcommand{\limtwo}[2]{\lim_{\substack{#1 \\ #2}}}     

\makeatletter
\def\captionfont@{\footnotesize}
\def\captionheadfont@{\scshape}

\long\def\@makecaption#1#2{%
  \vspace{2mm}
  \setbox\@tempboxa\vbox{\color@setgroup
    \advance\hsize-6pc\noindent
    \captionfont@\captionheadfont@#1\@xp\@ifnotempty\@xp
        {\@cdr#2\@nil}{.\captionfont@\upshape\enspace#2}%
    \unskip\kern-6pc\par
    \global\setbox\@ne\lastbox\color@endgroup}%
  \ifhbox\@ne 
    \setbox\@ne\hbox{\unhbox\@ne\unskip\unskip\unpenalty\unkern}%
  \fi
  \ifdim\wd\@tempboxa=\z@ 
    \setbox\@ne\hbox to\columnwidth{\hss\kern-6pc\box\@ne\hss}%
  \else 
    \setbox\@ne\vbox{\unvbox\@tempboxa\parskip\z@skip
        \noindent\unhbox\@ne\advance\hsize-6pc\par}%
\fi
  \ifnum\@tempcnta<64 
    \addvspace\abovecaptionskip
    \moveright 3pc\box\@ne
  \else 
    \moveright 3pc\box\@ne
    \nobreak
    \vskip\belowcaptionskip
  \fi
\relax
}
\makeatother
\def\writefig#1 #2 #3 {\rlap{\kern #1 truecm
\raise #2 truecm \hbox{#3}}}
\def\figtext#1{\smash{\hbox{#1}}
\vspace{-5mm}}
\psfigurepath{.:./pictures}

\newtheorem{lem}{Lemma}[section]
\newtheorem{pro}{Proposition}[section]
\newtheorem{thm}{Theorem}[section]
\newtheorem{cor}{Corollary}[section]

\newcommand{\cA}{\ensuremath{\mathcal A}}
\newcommand{\cB}{\ensuremath{\mathcal B}}
\newcommand{\cC}{\ensuremath{\mathcal C}}
\newcommand{\cD}{\ensuremath{\mathcal D}}
\newcommand{\cE}{\ensuremath{\mathcal E}}
\newcommand{\cF}{\ensuremath{\mathcal F}}

\newcommand{\cH}{\ensuremath{\mathcal H}}

\newcommand{\cK}{\ensuremath{\mathcal K}}

\newcommand{\cM}{\ensuremath{\mathcal M}}

\newcommand{\cO}{\ensuremath{\mathcal O}}
\newcommand{\cP}{\ensuremath{\mathcal P}}

\newcommand{\cV}{\ensuremath{\mathcal V}}
\newcommand{\cW}{\ensuremath{\mathcal W}}
\newcommand{\cX}{\ensuremath{\mathcal X}}

\newcommand{\frB}{\ensuremath{\mathfrak B}}

\newcommand{\frG}{\ensuremath{\mathfrak G}}

\newcommand{\frS}{\ensuremath{\mathfrak S}}

\newcommand{\frs}{\ensuremath{\mathfrak s}}

\newcommand{\frw}{\ensuremath{\mathfrak w}}

\newcommand{\bbB}{{\ensuremath{\mathbb B}} }

\newcommand{\bbD}{{\ensuremath{\mathbb D}} }
\newcommand{\bbE}{{\ensuremath{\mathbb E}} }

\newcommand{\bbH}{{\ensuremath{\mathbb H}} }

\newcommand{\bbL}{{\ensuremath{\mathbb L}} }

\newcommand{\bbN}{{\ensuremath{\mathbb N}} }

\newcommand{\bbP}{{\ensuremath{\mathbb P}} }

\newcommand{\bbR}{{\ensuremath{\mathbb R}} }
\newcommand{\bbS}{{\ensuremath{\mathbb S}} }
\newcommand{\bbT}{{\ensuremath{\mathbb T}} }

\newcommand{\bbZ}{{\ensuremath{\mathbb Z}} }

\newcommand{\gb}{\beta}
\newcommand{\gga}{\gamma}            
\newcommand{\gd}{\delta}
\newcommand{\gep}{\varepsilon}       

\newcommand{\gr}{\rho}
\newcommand{\gz}{\zeta}

\newcommand{\gD}{\Delta}
\newcommand{\gk}{\kappa}
\newcommand{\go}{\omega}
\newcommand{\gO}{\Omega}
\newcommand{\gl}{\lambda}
\newcommand{\gL}{\Lambda}
\newcommand{\gs}{\sigma}
\newcommand{\gS}{\Sigma}

\newcommand{\Ham}{{\bf H}}         
\newcommand{\PF}{{\bf Z}}          

\newcommand{\Is}{\mu}              
\newcommand{\Perc}{\Phi}           
\newcommand{\Joint}{\bbP}          

\newcommand{\Tor}[1]{\bbT_{#1}}    
\newcommand{\Boxx}[1]{\bbD_{#1}}   
\newcommand{\uTor}{\widehat \bbT^d}         
\newcommand{\uBox}{\widehat\bbD^d}         
\newcommand{\rBox}{{\uBox_r}}
\newcommand{\rBoxTwoD}{{\widehat\bbD^2_r}}
\newcommand{\norm}[1]{\|#1\|}
\newcommand{\normI}[1]{\|#1\|_{{\scriptscriptstyle 1}}}
\newcommand{\normII}[1]{\|#1\|_{{\scriptscriptstyle 2}}}
\newcommand{\normsup}[1]{\|#1\|_{{\scriptscriptstyle\infty}}}
\newcommand{\bdf}{\eta}
\newcommand{\intST}{\cW_\gb}
\newcommand{\intSTbd}{\cW_{\gb,\bdf}}
\newcommand{\intSTbdf}[1]{\cW_{\gb,#1}}
\newcommand{\intSTbdmin}{\cW^\star_{\gb,\bdf}}
\newcommand{\nnb}[2]{{\bk{#1,#2}}}
\newcommand{\gbc}{{\gb_{\rm\scriptscriptstyle c}}}
\newcommand{\hw}{{\bdf_{\rm\scriptscriptstyle w}}}
\newcommand{\bigO}{O}
\newcommand{\smallo}{o}
\newcommand{\bk}[1]{\langle#1\rangle}
\newcommand{\bnd}{\partial}
\newcommand{\abs}[1]{\lvert#1\rvert}
\newcommand{\bigabs}[1]{\Bigl\lvert#1\Bigr\rvert}
\newcommand{\ra}{\rightarrow}
\newcommand{\setof}[2]{\{#1\,:\,#2\}}
\newcommand{\bigsetof}[2]{\Bigl\{#1\,:\,#2\Bigr\}}
\newcommand{\BoxxNM}{{\Boxx{N,M}}}
\newcommand{\BoxxrN}{{\Boxx{N,rN}}}
\newcommand{\dBoxx}[1]{\bbD^\star_{#1}}
\newcommand{\dBoxxrN}{{\dBoxx{N,rN}}}
\newcommand{\halfspace}{{\bbL^d}}
\newcommand{\dhalfspace}{{\bbL^{d}_\star}}
\newcommand{\theplane}{\gS}
\newcommand{\thedplane}{\gS^*}
\newcommand{\thewall}{{\theplane_N}}
\newcommand{\thedwall}{{\theplane_N^\star}}
\newcommand{\taubd}{\tau_{\scriptscriptstyle\rm bd}}
\newcommand{\Isbd}[2]{\Is^{\beta,\bdf}_{#2,#1}}
\newcommand{\dIsbd}{\Is^{\beta^*,\bdf^*}_{\dBoxxrN}}
\newcommand{\Isbdf}[3]{\Is^{\beta,#3}_{#2,#1}}
\newcommand{\PFbd}[2]{\PF^{\beta,\bdf}_{#2,#1}}
\newcommand{\PFbdf}[3]{\PF^{\beta,#3}_{#2,#1}}
\newcommand{\Fbd}[2]{F^{\beta,\bdf}_{#2,#1}}

\newcommand{\Ebdf}[4]{\bk{#4}^{\beta,#3}_{#2,#1}}

\newcommand{\sgp}{\Isbd{+}{\halfspace}}
\newcommand{\sgm}{\Isbd{-}{\halfspace}}
\newcommand{\weightB}[1]{q_{N}^{\gb^*,#1^*}}
\newcommand{\weightL}[1]{q_\dhalfspace^{\gb^*,#1^*}}
\newcommand{\weight}{q^{\gb^*}}
\newcommand{\vol}{\text{vol}}
\newcommand{\Pol}{\text{Pol}}
\newcommand{\dH}{d_{\bbH}}
\newcommand{\Ism}[2]{\Is^{\gb}_{\gL,#2,#1}}
\newcommand{\Ismfree}[1]{\Is^{\gb}_{\gL,#1}}
\newcommand{\Ismiv}[2]{\Is^{\gb}_{#2,#1}}
\newcommand{\PFIs}{\PF^{\gb}_{\gL,\overline{\gs},h}}
\newcommand{\PFIsfree}{\PF^{\gb}_{\gL,h}}
\newcommand{\lgm}[2]{\nu^{\gb}_{\gL,#1,#2}}
\newcommand{\lgmfree}[1]{\nu^{\gb}_{\gL,#1}}
\newcommand{\lgmiv}[2]{\nu^{\gb}_{#1,#2}}
\newcommand{\PFlg}{\PF^{\gb}_{\gL,\mu,\overline{n}}}
\newcommand{\PFlgfree}{\PF^{\gb}_{\gL,\mu}}

\def\1{\ifmmode {1\hskip -3pt \rm{I}} \else {\hbox {$1\hskip -3pt \rm{I}$}}\fi}
\newcommand{\BV}{{ {\rm BV}(\uTor,\{\pm1\}) }}         
\def\lra{\leftrightarrow}
\def\nlra{\not \leftrightarrow}
\def\vol{{\rm vol}}
\newcommand{\df}{\stackrel{\gD}{=}}
\newcommand{\sTor}[1]{\widehat{\bbT}^d_{#1}}   
\newcommand{\sBox}[1]{\widehat{\bbB}_{#1}} 
\newcommand{\dBox}[1]{{\bbB}_{#1}} 
\def\top{{\rm top}}
\def\bot{{\rm bot}}
\newcommand{\st}{\tau_{\gb}}
\newcommand{\lb}{\left(}
\newcommand{\rb}{\right)}
\newcommand{\lbr}{\left\{}
\newcommand{\rbr}{\right\}}
\newcommand{\la}{\left\langle}
\newcommand{\ran}[2]{\right\rangle^{#1}_{#2}} 
\newcommand{\IsN}{\Is_{N,-}^{\gb}}

\begin{document}

\title[\,]{Rigorous probabilistic analysis of equilibrium crystal shapes}
\author{T.~Bodineau}
\address{
Universit\'e Paris 7, D\'epartement de Math\'ematiques, Case 7012, 2 place
Jussieu, F-75251 Paris, France}
\email{Thierry.Bodineau\@@gauss.math.jussieu.fr}
\author{D.~Ioffe}
\address{
Faculty of Industrial Engineering, Technion, Haifa 32000, Israel}
\email{ieioffe\@@ie.technion.ac.il}
\author{Y.~Velenik}
\address{
Fachbereich Mathematik, Sekr. MA 7-4, TU-Berlin, Stra\ss e des 17. Juni 136,
D-10623 Berlin, Germany
}
\email{velenik\@@math.tu-berlin.de}
\date{\today}
\begin{abstract}
The rigorous microscopic theory of equilibrium crystal shapes has made
enormous progress during the last decade. We review here the main results
which have been obtained, both in two and higher dimensions. In
particular, we describe how the phenomenological Wulff and Winterbottom
constructions can be derived from the microscopic description provided by
the equilibrium statistical mechanics of lattice gases. We focus on the
main conceptual issues and describe the central ideas of the existing
approaches.
\end{abstract}
\maketitle

\tableofcontents

\part{Introduction}
\label{part_Introduction}
\setcounter{section}{0}
\section{Phenomenological Wulff construction}

\subsection{Equilibrium crystal shapes}
\begin{figure}[t]
\centerline{
\psfig{file=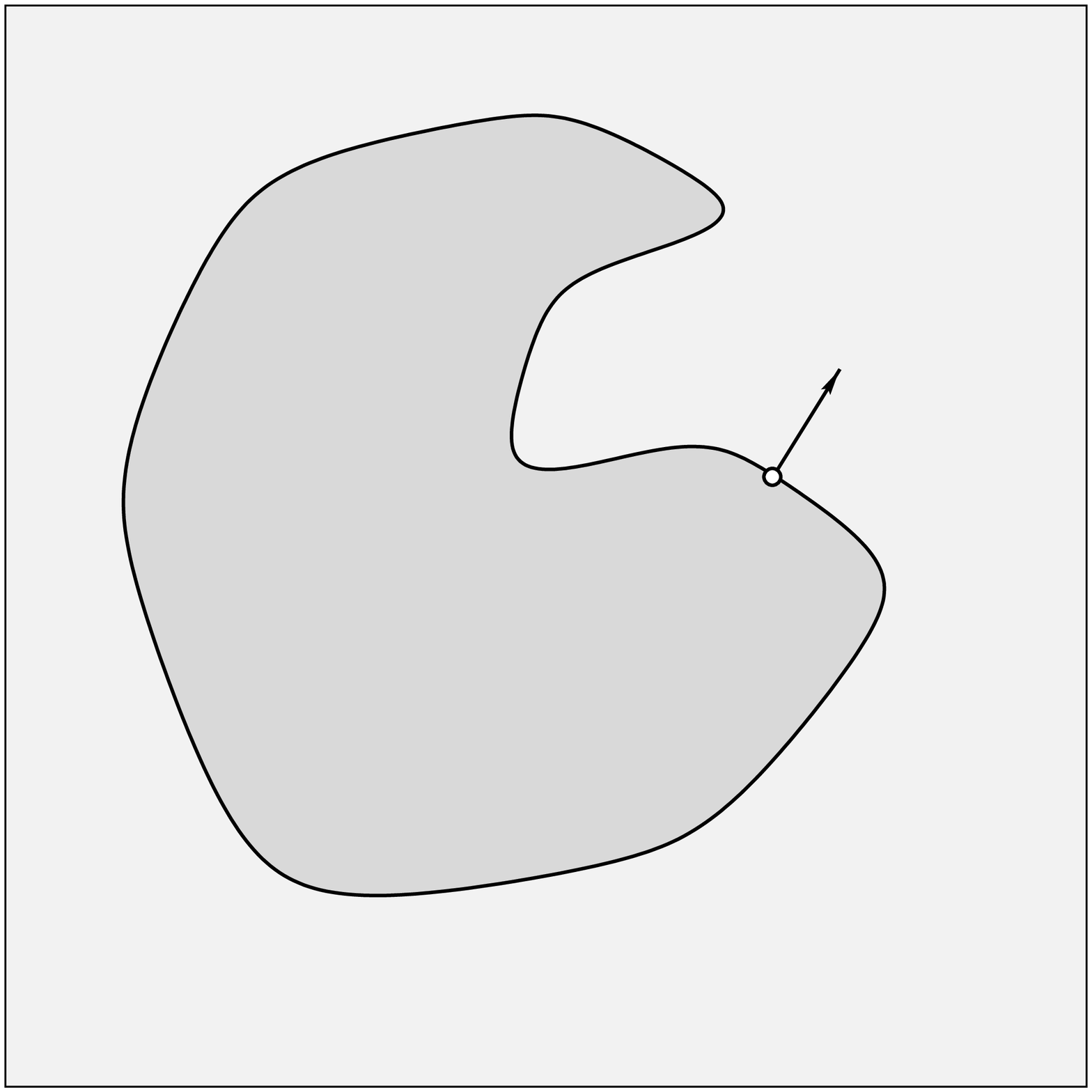,height=5cm}\hspace{1cm}
\raise 2.5 truecm \hbox{$
\cW_{\gb}\left( V\right)~=~\int_{\partial V}\st (\vec{n}_x )\,{\rm d}\cH^{(d-1)}_x
$}
}
\figtext{ 
\writefig       -2.36   3.10    {\footnotesize $x$} 
\writefig       -1.80   3.50    {\footnotesize $\vec{n}_x$} 
\writefig       -5.35   2.00    {\footnotesize $\partial V$} 
\writefig       -2.20   0.80    {\footnotesize \bf Vapor}
\writefig       -3.90   2.60    {\footnotesize \bf Crystal}
}
\caption{The free energy of the crystal-vapor interface is given by the 
integral of the anisotropic surface tension $\st$ over $\partial V$. 
$\cH^{(d-1)}$ is the $(d-1)$-dimensional Hausdorff measure.}
\label{wulff_tension}
\end{figure}

The phenomenological theory of equilibrated crystals  dates back at least to
the beginning of the century \cite{Wulff}. Suppose that two different 
thermodynamic phases (say crystal and its vapor) coexist at a certain 
temperature $T$. Assuming that the whole system is in equilibrium, in
particular that the volume $v$ of the crystalline phase is well defined,  what
could be said about the region this phase occupies? Of course, the issue cannot
be settled in the language of bulk free energies - these do not depend neither
on the shape, nor even on the prescribed volume $v$ of the crystal. Instead,
possible phase regions are quantified by the value of the free energy of the
crystal-vapor interface, or by the total surface tension between the crystal
and the vapor\footnote{In this review, our point of view is that of
mathematical physics; for an exposition of the problem from the viewpoint of
theoretical physics, we refer to \cite{RW} and references therein. }.
Equilibrium  shapes correspond, in this way, to the regions of minimal
interfacial energy. This is an isoperimetric-type problem: The surface tension
$\st$ (where, throughout the article, $\gb$ denotes the inverse temperature,
$\gb = 1/T$) is an anisotropic function of the local direction of the
interface. Thus, assuming that the crystal occupies a region  $V\subset\bbR^d$,
the corresponding contribution $\cW_{\gb}\left(V\right)$ to the free energy is
equal to the integral of $\st$ over the boundary  $\partial V$ of $V$
(Fig.~\ref{wulff_tension}).

The Wulff variational problem could then be formulated as follows:
\vskip 0.2cm
\noindent
$
\left({\rm {\bf WP}}\right)_v$\qquad\qquad $\cW_{\gb}
\left( V\right)~\longrightarrow~{\rm min}
\qquad\qquad{\rm Given}:\ {\rm vol}(V)~=~v
$
\vskip 0.2cm
\noindent
As in the usual isoperimetric case $\left({\rm WP}\right)_v $ is scale 
invariant,
\begin{eqnarray*}
\forall s >0 , \qquad 
\cW_\gb \big( \partial (sV) \big) = s^{d-1}  \cW_\gb \big( \partial V \big).
\end{eqnarray*}
Consequently, any dilatation of an optimal solution is itself optimal, and
one really talks here in terms of optimal shapes.

\begin{figure}[t]
\centerline{
\psfig{file=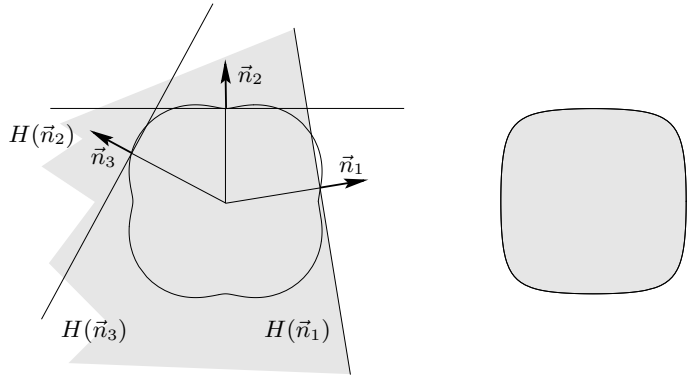,height=5cm}}
\figtext{ 
\writefig       -1.30   1.00    {\footnotesize $H(\vec{n}_1)$} 
\writefig       -4.7    3.60    {\footnotesize $H(\vec{n}_2)$} 
\writefig       -4.0    1.00    {\footnotesize $H(\vec{n}_3)$} 
\writefig       -0.30   3.20    {\footnotesize $\vec{n}_1$}
\writefig       -1.65   4.40    {\footnotesize $\vec{n}_2$}
\writefig       -3.60   3.30    {\footnotesize $\vec{n}_3$}
}
\caption{Function $\tau_\gb(\vec{n})$ (left)  with three half-spaces  $H(\vec{
n}_1 )$, $H(\vec{ n}_2 )$ and $H(\vec{ n}_3 )$ (for better visibility, only
$H(\vec{ n}_1 )$ has been shaded). The intersection  of {\bf all} such
half-spaces gives rise to the corresponding  Wulff shape (right).}
\label{fig_wulffconstruction}
\end{figure}

The canonical way to produce an optimal shape is given by the following 
Wulff construction (Fig.~\ref{fig_wulffconstruction}):
Define
\begin{eqnarray}
\label{Wulff shape}
\cK = \bigcap_{\vec{n} \in \bbS^{d-1}}  \left\{ x \in \bbR^d: \ 
x~\cdot ~\vec{n} \le \st (\vec{n})  \right\}\ \df\ 
\bigcap_{\vec{n} \in \bbS^{d-1}} H_{\gb}\left(\vec{n}\right) .
\end{eqnarray}
It would be convenient to normalize $\cK$ as
$$
\cK_1\ \df\ \sqrt[d]{\frac{1}{{\rm vol}(\cK )}}\cK .
$$
We refer to $\cK_1$ as to the normalized, or unit volume, Wulff shape.
The variational theory of $\left({\rm WP}\right)_v$, which we briefly address
in the subsequent subsection, states that any solution to 
$\left({\rm WP}\right)_v$ can be obtained by a shift of the 
corresponding dilatation $\cK_v\df \sqrt[d]{v}\cK_1$ of $\cK_1$.

\noindent
\subsection{Variational methods} 
\label{Variational methods}

The corresponding literature is rather
rich and diverse, here we merely attempt to facilitate the orientation of
the reader and to introduce some notations which will be useful in the 
sequel.

Since the half-spaces $H_{\gb}\left(\vec{n}\right)$ in \eqref{Wulff shape}
are convex, so is the Wulff shape $\cK$. Furthermore, in all the problems
we consider here, the surface tension $\st$ is bounded above and below,
\begin{equation}
\label{st_bound}
0~<~\min_{\vec{n} \in \bbS^{d-1}}\st (\vec{n})~\leq~
 \max_{\vec{n} \in \bbS^{d-1}}\st (\vec{n})~<~\infty .
\end{equation}
Accordingly, equilibrium crystal shapes are bounded and have non-empty
interiors, $0\in {\rm int}\big(\cK_v\big)$.

The fact that $\cK$ is optimal follows from the general Brunn-Minkowski
theory: Let $\tau_\gb^{**}$ be the support function of $\cK$, 
$\tau_\gb^{**}(x) = \sup \{ y \cdot x \ | \   y \in \cK \}$. Of course,
 if the homogeneous extension of $\st$
\begin{equation}
\label{tau tilde}
\st (\vec{x})\df \normII{ \vec{x}}\st \left( \frac{\vec{x}}{\normII{\vec{x}}}\right),
\end{equation}
 is
convex, then $\st$ and $\tau_\gb^{**}$ coincide. In general $\tau_\gb^{**}$
 is the convex lower-semicontinuous regularization of $\st$, in particular
 $\tau_\gb^{**}\leq\st$. Nevertheless, for the Wulff shape $\cK$,
\begin{eqnarray*}
\cW_{\gb}^{**}\left (\cK\right)~\df ~
\int_{\partial \cK} \tau_{\gb}^{**}(\vec{n}_x) \, d \cH^{(d-1)}_x ~=~
\int_{\partial \cK} \st (\vec{n}_x) \, d \cH^{(d-1)}_x .
\end{eqnarray*}
where, as before,  $\vec{n}_x$ is the outward normal to 
$\partial V$ in $x$ and $\cH^{(d-1)}$
is the $(d-1)$ dimensional Hausdorff measure in $\bbR^d$.

On the other hand, the action of the regularized functional 
$\cW_{\gb}^{**}$ could be extended to any compact set $V\subset\bbR^d$ in terms
of the mixed volume
$$
\cW_{\gb}^{**}\left( V\right)~=~\liminf_{\gep \to 0} \frac{1}{\gep} 
\left( \vol(V + \gep \cK) - \vol(V ) \right) ,
$$
the latter definition coincides with the integral definition of
$\cW_{\gb}^{**}$ for regular $V$. The Brunn-Minkowski inequality 
\cite{Schneider}
\begin{eqnarray*}
\label{Brunn-Minkowski}
\vol ( A+B) \geq 
\left( \vol(A)^{\frac{1}{d}} + \vol(B)^{\frac{1}{d}}
\right)^d \; ,
\end{eqnarray*}
implies that for any regular $V$ with ${\rm vol}\left( V\right) =
 {\rm vol}\left(\cK\right)$,
$$
\cW_{\gb}\left( V\right)~\geq ~
\cW_{\gb}^{**}\left( V\right)~
\geq ~ d \, \vol (\cK) = \cW_\gb (\cK) .
$$

Of course, we have been rather sloppy above, and we refer the reader to the
works  \cite{Taylor},~\cite{Fonseca} and   \cite{FonsecaMuller} for the
comprehensive discussion and results, including the history of the variational
Wulff problem. The language employed in the latter works is that of the
geometric measure theory, and we proceed with setting up some of the
corresponding notation which will also turn out to be useful for the
$\bbL_1$-approach to the microscopic justification of the Wulff construction,
as described  in Part~2 of this review. In the latter case,   the macroscopic
state of the system will be  determined by the value of an order parameter
which specifies the phase of the system. In the systems that we will consider,
the pure phases are characterized  by their averaged density, which are encoded
by two values $\rho_l(\gb)$ and $\rho_h(\gb)$,  for example $\rho_h$ for the
crystal and $\rho_l$ for the vapor. (In fact, we shall derive all the results
in the symmetrized spin language, in which case the two values will be $\pm
m^*(\gb)$, where $m^* (\gb )$ is the spontaneous magnetization (see Section~2)
at the inverse sub-critical temperature $\gb > \gb_c$).  For a given
temperature, it is convenient to replace this order parameter by a parameter
with values $\pm 1$. We suppose that the macroscopic region of $\bbR^d$ where
the system is  confined is the unit torus $\uTor = \left( \bbR / \bbZ
\right)^d$. The macroscopic  system is described by a function $v$ taking
values $\pm 1$  and the fact that $v_r = 1$ for some $r$ in $\uTor$ means that
locally at  $r$ the system is in equilibrium in the phase $m^*$.\\

\begin{figure}[t]
\centerline{
\psfig{file=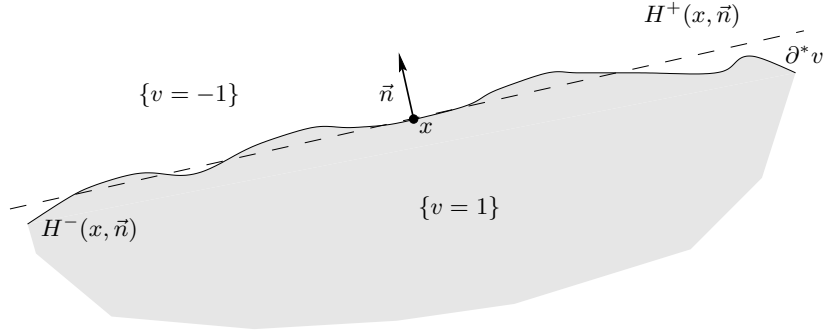,height=4cm}}
\figtext{ 
\writefig       -3.5    3.5     {\footnotesize $\{v=-1\}$} 
\writefig       0.22    2.00    {\footnotesize $\{v=1\}$}
\writefig       -0.3    3.5     {\footnotesize $\vec n$}
\writefig       3.22    4.55    {\footnotesize $H^+(x,\vec n)$} 
\writefig       -4.8    1.70    {\footnotesize $H^-(x,\vec n)$} 
\writefig       5.10    4.00    {\footnotesize $\partial^* v$} 
\writefig       0.23    3.07    {\footnotesize $x$} 
}
\caption{Measure theoretic unit normal to $\{ v = 1 \}$ at $x$} 
\label{fig_normext}
\end{figure}

For any measurable set $V$ in $\uTor$, the perimeter of $V$ is defined by
\begin{eqnarray}
\label{perimeter}
\cP(V) = \sup \left\{
\int_V {\rm div} \phi(x) \, dx \quad \big| \qquad 
\phi \in C^1( \uTor, \bbR^d), \  \ |\phi| \leq 1 
\right\} \; .
\end{eqnarray}
A function $v$ with values $\pm 1$ is said to be of bounded variation
in $\uTor$ if the perimeter of the set $\{ v = 1 \}$ is finite.
We denote by $\BV$ the set of functions of bounded variation in
$\uTor$ with values $\pm 1$ (see \cite{EG} for a review).
For any $v$ in $\BV$,
there exists a generalized notion of the boundary of $\{ v = 1 \}$
called reduced boundary and denoted by $\partial^* v$.
If $\{ v = 1 \}$ is a regular set, $\partial^* v$ coincides with the usual 
boundary $\partial v$.
Furthermore, a blow-up Theorem (see \cite{EG} p. 199) ensures that for all 
$x$ in $\partial^* v$ an approximate tangent plane can be defined locally.
This will imply the existence of a unit vector $\vec{n}_x$ called the
measure theoretic unit normal to $\{ v = 1 \}$ at $x$.
For any $x$ in $\bbR^d$ and any vector $\vec{n}$, we define the half spaces 
\begin{eqnarray*}
H^+(x, \vec{n}) & = & 
\{ y \in \bbR^d \ | \qquad (y-x) \cdot \vec{n} \geq 0 \} \; , \\
H^-(x, \vec{n}) & = & 
\{ y \in \bbR^d \ | \qquad (y-x) \cdot \vec{n} \leq 0 \} \; .
\end{eqnarray*}
Then for all $x$ in $\partial^* v$, there is a unit vector $\vec{n}_x$
such that
\begin{eqnarray*}
\lim_{r \to 0} \; \frac{1}{r^d} 
{\rm vol} \left( B(x,r) \inter \{ v = 1 \} \inter H^+(x, \vec{n}) \right) 
 & = & 0 \; , \\
\lim_{r \to 0} \; \frac{1}{r^d} 
{\rm vol} \left( B(x,r) \inter \{ v = - 1 \} \inter H^-(x, \vec{n}) \right) 
 & = & 0 \; , 
\end{eqnarray*}
where $B(x,r)$ is the ball of radius $r$ centered in $x$.
The previous property shows that the reduced boundary is not too wild
(see Fig. \ref{fig_normext}).
In fact, it is possible to prove that a set of finite perimeter has
``measure theoretically a $C^1$ boundary''. \\

The functional $\cW_\gb$ can be extended on 
$\bbL_1(\uTor ,[-\frac{1}{m^*},\frac{1}{m^*}])$ as follows
\begin{eqnarray}
\label{functional F}
\cW_\gb (v) =
\left\{
\begin{array}{l}
\int_{\partial^* v} \tau(\vec{n_x}) \, d \cH^{(d-1)}_x,
\qquad  {\rm if} \quad v \in \BV \; ,\\
\infty \; , \qquad \qquad \qquad \qquad  {\rm otherwise}.
\end{array}
\right.
\end{eqnarray}
Under the assumption that the homogeneous extension \eqref{tau tilde}
 of $\st$ 
is convex, a result by Ambrosio and Braides (see \cite{Ambrosio},  
Theorem 2.1) ensures that $\cW_\gb$ is lower 
semi-continuous with respect to $\bbL_1$ convergence. 
In certain cases (attractive interactions) the  convexity of $\st$  
can be derived from the properties of 
the corresponding microscopic system as will be explained later.

To any measurable subset $A$ of $\uTor$, we associate the function
$\1_A = 1_{A^c} - 1_A$ and simply write $\cW_\gb (A)=\cW_\gb (\1_A)$.
In this new setting, the isoperimetric problem is to find the minimizers of
\begin{eqnarray}
\label{variational}
\min \big\{ \cW_\gb (v) \ \big| \ v \in \BV, \qquad 
\big| \, \int_{\uTor} m^*  \, v_r \, dr \big|  \le m \big\},
\end{eqnarray}
where $m$ belongs to $]{\bar m}(\gb ) , m^*(\gb ) [$.
The parameter ${\bar m}$ is chosen such that the minima of
the variational problem above are translates of the set $\cK_m$
deduced from the Wulff shape $\cK$ by dilatation in order to satisfy the
volume constraint.
This restriction enables us to exclude pathological minimizers which
occur from the periodicity.
Nevertheless, notice that the precise shape or the uniqueness of the
minimizers of the variational problem will be irrelevant for the
microscopic derivation of the Wulff construction.

\subsection{Stability properties}
In two dimensions Wulff solutions to ${\rm (WP)}_v$ are stable in the 
metric of Hausdorff distance: let $V$ be a connected and 
simply connected subset
of $\bbR^2$ with a rectifiable boundary $\partial V$. Assume that 
${\rm Area} (V)\geq 1$. Then,
\begin{equation}
\label{stability}
\min_{x}{\rm d}_{\bbH}\left( V, x+\cK_1\right)~\leq ~c_1\sqrt{\intST (V) -
\intST (\cK_1 )} .
\end{equation}
This result has been established in \cite{DKS} as a generalization of the
classical Bonnesen inequality.

If $V$ consists of several connected and simply connected components,
 $V=\vee_{i=1}^{n} V_i$, and the total surface tension of $V$ is close to
the optimal,
$$
\intST (V)~=~\sum_{i=1}^{n}\intST (V_i )~\leq ~\intST (\cK_1 ) +\gep ,
$$
then, again assuming that ${\rm Area} (V) =\sum_{i=1}^{n}{\rm Area} (V_i) 
\geq 1$, an easy consequence of \eqref{stability} implies (see (2.9.7) and 
(2.9.8) in \cite{DKS}) that
 actually all but one components of $V$ are small, and that the only large
 component, say $V_1$, is close to a shift of $\cK_1$. Namely
$$
\sum_{i=2}^n{\rm Area} (V_i)~\leq~c_2\gep^2\qquad{\rm and}\qquad
\sum_{i=2}^n\intST (V_i)~\leq ~c_3\gep ,
$$
and $V_1$ satisfies \eqref{stability}.

These stability properties are indispensable for a sharp justification of
the phenomenological Wulff construction directly from the microscopic
 assumptions on the local inter-particle interactions (see 
Section~\ref{dima_structure} of
 Part~\ref{part_strongWulff}).

As far as we understand, stability properties of higher dimensional
isoperimetric problems are much less studied. Already in three dimensions
the Hausdorff distance is, of course, not an adequate measure of stability.
Trivial rate-free stability properties in $\bbL_1$ simply follow from 
the uniqueness of 
Wulff solutions and the compactness of BV-balls in $\bbL_1$. On a more 
qualitative side there are well studied stability properties in the 
class of convex sets \cite{Schneider} and, also, for sets with a smooth
boundary \cite{Hall}. We feel, however, that the statistical stability
 under the microscopic approximations in the problems we consider here
might be better than the impartial stability of the corresponding
variational problems. A result of this sort is supposed to appear in 
\cite{BodineauIoffeVelenik99}.

\subsection{Winterbottom problem}
The Wulff variational problem provides a description of an
equilibrium crystal shape deep inside a region filled with gas phase. If,
however, the spatial extent of the system is finite, it may happen that the
boundary of the surrounding vessel exhibits a preference toward the crystal
phase. In such a situation, the equilibrium state may not be given by the Wulff
shape anymore, but may have the crystal attached to the boundary. We discuss
briefly the simplest model of such an interaction between an equilibrium
crystal and an attractive substrate. Suppose, for simplicity, that our system
is contained in the half-space $H = \setof{x\in\bbR^d}{x(d)\geq 0}$; the
boundary of this half-space, the hyperplane $\frw = \setof{x\in\bbR^d}{x(d)= 0}$
represents the boundary of the vessel and is called the {\em wall}. We also
suppose to simplify the analysis, and because these assumptions will always be
satisfied, that $\tau_\gb(\vec n)=\tau_\gb(-\vec n)$, and that the homogeneous
extension of $\tau_\gb$ is convex\footnote{In the
models we consider in this paper, this is a consequence of FKG inequality.}.

To model the degree of attractiveness of the wall, we introduce a new
thermodynamical quantity, the {\em wall free energy} $\taubd(\gb,\bdf)$, which
depends on both the inverse temperature $\gb$ and the ``chemical structure'' of
the wall $\bdf$, and modify the free energy functional accordingly,
\begin{equation*}
\intSTbd(V) \df \intST(V) + (\taubd(\gb,\bdf) - \tau^*_\gb)\,\cH^{(d-1)}(\bnd
V\cap \frw)\,,
\end{equation*}
where $\tau^*_\gb \df \tau_\gb(\vec{e}_d)$, $\vec{e}_d\in\bbR^d$
with $\vec{e}_d(k)=\gd_{kd}$.
The wall free energy replaces therefore the surface tension $\tau_\gb$ along the
wall. At equilibrium, a thermodynamical stability argument shows that
$\taubd(\gb,\bdf)\leq \tau^*_\gb$ (this can also be proved in some microscopic
models, see Part~\ref{part_boundary}), so that this last term is always
non-positive. The new variational problem is

\vskip 0.2cm
\noindent
$\left({\rm {\bf WBP}}\right)_v$\qquad\qquad $\intSTbd
( V) \longrightarrow
{\rm min}$\qquad  Given: $V\subset H$, ${\rm vol}(V) = v$ .
\vskip 0.2cm
\noindent
It has first been studied in~\cite{Winterbottom67} and is called the
Winterbottom variational problem. Let us now discuss what its solution looks
like. It turns out that there are three cases to consider:
\begin{enumerate}
\item $\taubd(\gb,\bdf)= \tau^*_\gb$\hfill

\vspace*{1mm}
In this case, $\intSTbd(V) = \intST(V)$ and therefore the solution is the Wulff
shape associated to $\tau_\gb$. The equilibrium crystal is not attached to the
wall. This can happen even if {\it a priori} the chemical structure of the wall
is such that it is energetically favorable for the crystal to lay on the
wall, see Part~\ref{part_boundary} for a discussion from a microscopic point of view.

\vspace*{1mm}
\item $\abs{\taubd(\gb,\bdf)}< \tau^*_\gb$\hfill

\vspace*{1mm}
\begin{figure}[t]
\centerline{
\psfig{file=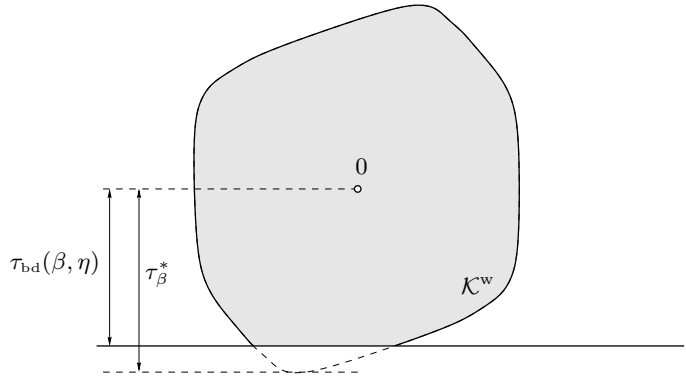,height=5cm}}
\figtext{ 
\writefig       -0.40   3.20    {\footnotesize $0$} 
\writefig        1.00   1.60    {\footnotesize $\cK^{\rm w}$} 
\writefig       -5.00   1.93    {\footnotesize $\taubd(\gb,\bdf)$} 
\writefig       -3.20   1.80    {\footnotesize $\tau^*_\gb$} 
}
\caption{The Winterbottom shape is obtained by taking the intersection between
the Wulff shape and the half-space $\{x(d)\geq-\taubd(\gb,\bdf)\}$, and rescaling
the obtained body.}
\label{fig_winterbottom}
\end{figure} 
Now the wall is really attractive for the crystal shape. The solution of the
variational problem is given by a suitably rescaled version of the following
set (see~Fig.~\ref{fig_winterbottom}),
\begin{equation*}
\cK^{\rm w} \df \cK \cap \setof{x\in\bbR^d}{x(d)\geq -\taubd(\gb,\bdf)}
\end{equation*}
so that the volume constraint is satisfied (notice that this variational problem
is still scale invariant); see~\cite{KoteckyPfister94} for a simple proof.
 
\vspace*{1mm}
\item $\taubd(\gb,\bdf) = -\tau^*_\gb$\hfill

\vspace*{1mm}
This is a somewhat pathological case. Indeed, the solution of the variational
problem is completely degenerate, the solution being unbounded. A minimizing
sequence is, for example,
\begin{equation*}
R_n = \setof{x\in H}{\abs{x(k)}\leq n,\, k=1,\dots,d-1,\, x(d)\leq n^{1-d}\,
v}\,.
\end{equation*}
As $n\ra\infty$, $R_n$ covers the whole wall with a film of vanishingly small
width; the limiting value of the surface free energy functional is $0$. This
describes the regime of so-called {\em complete wetting} where the wall so
strongly prefers the crystal that it wants to prevent any contact with the gas
phase.
\end{enumerate}

\subsection{Microscopic justification}
\label{ss_intro_pheno_sloppy}
 Microscopic models we consider here
are simple lattice gas type models (in the magnetic interpretation), which
are going to be defined precisely in the next section. The prototype situation
when the Wulff construction is thought to be recovered as a law of large 
numbers as the size of the microscopic system tends to infinity could be
loosely described as follows: Suppose that the particles of a certain substance
live on the vertices of the integer lattice $\bbZ^d$, so that each vertex of
 $\bbZ^d$ could be either occupied by a particle or remain vacant. Thus,
 various particle configurations $n$ could be labeled by  points of 
$\{ 0,1\}^{\bbZ^d}$, where one puts $n_i =1$ if there is a particle at
site $i\in\bbZ^d$, and $n_i =0$, otherwise. These random configurations 
 are sampled from a Gibbs distribution
 $\bbP$, which takes into account the assumptions on the microscopic 
interactions between the particles. The strength of the interaction
is quantified by the value $\gb =1/T$ of the inverse temperature; the larger
 $\gb$ (respectively the smaller the temperature $T$) is, the stronger is
the interaction. In many instances sufficiently low temperatures give rise
to two stable phases - the low density phase (which we call
 vapor) with an average particle density per site $\rho_l$ and the
high density phase (crystal) with a corresponding average density
 $\rho_h$, $0 <\rho_l <\rho_h <1$. 

Suppose now that all the particles are confined to a large finite volume vessel
$\gL_N\subset\bbZ^d$, where the subindex $N$ indicates the linear size
of $\gL_N$;  we put for simplicity $|\gL_N|=N^d$. 
Let us fix $\rho\in (\rho_l ,\rho_h )$ and ask what are the typical 
geometric properties of particle configurations $n$ under the conditional
 measure $\bbP\left(~\cdot ~\big| \sum_{i \in \gL_N}n_i =\rho N^d\right) $. 
In other words, we fix the total number of particles $\rho N^d$ in such a 
way that it falls in-between the two stable values $\rho_l N^d$ and
$\rho_h N^d$.

The prototype law of large numbers result we have in mind is schematically:
\begin{equation*}
\bbP\left(\;
\raise -2.4 truecm \hbox{\psfig{file=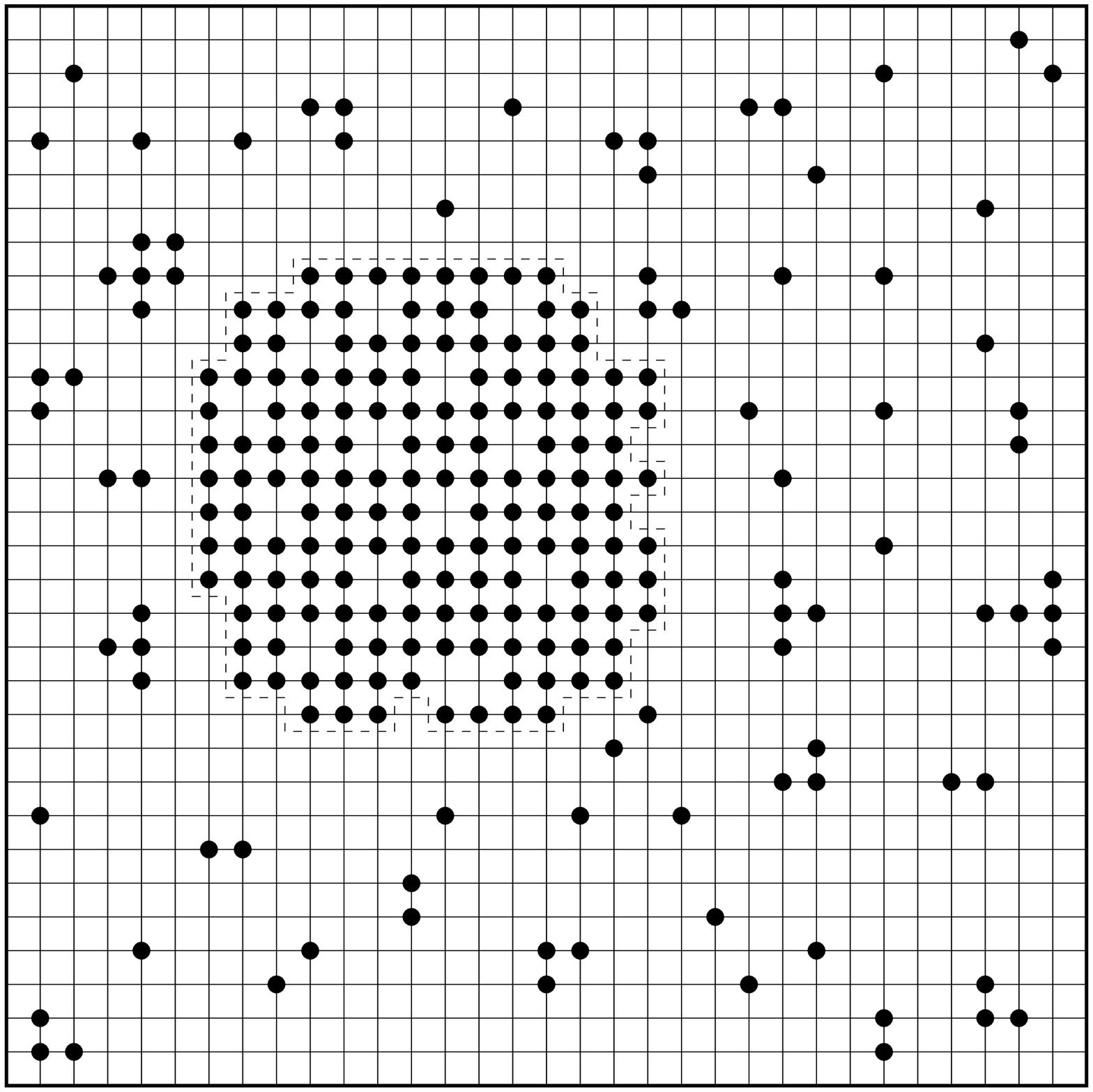,height=5cm}}\quad
\left| \;\;\sum_{i\in\gL_N}n_i =\rho N^d\right.\right)~\longrightarrow ~1\,.
\end{equation*}

Thus, with an overwhelming  $\bbP\left(~\cdot ~\big| \sum_{i \in \gL_N}n_i
=\rho N^d\right) $-probability particle configurations $n$ on $\gL_N$, $n\in
\{0,1\}^{\gL_N}$, obey the following phase segregation pattern: $\gL_N$ splits
into  two regions, $\gL_N =\gL_N^h\vee \gL_N^l$, where $\gL_N^h$ is occupied by
the high density phase, and, respectively, $\gL_N^l$ by the low density one.
The relative volume of $\gL_N^h$ can be recovered from the  canonical
constraint
$$
\rho_h\big|\gL_N^h\big|~+~\rho_l\big|\gL_N^l\big|~=~\rho N^d
$$
and the shape of $\gL_N^h$ is asymptotically Wulff.

There is a long way even towards making the above statement precise - 
we should define the microscopic models, quantify the notion of phases, in
particular of phases over finite volumes, and explain how the surface
 tension is produced in the large $N$ limit.

\section{Microscopic Models}\label{ssec_intro_models}
\setcounter{equation}{0}
\subsection{Models with finite-range ferromagnetic 2-body interactions}
We want to introduce mathematically precise realizations of the models discussed
in~subsection \ref{ss_intro_pheno_sloppy}. As described there, our interest
lies in models of lattice gases. For simplicity we restrict our attention to a
particular subclass of such models, which enjoy several nice properties, the
Ising models with finite-range ferromagnetic 2-body interactions.

We consider a family of random variables $n_i$, $i\in\bbZ^d$, taking values $0$
and $1$. Any site $i$ of the lattice $\bbZ^d$ is either occupied by a particle,
in which case $n_i=1$, or empty, in which case $n_i=0$. The random variables
$n_i$ are called {\em occupation numbers} and they completely describe a
configuration of the lattice gas. We consider a formal Hamiltonian of the form
\begin{equation*}
\tfrac 12 \sum_{i,j} K_{ij}\,n_i n_j\,,
\end{equation*}
the 2-body interactions are such that $K_{ij}=K_{\normI{j-i}}$, $K_{ij}\geq 0$
and $K_{ij}=0$ if $\normI{i-j}>r$, where $r$ is the {\em range} of the
interaction. We introduce two parameters, the {\em chemical potential} $\mu$
and the {\em inverse temperature} $\gb$, and set $\gL\Subset\bbZ^d$. The Gibbs
measure in $\gL$ with boundary condition $\overline{n}\in\{0,1\}^{\bbZ^d}$ is
the probability measure on $(\{0,1\}^{\bbZ^d}, \cA)$, with $\cA$ the usual
product $\gs$-field, defined by
\begin{equation*}
\lgm{\mu}{\overline{n}}(n) = 
   \begin{cases}
      \frac 1{\PFlg} \exp\bigl( \gb\mu\sum_{i\in\gL} n_i +
      \gb\displaystyle\sum_{\{i,j\}\cap\gL\neq\eset}K_{ij}\, n_i n_j
      \bigr) &
         \text{if $n_i = \overline{n}_i$, for all $i\not\in\gL$,} \\
      0 &
         \text{otherwise,}
   \end{cases}
\end{equation*}
where
\begin{equation*}
\PFlg = \sum_{N\geq 0} e^{\gb\mu N} \sumtwo{n\,:}{\sum_{i\in\gL}n_i=N}
e^{\gb\sum_{\{i,j\}\cap\gL\neq\eset} K_{ij}\, n_i n_j}\,.
\end{equation*}
Two types of boundary conditions are particularly relevant for us, the $\mathbf
1$ b.c., corresponding to setting $n\equiv 1$, and the $\mathbf 0$ b.c.,
$n\equiv 0$. We also need a different kind of boundary conditions: The Gibbs
measure in $\gL$ with {\em free} boundary conditions is the probability measure
on $(\{0,1\}^\gL, \cF_\gL)$ defined by
\begin{equation*}
\lgmfree{\mu}(n) = \frac 1{\PFlgfree} \exp\bigl( \gb\mu\sum_{i\in\gL} n_i +
\gb\sum_{\{i,j\}\subset\gL} K_{ij}\, n_i n_j
\bigr)\,.
\end{equation*}
These measures describe the lattice gas in the {\em Grand Canonical Ensemble},
in which the total number of particles, or equivalently the {\em density}
$\gr(n)=\frac 1{\abs\gL} \sum_{i\in\gL} n_i$, is not fixed. The description of
a gas in the {\em Canonical Ensemble} corresponds to the conditioned measure
\begin{equation*}
\lgm{\mu}{\overline{n}} (\,\cdot\,|\,\gr(n)=\widetilde\gr)\,,
\end{equation*}
with $\widetilde\gr\in {\rm Range}(\gr)$ (this measure is obviously independent of $\mu$). 
The existence of the Gibbs states
$\lgmiv{\mu}{\overline{n}}=\lim_{\gL\nearrow\bbZ^d} \lgm{\mu}{\overline{n}}$,
for $\overline{n}={\boldsymbol 0}$, $\boldsymbol 1$ or free, can be easily
proved using correlations inequalities; moreover, it is unique if $\mu\neq
-\tfrac12 \sum_j J_{0j}$. Restricting the chemical potential to the particular
line $\mu=-\tfrac12 \sum_j J_{0j}$, it can be proved that there exists a
critical value $\infty>\gbc>0$ such that
\begin{itemize}
\item For all $\gb<\gbc$, there is a unique Gibbs state and
$\lgmiv{\mu}{\overline{n}}(\gr) = 1/2$.
\item For all $\gb>\gbc$, $\gr_h(\gb) \equiv \lgmiv{\mu}{\mathbf 1}(\gr) > 1/2
> \lgmiv{\mu}{\mathbf 0}(\gr) \equiv \gr_l(\gb)$.
\end{itemize}

\medskip
It is rather convenient to work with another, equivalent, formulation of these
models, in which the symmetries present when $\mu=-\tfrac12\sum_j J_{0j}$ are more
transparent; this is the {\em magnetic interpretation}. To do this, we
introduce a new family of random variables $\gs_i$, $i\in\bbZ^d$, defined by
\begin{equation*}
\gs_i = 2n_i-1\,.
\end{equation*}
The random variables $\gs_i$ therefore take values in $\{-1,1\}$; $\gs_i$ is
called the {\em spin} at the site $i$. Expressed in these variables, the model
is defined through the following Gibbs measure in $\gL$ with boundary conditions
$\overline{\gs} \in \{-1,1\}^{\bbZ^d}$,
\begin{equation*}
\Ism{\boldsymbol h}{\overline{\gs}}(\gs) = 
   \begin{cases}
      \frac 1{\PFIs} \exp\bigl( \gb \displaystyle\sum_{i\in\gL} h_i\, \gs_i +
      \gb\displaystyle\sum_{\{i,j\}\cap\gL\neq\eset} J_{ij}\, \gs_i \gs_j
      \bigr) &
         \text{if $\gs_i = \overline{\gs}_i$, for all $i\not\in\gL$,} \\
      0 &
         \text{otherwise,}
   \end{cases}
\end{equation*}
where $h_i\in\bbR$ are called the {\em magnetic fields} and the {\em coupling
constants} $J_{ij}=J_{\normI{i-j}}$ satisfy $J_{ij}\geq 0$ and $J_{ij}= 0$ if
$\normI{i-j}>r$. A configuration $\gs$ such that $\gs_i = \overline{\gs}_i$, for
all $i\not\in\gL$, is said to be {\em compatible with b.c. $\overline\gs$ in
$\gL$}; the set of all such configurations is denoted by
$\gO_{\gL,\overline\gs}$. We are particularly interested in the $+$ and $-$
b.c. corresponding respectively to $\overline{\gs}\equiv 1$ and
$\overline{\gs}\equiv -1$. The Gibbs measure in $\gL$ with free b.c. is the
probability measure on $(\{-1,1\}^\gL,\cF_\gL)$ defined by
\begin{equation*}
\Ismfree{\boldsymbol h}(\gs) = \frac 1{\PFIsfree} \exp\bigl( \gb\sum_{i\in\gL}
h_i\, \gs_i + \gb\sum_{\{i,j\}\subset\gL} J_{ij}\, \gs_i \gs_j \bigr)\,.
\end{equation*}
Expected value w.r.t. these measures are denoted with brackets notations,
$\bk{\,\cdot\,}^{\gb}_{\gL,\overline\gs,\boldsymbol h}$, ...

In the magnetic formulation, the {\em Canonical Ensemble}
corresponds to fixing the value of the {\em magnetization} (density) $m(\gs) =
\frac1{\abs\gL} \sum_{i\in\gL}\gs_i$,
\begin{equation*}
\Ism{\boldsymbol h}{\overline{\gs}} (\,\cdot\,|\,m(\gs)=\widetilde m)\,,
\end{equation*}
where $\widetilde m\in {\rm Range}(m)$. If $h_i\equiv h$ for all $i$, then the
(infinite-volume) Gibbs states $\Ismiv{h}{\overline\gs}$ for $+$, $-$ and free
b.c. can be shown to exist; it is always unique when $h\neq 0$. The phase
transition statement takes now the following (simpler) form: There exists
$\infty>\gbc>0$ such that
\begin{itemize}
\item For all $\gb<\gbc$, the Gibbs state is unique and
$\bk{m}^\gb_{\overline{\gs},0} = 0$.
\item For all $\gb>\gbc$, $m^*(\gb) \equiv \bk{m}^\gb_{+,0} > 0
> \bk{m}^\gb_{-,0} = -m^*(\gb)$.
\end{itemize}
We will use the terminology {\em Ising models} to refer to the lattice gases in
the magnetic formulation. When $h=0$, we will generally omit it from the
notations.

\medskip
Ferromagnetic models are particularly well-suited for non-perturbative
analyses. Indeed, they enjoy several very useful qualitative properties, most
of which taking form of correlation inequalities. Of particular importance for
us are the following statements ($\gs_A \df \prod_{i\in A}\gs_i$):
\begin{align*}
\bk{\gs_A}^{\gb}_{\gL,\boldsymbol h} &\geq 0\,, \\
\bk{\gs_A\gs_B}^{\gb}_{\gL,\boldsymbol h} &\geq \bk{\gs_A}^{\gb}_{\gL,\boldsymbol h}
\bk{\gs_B}^{\gb}_{\gL,\boldsymbol h}\,, \\
\intertext{provided $h_i\geq 0$ for all $i$ (1st and 2nd Griffiths', or GKS,
inequalities~\cite{Griffiths72,KellySherman68}); also,}
\frac{\partial^2}{\partial h_i\partial h_j}\,\bk{\gs_k}^{\gb}_{\gL,\boldsymbol
h} &\leq 0\,,
\intertext{for all $i$, $j$ and $k$, provided $h_l\geq 0$ for all $l$ (GHS
inequalities~\cite{GriffithsHurstSherman70}); finally}
\bk{fg}^{\gb}_{\gL,\boldsymbol h} &\geq \bk{f}^{\gb}_{\gL,\boldsymbol h}
\bk{g}^{\gb}_{\gL,\boldsymbol h}\,,
\end{align*}
for any increasing\footnote{A function $f:\{-1,1\}^{\bbZ^d}\ra\bbR$ is {\em
increasing} if $f(\gs)\geq f(\gs')$ as soon as $\gs_i\geq \gs'_i$, for all $i$;
it is called {\em decreasing} if $-f$ is increasing.} functions $f$ and $g$,
and any ${\boldsymbol h}\in\bbR^\gL$ (FKG
inequality~\cite{FortuinKasteleynGinibre71}). Observe that any b.c. can be
obtained starting with free b.c. and applying suitable magnetic fields on the
spins on the inner boundary of $\gL$, where the {\em inner boundary} of a set
$A\subset\bbZ^d$ is defined as
\begin{equation*}
\bnd_{\rm in} A \df \setof{i \in A}{\exists j \not\in A,\, i \sim j}\,,
\end{equation*}
where $i \sim j$ means that $J_{i,j} \not = 0$. Similarly, we define the
{\em (exterior) boundary} of $A$ by
\begin{equation*}
\bnd A \df \setof{i \not\in A}{\exists j \in A,\, i \sim j}\,.
\end{equation*}
\subsection{2D nearest-neighbors ferromagnetic Ising model}
\label{ssec_2dIsing}
A particularly simple member of the above-mentioned class of models is
the two-dimensional nearest-neighbors Ising model, in which $J_{ij}=0$
if $i$ and $j$ are not nearest-neighbors, and $J_{ij}=1$ if they are.
This model has still
additional remarkable features. First, even though this only plays a very
marginal role in this review, it is the only one for which it is possible to
compute explicitly various quantities (free energy, surface tension,
correlations, ...). Of more importance for our purposes is the property of {\em
self-duality}\footnote{The fact that this model is {\em self}-dual is very
convenient, but is not required anywhere. What we need is to be able to control
precisely the dual of the model; for example, the Ising model on the hexagonal
lattice is not self-dual, but it would be possible to prove the same kind of
statements for this model as for the one on the square lattice.} that it
enjoys.

The nearest-neighbors model admit a geometric description in terms of very
simple objects, the contours. To define contours in the present context, it is
useful to introduce the notion of the dual of the lattice $\bbZ^2$. The {\em
dual lattice} is the set of dual sites
\begin{equation*}
\bbZ^2_\star = \setof{x\in\bbR^2}{x+(\tfrac12,\tfrac12) \in \bbZ^2}\,.
\end{equation*}
To each edge $e=\bk{x,y}$, $x,y\in\bbZ^2$, we associate a dual edge
$e^*$ connecting nearest-neighbors dual sites, which is the unique such edge
intersecting $e$ (as subset of $\bbR^2$).
 
Now, if we consider the Ising model in $\gL\Subset\bbZ^2$ with b.c.
$\overline\gs$, a configuration $\gs\in\gO_{\gL,\overline\gs}$ is entirely
determined by giving the following set of dual edges,
\begin{equation*}
\setof{e^*}{e^* \text{ dual to } e=\bk{i,j},\,
\{i,j\}\cap\gL\neq\eset,\,\gs_i\gs_j=-1}\,.
\end{equation*}
The maximal connected components of these dual edges, seen as closed line
segments in $\bbR^2$, are called {\em contours}.
We denote by $\boldsymbol\gga(\gs)$ the contours of the configuration $\gs$.
The {\em boundary} $\bnd\gga$ of a contour $\gga$ is the set of all dual sites
belonging to an odd number of the dual edges composing $\gga$. A contour is
said to be {\em closed} if $\bnd\gga=\eset$, otherwise it is {\em open}.

A set $\gL\Subset\bbZ^2$ is {\em simply connected} if
$\union_{i\in\gL}\setof{x\in\bbR^2}{\normsup{x-i}\leq 1/2}$ is a simply
connected subset of $\bbR^2$.

Given $\gL\subset\bbZ^2$, its {\em dual} is
$\gL^*=\setof{i\in\bbZ^2_\star}{\exists j\in\gL,\, \normsup{j-i}=1/2}$. A
family of contours is said to be {\em $\gL^*$-compatible} if they are disjoint
(as sets of bonds and sites) and are included in $\gL^*$. A family of contours
$\boldsymbol\gga$ is said to be {\em $(\gL,\overline\gs)$-compatible} if there
exists a configuration $\gs\in\gO_{\gL,\overline\gs}$ such that
$\boldsymbol\gga(\gs) = \boldsymbol\gga$. It is easy to show that for simply
connected $\gL$, $\gL^*$-compatibility of a family of closed contours is
equivalent to $(\gL,+)$-compatibility.

The measure $\Is^{\gb}_{\gL,\overline\gs}$ can be easily written in terms of
these objects; for any $\gs\in\gO_{\gL,\overline\gs}$,
\begin{equation}\label{eq_contours}
\Is^\gb_{\gL,\overline\gs}(\gs) = \frac 1{Z^\gb_{\overline\gs}(\gL)} \exp\{
-2\gb\sum_{\gga\in\boldsymbol\gga(\gs)} \abs\gga \}\,,
\end{equation}
where $\abs{\gga}$ is the number of edges in $\gga$ and
\begin{equation}\label{eq_dualLT}
Z^\gb_{\overline\gs}(\gL) = \sum_{\boldsymbol\gga\text{
$(\gL,\overline\gs)$-comp.}} \exp\{ -2\gb\sum_{\gga\in\boldsymbol\gga}
\abs\gga \}
\equiv \sum_{\boldsymbol\gga\text{
$(\gL,\overline\gs)$-comp.}} \prod_{\gga\in\boldsymbol\gga}w(\gga;\gb) \,.
\end{equation}

\medskip
We now discuss the property of self-duality.
Let $\gL\Subset\bbZ^2$ be simply connected. We consider the model at inverse
temperature $\gb^*$ in the box $\gL^*\Subset\bbZ^2_\star$, with free boundary
conditions. There exists another graphical representation for this model, the
{\em high-temperature} representation, which results from writing
\begin{equation*}
e^{\gb^*\gs_i\gs_j} =  \cosh\gb^* (1+\gs_i\gs_j\tanh\gb^*)\,,
\end{equation*}
opening all the brackets and expanding. After a simple summation over $\gs$,
this yields
\begin{align}\label{eq_dualHT}
Z^{\gb^*}_{\gL^*} = C(\gL)\,
\sum_{\boldsymbol\gga\text{ $\gL^*$-comp.}}
(\tanh\gb^*)^{\sum_{\gga\in\boldsymbol\gga} \abs\gga}
&\equiv C(\gL)\, \sum_{\boldsymbol\gga\text{ $\gL^*$-comp.}}
\prod_{\gga\in\boldsymbol\gga} w^*(\gga;\gb^*)\nonumber\\
&\equiv C(\gL)\,Z^{\gb^*}(\gL^*) \,,
\end{align}
where $C(\gL)$ is some constant which only depends on the set $\gL$. Setting
$\tanh\gb^*=e^{-2\gb}$, we see from \eqref{eq_dualLT} and \eqref{eq_dualHT}
that $Z^\gb_+(\gL)=Z^{\gb^*}(\gL^*)$, since $\gL$ is simply connected. In the
same way, we can expand the 2-point function, for example, and get the
following very useful identity
\begin{equation}\label{eq_randomline}
\bk{\gs_i\gs_j}^{\gb^*}_{\gL,+} = \sum_{\gl:i\ra j}q^{\gb^*}_{\gL^*}(\gl)\,,
\end{equation}
where the sum is over all open contours $\gl$ such that $\partial\gl=\{i,j\}$,
and
\begin{align*}
q^{\gb^*}_{\gL^*}(\gl) &= w^*(\gl;\gb^*) \, \frac {Z^{\gb^*}(\gL^*\,|\,\gl)}
{Z^{\gb^*}(\gL^*)}\,,\\
Z^{\gb^*}(\gL^*\,|\,\gl) &= \sumtwo{\boldsymbol\gga \text{
closed}}{(\boldsymbol\gga, \gl) \text{ $\gL^*$-comp.}}
\prod_{\gga\in\boldsymbol\gga} w^*(\gga;\gb^*)\,.\nonumber
\end{align*}
Identity \eqref{eq_randomline} is the so-called {\em random-line representation}
for the 2-point function of the Ising model, and plays a basic role in the
approach to the DKS theory of Part~\ref{part_strongWulff} (see
\cite{PfisterVelenik97, PfisterVelenik98} for much more details on this topic).
What is particularly useful is that the weights $q^{\gb^*}_{\gL^*}$, which we
have defined for an open contour, can be immediately extended to any family of
$\gL^*$-compatible contours (closed or open). In particular, if
$\boldsymbol\gga$ is a family of $\gL^*$-compatible {\em closed} contours, then
the following identity holds
\begin{equation*}
q^{\gb^*}_{\gL^*}(\boldsymbol\gga) = \Is^\gb_{\gL,+}(\boldsymbol\gga \subseteq
\boldsymbol\gga(\,\cdot\,))\,.
\end{equation*}
Applications and further results about the random-line representation are given
in Section~\ref{dima_skeletons} and in Part~\ref{part_boundary}.
The results stated
above also hold when the coupling constants are allowed to vary from edge to
edge, provided they remain ferromagnetic; if we denote by $J(e)$ the coupling
constant at edge $e$, then the duality relation takes the form
\begin{equation}\label{eq_dualitygeneral}
\tanh(\gb^*J^*(e^*)) = e^{-2\gb J(e)}\,.
\end{equation}

\subsection{Kac models}
In the original van der Waals Theory, the occurrence of phase
transitions is due to long range attractive forces between 
molecules.
In its statistical mechanics formulation, these forces are 
described by Kac potentials that depend on a positive 
scaling parameter $\gep$ which controls the strength and the
range of the potential (see \cite{KUH}).
The first probabilistic approach to this model was made
in the celebrated paper of Lebowitz and Penrose \cite{LebPenrose}.\\

In dimension $d$,
Ising systems with Kac potentials are defined by Gibbs measures  
with potentials depending on a scaling parameter  
$\gep>0$ 
   $$
\forall \, i,j \in \bbZ^d, \qquad
J^\gep_{i,j} = \gep^d J(\gep \normII{i-j} ) \; ,
   $$
and $J$ is a non-negative, smooth function supported by $[0,1]$
and normalized so that
   $$
\int_{\bbR^d} \! dr \, J(\normII{ r }) = 1.
   $$
The Gibbs measure on the domain $\gL$ is denoted by
$\Is^\gb_{\gep,\gL}$.
The constant $\gep$ will be so that the system
has finite but long range interaction. 
It is convenient to consider interaction parameters of the
form $\gep =2^{-m}$ ($m$ is typically assumed to be large but fixed).

This model bridges the finite range models and the mean field models.
In particular, if the range of the interaction, i.e. $\gep^{-1}$, 
is scaled proportionally to the number of spins
then the statistical properties of the 
system can be recovered from a mean field functional. 
In the true thermodynamic limit, when $\gep$ is kept fixed while
the number of spins goes to infinity,  
the behavior of the system cannot be described by 
the mean field continuum limit.
Nevertheless, by localizing in finite size regions it is possible to 
derive some informations from the mean field functional. 
This strategy was used to recover the phase diagram of the model and to
prove that it is arbitrarily close to the one of the mean field model when 
$\gep$ goes to 0.
More precisely, let us recall the following result which 
has been proven by Cassandro, Presutti
\cite{CP} and by Bovier, Zaharadnik \cite{BZ} 
(see also \cite{BP})
\begin{thm}
For any $\gb >1$, there is $\gep_0 >0$ such that for any $\gep$ smaller
than $\gep_0$ a phase transition occurs and there are at least 2
distinct pure phases $\Is_{\gep}^+$ and $\Is_{\gep}^-$.
\end{thm}

\noindent
If $\gb >1$, there is a breaking of symmetry and the spontaneous
magnetization is denoted by $\Is_{\gep}^+ (\gs_0) = m^*_\gep$.
Define $m^* = \lim_{\gep \to 0} m^*_\gep$.
This Theorem was proven via a renormalization procedure which
we shall describe in Subsection \ref{subsection Kac potentials}.

\subsection{Surface tension}
\label{subsection Surface tension}
\begin{figure}[t]
\centerline{
\psfig{file=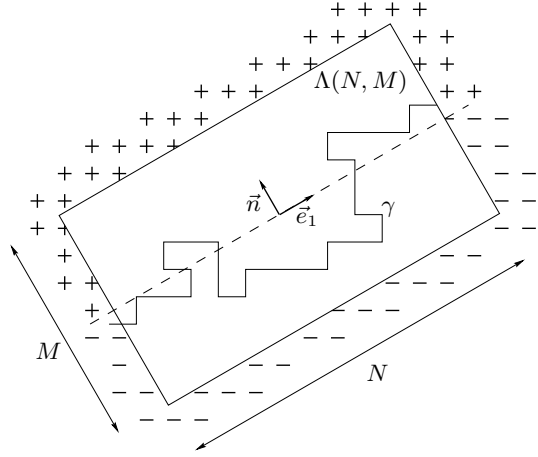,height=6cm}}
\figtext{ 
\writefig       -3.10   1.70    {\footnotesize $M$} 
\writefig        1.30   1.40    {\footnotesize $N$} 
\writefig        0.60   5.30    {\footnotesize $\gL(N,M)$} 
\writefig        0.36   3.50    {\footnotesize $\vec{e}_1$} 
\writefig       -0.30   3.66    {\footnotesize $\vec{n}$} 
\writefig        1.50   3.65    {\footnotesize $\gga$} 
}
\caption{Definition of the surface tension.}
\label{surface_tension}
\end{figure} 
We fix $\vec{n}$ a vector in $\bbS^{d-1}$ and consider an orthonormal
basis $(\vec{e}_1, \dots, \vec{e}_{d-1}, \vec{n})$.
Let ${\widehat \gL}(N,M)$ be the parallelepiped of $\bbR^d$ centered at 0 
with side length $N$ for the sides  parallel to $(\vec{e}_1, \dots, 
\vec{e}_{d-1})$ and side length $M$ for the sides  parallel to 
$\vec{n}$.
The microscopic counterpart of ${\widehat \gL}(N,M)$ is denoted by
${\gL}(N,M)$. 
The boundary $\partial \gL(N,M)$ is split into 2 sets
\begin{eqnarray*}
\partial^+_{\vec{n}} \gL(N,M) & = & \{ i \in \partial \gL(N,M) \; | \; 
\vec{i}.\vec{n} \geq 0\},\\
\partial^-_{\vec{n}} \gL(N,M) & = & \{ i \in \partial \gL(N,M) \; | \; 
\vec{i}.\vec{n} < 0\}.
\end{eqnarray*}
We fix the boundary conditions outside $\gL(N,M)$ to be equal to
1 on $\partial^+_{\vec{n}} \gL(N,M)$ and to $-1$ on 
$\partial^-_{\vec{n}} \gL(N,M)$.
The corresponding partition function on $\gL(N,M)$ is denoted by 
${\bf Z}^\gb_{\gL(N,M),\vec{n}, \pm}$.

Notice that any configuration $\gs$ contributing to the partition function
${\bf Z}^\gb_{\gL(N,M),\vec{n}, \pm}$ contains a $\pm$-contour $\gamma$ which 
crosses $\gL(N,M)$ under the ``averaged'' direction orthogonal to $\vec{n}$ 
(Fig.~\ref{surface_tension}). Such a contour is absent in the configurations
$\gs$ contributing to partition functions ${\bf Z}^\gb_{\gL(N,M), + }$ with 
pure boundary conditions on $\partial \gL(N,M)$. This contour represents 
the microscopic $\pm$-interface under the direction $\vec{n}$.

\medskip
\noindent
{\bf Definition :}
The surface tension in the direction $\vec{n} \in \bbS^{d-1}$  is 
defined\footnote{Notice that surface tension is sometimes defined with an
extra multiplicative factor $\frac{1}{\gb}$.}  by
\begin{equation}
\label{tau}
\tau_\gb (\vec{n}) = \lim_{N \to \infty} \; \lim_{M \to \infty} \;  
- {1 \over N^{d-1}}
\log { {\bf Z}^\gb_{\gL(N,M),\vec{n}, \pm } \over {\bf Z}^\gb_{\gL(N,M), + }}.
\end{equation}
\qed
\noindent

The proof the existence of the surface tension can be found in many papers (
\cite{Ab},~\cite{Pfister} to mention a few). A general approach has been
developed  by Messager, Miracle-Sole and Ruiz \cite{miracle}.  The core of
their proof is the sub-additivity of the sequence of finite-volume approximation
to $\tau_\gb (\vec{n})$ which is obtained by
means of FKG inequality. The proof is also valid for a wide range of models
like Ising models with  finite range interactions, Potts and SOS  models.
Furthermore, they showed that surface tension can be defined with 
parallelepipeds $\gL(N,M_N)$, where $M_N$ is a function of $N$ which diverges
as $N$ goes to infinity. More general domains can also be considered provided
they contain a parallelepiped of the type $\gL(N,M_N)$.

The convexity of the homogeneous extension of 
$ \tau_\gb$ (see (\ref{tau tilde})) is a consequence
of the pyramidal inequality  proven in Theorem 3 of \cite{miracle} :
Let $A_0, \dots, A_d$ be $d+1$ points of $\bbR^d$ and denote by
$(\gD_i)_{i \le d}$ the  simplex defined by these points.
Let $\vec{n}_i$ be the unit normal to $\gD_i$ and $|\gD_i|$ its area. Then,
the pyramidal inequality says
\begin{eqnarray*}
|\gD_0| \, \tau_\gb (\vec{n}_0) \leq 
\sum_{i =1}^d |\gD_i| \, \tau_\gb (\vec{n}_i).
\end{eqnarray*}

Note also that the homogeneous extension of $\tau_\gb$ is 
 continuous because it is locally 
bounded and convex.
Furthermore, $\tau_\gb$ is uniformly positive on $\bbS^{d-1}$.
This follows from the fact that the surface tension $\tau_\gb (\vec{n}_0)$ 
in the direction $\vec{n}_0 = (1,0, \dots ,0)$ is strictly positive as $\gb$ is 
larger than $\gb_c$ (see Lebowitz and Pfister \cite{LebPfister}).

\section{Scope of the theory}
\setcounter{equation}{0}

The key notion behind the attempts to give a rigorous meaning to the type  of
the phase segregation phenomena, which have been vaguely discussed in 
Subsection~\ref{ss_intro_pheno_sloppy}, is that  of {\bf renormalization} or
{\bf coarse graining}. The energy (probability) competes with the entropy
(number) of microscopic configuration in the corresponding energy shells.
Macroscopic quantities like surface tension are produced in the aftermath of
the entropy/energy cancelation, which is to say that in order to derive
large-$N$ ($N$-linear size of the system) asymptotics one should renormalize
appropriate microscopic objects. The appropriate objects here are, of course,
microscopic phase boundaries, which decouple between different ``large''
microscopic phase regions. These renormalization procedures could follow  two
different trends, depending on whether the renormalized (mesoscopic) 
structures keep track of the microscopic or macroscopic state of the system.

\subsection{Dobrushin-Koteck\'{y}-Shlosman Theory} The coarse graining
 of the DKS theory closely follows microscopic phase segregation patterns.
 Basic tools comprise a fluctuation analysis of the microscopic phase 
boundaries and sharp uniform local limit estimates over domains encircled
by such boundaries. Thus, the notion of finite volume phases is quantified 
 by the rate of the relaxation of the statistics of microscopic 
observables inside the microscopic phase regions  towards the corresponding
equilibrium values.

The theory has been developed using the low-temperature cluster expansions
in the seminal monograph \cite{DKS}. Our exposition in Part~3 is 
non-perturbative and follows the works \cite{Pfister}, \cite{I1},
\cite{I2}, \cite{PfisterVelenik97},
\cite{ScS2} and \cite{IS}. By and large the existing results are confined
to the simplest two-dimensional models (percolation and nearest neighbor 
Ising).

\subsection{$\bbL_1$-Theory} The renormalization approach of the 
$\bbL_1$-theory is, in a sense, opposite to that of DKS. In the latter case
 the principal coarse grained objects (skeletons, see Part~\ref{part_strongWulff}) 
are built
upon underlying families of large {\bf microscopic} contours. Such information
 is waved out in the $\bbL_1$-approach, and the basic renormalization objects
 here are the local (mesoscopic) order parameters or, in the spin language,
 locally averaged magnetization on various length scales.
 The idea is
that on sufficiently large scales local averages of the magnetization
are, with an overwhelming probability, close to one of the two equilibrium
values $\pm m^*$. Thus, under the renormalization, configurations are
characterized by their phase labels on different mesoscopic blocks. 
The objective of the ${\bbL}_1$-theory is to describe typical
mesoscopic magnetization profiles (or their phase labels) under a relaxed
canonical constraint of shell type. Unlike in the DKS case, the mesoscopic
phase labels are classified by their proximity to various  {\bf macroscopic}
 states. 
Combinatorial complexity of this
approximation is reduced by an exponential tightness property of the
mesoscopic phase labels (for a general claim of this sort 
see Theorem~\ref{thm Compactness}), which enables 
to restrict
 attention only to  $\bbL_1$-compact subsets of feasible macroscopic 
states, namely to the phase-sets of finite perimeter.
The core of the compactness estimates is based on the renormalization 
decoupling techniques introduced in \cite{Pisztora1} and 
 on the methods developed
to control the phase of small contours by
\cite{I2}, \cite{PfisterVelenik97}, \cite{ScS2} and \cite{IS}.
These techniques are robust enough to be applied on a renormalized
scale in any dimensions in a non perturbative setting.

Our exposition in this review is based on the work of \cite{Bo} with, though,
one  exception -- we specifically stress that all the relevant 
estimates of the $\bbL_1$-theory are obtained on appropriate
 {\bf finite} scales.
The validity of Lemma \ref{lem surface tension} up to the slab percolation
threshold follows from the results of \cite{CePi}.

\medskip
\subsection{Boundary Phenomena}
Parts~\ref{part_weakWulff} and \ref{part_strongWulff} provide a derivation of
Wulff construction from the basic principles of Equilibrium Statistical
Mechanics. Part~\ref{part_boundary} is concerned with a study of the effect of
the boundary conditions on the macroscopic geometry of the phase separation. In
particular, it is shown how the interaction with the boundary of the vessel can be
analyzed, and used to provide a derivation of Winterbottom construction. The
relationship between the macroscopic geometry in this case and the wetting
transition is also discussed. The presentation follows~\cite{PfisterVelenik97}
for the 2D case, and~\cite{BodineauIoffeVelenik99} for the higher-dimensional
ones.

\medskip
\subsection{Bibliographical review}
%
%
The rigorous 
investigation of the macroscopic geometry of phase separation under a
canonical constraint certainly started with two seminal papers of Minlos and
Sinai in 1967-68~\cite{MinlosSinai67, MinlosSinai68}. In these papers, the
authors considered nearest-neighbor very low temperature 
Ising models in arbitrary dimensions $d\geq
2$, even though they only wrote down the proof explicitly in the case $d=2$.
Their results
 could be roughly stated in the following way: At sufficiently low
temperatures, typical configurations of the Ising model in the exact canonical
ensemble over finite vessels of linear size $N$,  consist of a
 single large contour whose shape is ``nearly a square'', whereas the rest of the 
contours are small, that is at most of the order $\log N$. 
This is the picture of low temperature excitations of canonical ground states,
 and it has been treated by the authors as such. In particular, the entropic factor
has been frequently suppressed by the microscopic energy cost. 
However, exact asymptotic results
on the level of a microscopic justification of the Wulff construction depend, even
at very low but still non-zero temperatures, on a non-trivial entropy/energy
 competition, and, hence, could not be derived in this way.

Then there followed 15-20 years of a relative stagnation,
  the only
contributions to the area being confined to generalizations of~\cite{MinlosSinai67,
MinlosSinai68} to more complicated models~\cite{Kuroda82}. A popular interest to
the problem has been revived towards mid-eighties in the framework of an 
on-going mingle between probability and statistical mechanics \cite{Roberto},
\cite{FollmerOrt}, \cite{JoelRoberto}, \cite{CCSc}. 

 A breakthrough
occurred around 1989, when Dobrushin, Koteck\'y and Shlosman 
 found a way to derive the Wulff shape in a scaling limit of the low temperature
 2D Ising model. They found much more: Essentially the monograph \cite{DKS}
 sets up a comprehensive mathematical theory of phase segregation. This theory
 happened to be an intrinsically  probabilistic one. The DKS approach is, above
 all, to quantify the phenomenon of phase separation in terms of probabilistic
limit theorems and, accordingly, to study the probabilistic structures related
to the canonical states. Thus, in a sharp contrast with most of the 
preceding works, the ideology of \cite{DKS} has been from the start a very
 robust one and, actually, pertained to the whole of the phase transition
 region. It could be implemented, however, only at very low temperatures, since
 the authors used low temperature cluster expansions as the principal tool for
proving the corresponding probabilistic theorems. 

The ideas of \cite{DKS} did not wait long to inspire a wave of investigations,
even before the draft of the work started to circulate.  Two subsequent works
of a fundamental importance are \cite{Pfister}, where an alternative simplified
proof of parts of the DKS results has been given using  techniques, which are
specific to the 2D Ising model, like self-duality, and \cite{ACC}, where the
Wulff construction has been derived in the context of the 2D Bernoulli
percolation, but in a completely non-perturbative fashion, that is down to the
percolation threshold $1/2$. In both instances the exact canonical setting has
been  substituted by shell-type integral constraints, and, respectively,
softer integral type limit results have been used instead of the local
estimates of the original DKS theory. 

The results and techniques of \cite{ACC} and \cite{Pfister} have been combined
 with  profound renormalization ideas of \cite{Pisztora1} and 
lead to an extension
of this weak integral approach to the Wulff construction in the whole of
the 2D Ising phase coexistence limit \cite{I1}, \cite{I2}. Simpler proofs of
some of the basic estimates of these two works (e.g estimates in the phases
 of small contours or skeleton lower bounds) have been found in 
 \cite{SS}, \cite{ScS1}, and the integral version of the two-dimensional
 DKS theory has been essentially completed in \cite{PfisterVelenik97}, the estimates
of the latter work being already optimal along the lines of the integral 
approach. Furthermore, 
 Pfister and
Velenik ~\cite{PfisterVelenik96,
PfisterVelenik97} 
investigated the effect of boundary conditions, and in particular
studied the effect of an arbitrary boundary magnetic field, thus providing a
derivation of the Winterbottom construction. 

In spite of these successes, a non-perturbative treatment of
the full DKS theory was still out of reach, because a key ingredient was
missing: only rough estimates were available in the phase of small contours. By
proving a local limit theorem in the phase of small contours, Ioffe and
Schonmann were finally able to provide a non-perturbative version of the strong
Wulff theory \cite{IS}. The techniques of \cite{IS} are based on improved 
versions of asymptotic expansions in metastable cutoff phases developed in
\cite{ScS2}. 

In principle, the two-dimensional DKS theory should lead to exact expansions of
canonical partition functions up to zero-order terms. This, however, requires a
superb control over the statistical behavior of microscopic phase boundaries,
which is currently beyond the reach for the Ising model at moderately low 
temperatures. A certain progress, though, has  been reported  at very low
temperatures \cite{DH}, \cite{H} or either in the case of simplified models
\cite{HI}. Finally, it should be noted that at moderately low temperatures the
success of the DKS theory in two dimensions has been by and large confined to
the Ising and percolation models, and that there are serious technical and
possibly theoretical challenges to extend it to more general two-dimensional
models (see Section~\ref{dima_problems} for more on this).

On the other hand, as it has been communicated to us, 
an appropriate version of the low temperature DKS
theory (as originally developed in \cite{DKS}), should apply to any 
2-phase model in the realm of the Pirogov-Sinai theory~\cite{Senya}.

There is a strong interplay between dynamical properties of the Ising model and
its behavior in equilibrium : in absence of phase transition, the correlations
at equilibrium  are related to the exponential relaxation of the system;
instead as a  phase transition occurs, the dynamics is driven by the evolution
of droplets (nucleation, motion by mean curvature ...). We will not enter into
details and simply refer to the seminal paper on metastability by Schonmann and
Shlosman \cite{ScS2} and to the lecture notes by   Martinelli \cite{Martinelli}
(and references therein) for a survey of the recent works. Let us just mention
that, as far as phase coexistence is considered, many dynamical results are
only valid in  dimension 2 because of the absence of a precise description of
the equilibrium properties in higher dimensions.

\bigskip
If the 2D case was subject to rapid progress, the best results for higher
dimensions remained for a long time those of Minlos and Sinai.

The turning point of the latest developments should be traced back to the
seminal works by Pisztora~\cite{Pisztora1} and by Cassandro and Presutti~\cite{CP},
 where crucial renormalization decoupling estimates have been established 
in the case of the nearest neighbour Ising and, respectively, 
Kac interactions. 

The basic philosophy of the $\bbL_1$-approach has been originally developed 
in the 
works \cite{ABCP}, \cite{BCP}, \cite{BBBP}, \cite{BBP} in the 
context  of the  Ising systems with Kac potentials, and, in a 
less explicit way, elements and ideas of the theory already appeared 
in \cite{ACC}, \cite{Pisztora1}, \cite{I2} and \cite{PfisterVelenik97}. 

Using an embedding of the renormalized observables into a continuum setting, 
Alberti, Bellettini, Cassandro and Presutti
\cite{ABCP}, \cite{BCP} 
 emphasized the appropriateness of geometric measure theory setting,
introduced relevant analytic approximation procedures (see Subsection~2.6.1)
 and 
proved large deviation bounds for the appearance of a droplet of the minority
phase in a scaling limit when the size of the domain diverges not much faster
than the range of the Kac potentials.
In this scaling the system can be controlled by a continuum limit 
via  the $\Gamma$-convergence of functionals associated to the spins 
system \cite{ABCP} and 
by compactness arguments~\cite{BCP}.
 
The approach of \cite{ABCP} and \cite{BCP} has been  extended by Benois,
Bodineau, Butta and Presutti \cite{BBBP}, \cite{BBP} to the case when the range
of the  interaction remains fixed and does not change with the size of the
system.  The latter works are, already, structured in a way very similar to the
one  we expose here. Thus the main steps of \cite{BBBP} and \cite{BBP}
comprise  the coarse-graining of the rescaled magnetization profiles by the 
${\bbL_1}$-proximity to various continuum sets of finite perimeter, surgery
procedures to confine interfaces to tubes around the boundaries of such sets
and exponential tightness arguments to reduce the combinatorial  complexity of
the rescaled problem. The essential model-related input has been provided by
the decoupling estimates on the renormalized magnetization  \cite{CP},~
\cite{BZ} and by the result on the instanton structure of  Kac interfaces 
\cite{DOPT1, DOPT2}. The latter structure, however, yields only approximate 
bounds at each fixed finite interaction range. Consequently, the exact  (van
der Waals) surface tension could be recovered only when the range of the
interaction tends to infinity, that is only in the Lebowitz-Penrose  limit.
Nevertheless, at long but finite range interactions one could say that the
typical mesoscopic configurations  concentrate on  droplets with
$\bbL_1$-almost spherical shapes.

A complete picture of the higher-dimensional $\bbL_1$-Wulff construction has
been, for the first time, grasped and worked out in a recent remarkable  work
\cite{Cerf}, where the corresponding results have been established  in the
context  of the super-critical 3-dimensional Bernoulli bond percolation.  Using
novel and unusual renormalization procedures based on the decoupling results of
\cite{Pisztora1}, he has  essentially rediscovered all the main steps of the
$\bbL_1$-approach as  described above. The main turning point of \cite{Cerf}
was the  introduction of an alternative ingenious definition of the surface
tension which happened to be compatible with the setup of
$\bbL_1$-renormalization procedures \footnote{It should be noted, though, that
despite relative technical simplicity of this  observation,  the work
\cite{Cerf} most certainly prompted the  completion of the $\bbL_1$-theory  by
many years.}.   

The work of \cite{Cerf} triggered a wave of new investigations. In \cite{Bo}
his ideas on how to define and treat the surface tension have been combined
with an appropriate adjustment of the renormalization approach of  \cite{BBBP}
and \cite{BBP}, which lead to a relatively short proof of the  $\bbL_1$-Wulff
construction for the nearest neighbour Ising model in three and higher
dimensions and at sufficiently low temperatures. Most recently, a similar
construction has been established up to the FK slab  percolation  threshold in
\cite{CePi}. In the latter article new and important techniques have been
developed in order to go around mixed boundary conditions via  bulk relaxation
properties of the FK-measures.

Although the techniques of the $\bbL_1$-theory might look ``soft'' when
compared to the local limit setting of the DKS approach, one should bear in
mind that there is always a ``hard'' step needed to initialize the 
$\bbL_1$-machinery: The renormalized mesoscopic phase labels have to
possess sufficiently good decoupling properties. For the case of Kac models
the corresponding estimates have been established in \cite{CP}, \cite{BZ},
\cite{BMP},
 and in the case of percolation (including FK for the nearest neighbor
Ising model) models in dimension $d\geq 3$ in \cite{Pisztora1}, on which 
both \cite{Cerf},\cite{CePi} and \cite{Bo} rely in a fundamental way.

\medskip
Higher dimensional Winterbottom type shapes have been recovered  in the context
of effective interface models \cite{BolthausenIoffe}, \cite{BAD}, \cite{DGI},
\cite{DunlopMagnen98} following the original two-dimensional model defined and
studied in \cite{Dunlopetall}.

The results of these works have been also formulated in terms of $\bbL_1$
concentration properties, but the corresponding approach is quite different
from the one we expose here. Thus, the analysis of \cite{BolthausenIoffe}
heavily relies on specific properties of  Gaussian interactions.  It should
be noted, though, that, unlike in the nearest neighbour higher dimensional
Ising case,  there is better  insight into  the fluctuation and relaxation
properties of higher dimensional microscopic  interfaces
\cite{FunakiSpohn}, \cite{DGI}. On the other hand, the shapes produced by the
effective interface models are much less ``physical'', in particular the
equilibrium shapes are not scale invariant, and the corresponding surface
tension is not convex.

\part{$\bbL_1$-Theory}
\label{part_weakWulff}
\setcounter{section}{0}

On the macroscopic level the phenomenon of phase segregation is studied in
terms of concentration properties of the locally averaged magnetization. 
Statistical properties of the microscopic phase boundaries are waved out,
 and the backbone of the $\bbL_1$-theory are hard model-oriented 
renormalization estimates, which enable a sharp surface order analysis of
the mesoscopic  magnetization profiles. Example of such coarse graining 
procedures in the case of Kac, percolation and Ising models are given 
in Section~\ref{Examples of mesoscopic phase labels}.

\noindent
The averaging is performed on various mesoscopic scales:
\vskip 0.2cm
\noindent
{\bf Mesoscopic Notation.}  
 All the intermediate scales are of the 
form $2^k ,k\in\bbN$.  For any $M=2^k$ fixed 
we split the unit torus $\uTor$ into the disjoint union of
the corresponding mesoscopic boxes,
\begin{equation}
\label{ksplit}
\uTor ~=~\bigvee_{x\in\sTor{k}}\sBox{k}(x) ,
\end{equation}
where $\sTor{k}$ is the scaled  
 embedding of the discrete torus $\Tor{M} = \{1, \dots, M\}^d$ into
$\uTor$ as 
$$
\sTor{k} ~\df ~ \uTor\cap\left(\frac1M\Tor{M}\right) ,
$$
and, given $x\in\uTor$  the box $\sBox{k}(x)\subset\uTor$ 
is defined via
$$
\sBox{k}(x)~\df ~ x+\Big[-\!\frac{1}{2^{k+1}}~,~\frac{1}{2^{k+1}}~\Big)^d .
$$
Let us use $\cF_k$ to denote the (finite) algebra of the subsets of $\uTor$
 generated by the partition \eqref{ksplit} .
Given the size of the system $N = 2^n$, the local magnetization $\cM_k$ on
the $M=2^k\leq N$ scale is always an $\cF_{n-k}$-measurable function. This
notation should not be confusing: the subindex $k$ in $\cM_k$ measures the 
``coarseness'' of the mesoscopic magnetization profile. Thus, $\cM_0$ 
corresponds to the microscopic configuration, and $\cM_n$ identically
 equals to the averaged total magnetization. In general the local magnetization
$\cM_k$ is a piecewise constant function on $\uTor$ defined as
$$
\forall x \in \sTor{n-k}, \forall y \in \sBox{n-k} (x),
\qquad 
\cM_k(\gs,y) ~=~\frac{1}{M^d}\sum_{j\in\dBox{M}(2^nx )}\gs_j \, .
$$
Notice that the microscopic counterpart of the box $\sBox{n-k}(x)$
is the box $\dBox{M}(2^nx)$ of side length $M$ centered in $2^n x$.\\

We formulate all the results of Section~\ref{Results and the strategy of the proof}
for the nearest neighbor Ising model.
Along with the super-critical Bernoulli percolation this
is the only instance when a relatively complete $\bbL_1$-theory has been 
developed. 
In both instances, the validity of the $\bbL_1$-Theory hinges in a crucial way 
on the validity of Pisztora's coarse graining \cite{Pisztora1}, which is
by far the most profound model related fact employed.
Nevertheless, the approach itself is rather robust, and in 
subsequent Subsections we shall try to distinguish between specific model
dependent properties and more general results. In particular, 
compactness properties of local magnetization profiles are discussed in
Section~\ref{Coarse graining and mesoscopic phase labels} 
without any reference to specific models. Instead we briefly
indicate how the conditions of the corresponding general exponential 
tightness Theorem could be verified in several particular cases.

\section{Results and the strategy of the proof}
\label{Results and the strategy of the proof}
\setcounter{equation}{0}

\subsection{Main results}
\label{subsection Main results}

For simplicity, we restrict to the case of the torus $\Tor{N}$
and denote by $\Is_N$ the Gibbs measure with periodic boundary
conditions.

Define the  total magnetization ${\bf M}_{\Tor{N}}$ as 
$$
{\bf M}_{\Tor{N}}~\df~
\frac{1}{N^d} \sum_{i \in \Tor{N}} \gs_i .
$$
Let us define also the set $\frB_p$ as
\begin{align*}
\begin{split}
\frB_p = \{ \gb \; : \; {\text{Pisztora's coarse-graining hold for the
Ising model at inverse temperature $\gb$}} \}.
\end{split}
\end{align*}
We refer to the original article \cite{Pisztora1} and \cite{CePi} for the
precise relevant definitions (see also remark at the end of the Subsection 
\ref{subsection : Ising nearest neighbor}).
It is known that $\frB_p$ contains all except for at most countably
many points of the interval $]\tilde \gb_c,\infty[$, where $\tilde \gb_c$ 
is the so called slab percolation threshold, which is conjectured to 
coincide with $\gb_c$.\\

A compact way to state the main result of the $\bbL_1$-theory is:
\begin{thm}
\label{theo 1}
For any $\gb \in \frB_p$ and $m$ in $] {\bar m}, m^*[$ 
\begin{eqnarray*}
\lim_{N \to \infty} \;   \frac{1}{N^{d-1}}
\log \Is_{N} \big( \big| {\bf M}_{\Tor{N}} \big| \leq m \big) = - 
\cW_\gb (\cK_m ),
\end{eqnarray*}
where ${\bar m} ={\bar m}(\gb )$ and $\cK_m$ were defined in 
Subsection \ref{Variational methods}.
\end{thm}

\noindent
{\bf Remark.} 
The above Theorem has been established for $\gb \gg 1$ in \cite{Bo}. The only
additional ingredient required for an extension of the results of  the latter
paper to the whole of the temperature range $\tilde \gb_c$ was the validity of
the Lemma \ref{lem surface tension}. Such a statement happens to be highly
non-trivial, and it has been proven in \cite{CePi} along with  an alternative
derivation of the claim of Theorem \ref{theo 1}. \qed\\

\noindent
Theorem~\ref{theo 1} looks like a surface order large deviation principle.
 Such an appellation, however, would not help to explain the structure of
the underlying phenomena. In fact Theorem~\ref{theo 1} is essentially 
equivalent to
a seemingly  stronger statement on the macroscopic  geometry of the phase 
segregation of local magnetization profiles under the conditional measure 
$\Is_N\left( ~\cdot~\Big|\big|{\bf M}_{\Tor{N}} \big| \leq m \right)$:

For any function $v$ in  $\bbL^1( \uTor, [- \frac{1}{m^*}, \frac{1}{m^*}])$, 
the $\gd$-neighborhood of $v$ is denoted by $\cV(v,\gd)$
\begin{eqnarray*}
\cV(v,\gd) \df \left\{ v' \in \bbL^1 \big( \uTor, [- \frac{1}{m^*}, \frac{1}{m^*}] 
\big) 
\ \big| 
\qquad \int_{\uTor}  |v_x ' - v_x| \, dx \leq \gd \right\}.
\end{eqnarray*}
The $\bbL_1$-Theorem on the phase separation says that for
$\gb$ large enough with $\Is_{N} \left( \, . \; \Big| \, 
\big| {\bf M}_{\Tor{N}}\big| \le m 
\right)$-probability converging to 1, the function
$\cM_k$ is close to some translate of the Wulff shape 
$m^* \1_{\cK_m}$. 

\noindent
More precisely, fix a number $\nu <1/d$.
\begin{thm}
\label{theo 2}
For any $\gb \in \frB_p$ and $m$ in $]{\bar m}, m^*[$ the following holds:

\noindent
For every $\gd >0$, one can choose a scale $k_0 =k_0 (\gb ,\gd )$, such that
\begin{eqnarray*}
\lim_{N \to \infty} \;  \min_{k_0\leq k\leq \nu n} \; 
\Is_{N} \left(
\frac{\cM_k}{m^*} \in \bigcup_{x \in \uTor} \cV(\1_{\cK_m + x},\gd)
\; \Big| \, \big| {\bf M}_{\Tor{N}} \big| \leq m \right) = 1,
\end{eqnarray*}
where ${\bar m}$ and $\cK_m$ were defined in Subsection 
\ref{Variational methods}.
\end{thm}

The proofs of Theorems \ref{theo 1} and \ref{theo 2} are similar and
are divided into 2 steps. 
The first step amounts to prove a compactness Theorem and the second one
to derive precise logarithmic asymptotics.

\subsection{Exponential tightness}

Recall \cite{EG} that for any $a$ positive, the set
$$
K_a \df \big\{ v \in \BV \; |\quad  \cP (\{v=1\}) \leq a \big\},
$$
is compact with respect to convergence in $\bbL^1(\uTor)$.
\begin{pro}
\label{prop 1}
Let $\gb$ be in $\frB_p$.
Then there exists a constant $C(\gb)>0$ such that for all $\gd$ 
positive one can find $k_0(\gd)$
\begin{eqnarray*}
\forall a >0, \qquad 
\limsup_{N \to \infty} \;   \frac{1}{N^{d-1}}
\max_{k_0(\gd) \leq k\leq \nu n}
\log 
\Is_{N} 
\left( \frac{\cM_k}{m^*}  \in \cV(K_a,\gd)^c \right) 
\leq - C(\gb) \, a,
\end{eqnarray*}
where $\cV(K_a,\gd)$ is the $\gd$-neighborhood of $K_a$ in
$\bbL^1( \uTor, [-\frac{1}{m^*},\frac{1}{m^*}])$.
\end{pro}
This proposition tells us that only the configurations
close to the compact set $K_a$ have a contribution which is
of the surface order. 
This statement reduces the complexity of the problem : as
$K_a$ is compact, it is enough to derive the leading
terms in the logarithmic asymptotics  for the probability of a finite 
number of events.

In Section~\ref{Coarse graining and mesoscopic phase labels}, 
we prove that the analog of Proposition \ref{prop 1} holds for a broad 
class of models.

\subsection{Precise logarithmic asymptotics}

As the minimizers are known, it is sufficient to derive a lower
bound for configurations  concentrated close to $\cK_m$.
\begin{pro}
\label{prop 2}
Let $\gb$ be in $\frB_p$ and let $m$ be in $]{\bar m}, m^*[$
\begin{eqnarray*}
\liminf_{N \to \infty} \;   \frac{1}{N^{d-1}}
\min_{k_0(\gd) \leq k \leq \nu n}
\log \Is_{N} 
\left( \frac{\cM_k}{m^*}  \in \cV( \1_{\cK_m},\gd) \right) 
\geq - \cW_\gb(\cK_m ) - o(\gd) \, ,
\end{eqnarray*}
where the function $o(\cdot)$  depends only on $\gb$ and vanishes 
as $\gd$ goes to 0.
\end{pro}

According to proposition \ref{prop 1}, we will prove
the upper bound only for a restricted class of events
\begin{pro}
\label{prop 3}
Let $\gb$ be in $\frB_p$. Then for all $v$ in $\BV$ such that 
$\cW_\gb(v)$ is finite, one can choose $\gd_0 = \gd_0 (v)$,
such that uniformly in $\gd < \gd_0$
\begin{eqnarray*}
\limsup_{N \to \infty} \; \frac{1}{N^{d-1}}
\max_{k_0(\gd)  \leq k\leq \nu n}
\log \Is_{N}  
\left( \frac{\cM_k}{m^*}  \in \cV(v,\gd) \right) \leq - \cW_\gb(v)
+ o(\gd) \, .
\end{eqnarray*}
where the function $o(\cdot)$  depends only on $\gb$ and $v$ and vanishes 
as $\gd$ goes to 0.
\end{pro}

The Propositions above ensure that given a precision $\gd$, there is
a finite scale $k_0(\gd)$ after which the phases are uniformly segregated 
with this precision.

\subsection{Scheme of the proof}

The scheme of the proof is well known in the soft context of large 
deviations: one first proves an exponential tightness property and then
a weak large deviation principle (Proposition \ref{prop 2} holds
also for any bounded variation function with finite perimeter).
To be sure, the proof itself has nothing to do with the theory of large 
deviations: the
central tools here are the renormalization estimates leading to 
Peierls type bounds and estimate 
in the phase of small contours, and, of course,  the identification 
methods to produce  the macroscopic  surface 
tension in the precise logarithmic asymptotics.

Thus, Proposition \ref{prop 1} tells us that, under the appropriate 
renormalization, the occurrence of many small 
contours or of very large contours is unlikely.
It is a straightforward consequence of the general exponential
tightness Theorem \ref{thm Compactness}, which we state in 
Section~\ref{Coarse graining and mesoscopic phase labels}.
The statement is reminiscent to the results proven in \cite{BBP},
but the proof itself is  based on the analysis of the phase of small contours
developed in \cite{I2}, \cite{SS}, \cite{PfisterVelenik97}.

To prove Propositions \ref{prop 2} and \ref{prop 3}, we first consider
the macroscopic event $\big\{  \frac{\cM_k}{m^*} \in \cV(v,\gd) \big\}$
and by using several localization procedures, we reduce to compute the
probability of microscopic events from which, adopting the procedure
 developed in \cite{Cerf}, we can derive the exact 
surface tension factor.
This enables us to avoid the computations related to the
 microscopic phase boundaries at, however, a principal cost of loosing track
of the latter.

Since the most likely configurations in 
$\big\{  \frac{\cM_k}{m^*} \in \cV(v,\gd) \big\}$ are those for which
both phases coexist along the boundary of $\partial^* v$, we would like 
to prove that a microscopic interface is localized close to the 
boundary.
To derive the lower bound (Proposition \ref{prop 2}), one can enforce  
such a microscopic interface and then recover the surface tension factor.

This is not the case for the upper bound (Proposition \ref{prop 3})
because the $\bbL_1$ constraint $\big\{ \frac{\cM_k}{m^*} \in \cV(v,\gd) \big\}$
imposed on the magnetization is not strong enough to localize
the interface close to $\partial^* v$ : there might be 
mesoscopic fingers of one phase percolating into the other. 
To circumvent this problem, we follow an argument developed in \cite{BBBP}
and first prove a weak localization on a mesoscopic 
level. 
This involves a surgery  procedure called the minimal section argument.
This procedure ensures that one can chop off the mesoscopic fingers 
without changing too much the probability of the event and therefore
localize the interface on a mesoscopic level.
The renormalization is an essential feature of this proof.
Once the interface is localized on the mesoscopic level,
it remains to identify surface tension.

\medskip

We now proceed by first defining a coarse graining and deducing
the exponential tightness from Theorem \ref{thm Compactness}.
Then we compute the logarithmic asymptotics.

\section{Coarse graining and mesoscopic phase labels}
\label{Coarse graining and mesoscopic phase labels}
\setcounter{equation}{0}

At every mesoscopic scale $M=2^k$ the local magnetization $\cM_k$ gives a 
coarse grained representation of the system. Statistical properties of the 
microscopic configurations are waved out, and instead one keeps track only 
of the local order parameters over the corresponding mesoscopic blocks. These
are quantified by three values $\pm 1$ and $0$ according to whether they are 
sufficiently close to one of the two equilibrium values $\pm m^*$ or 
not. $0$-blocks play the role of the mesoscopic phase boundaries, and the 
$\pm 1$ blocks of the corresponding mesoscopic phase regions. Thus, the outcome
of the renormalization could be schematically represented as the 
following two-step diagram :
\begin{equation*}
\left\{
\begin{split}
&\text{Microscopic}\\
&\text{configurations}
\end{split}
\right\}
\ \longrightarrow\ 
\left\{
\begin{split}
&\text{Local}\\
&\text{magnetization}
\end{split}
\right\}
\ \longrightarrow\ 
\left\{
\begin{split}
&\text{Mesoscopic}\\
&\text{phase labels}
\end{split}
\right\} .
\end{equation*}

\noindent
There are two principal results to be discussed in this Subsection: we show 
that the $\bbL_1$-difference between the local magnetization and the 
corresponding phase labels vanishes on the exponential scale, and we give a 
general exponential tightness criterion for families of 
$\{\pm 1, 0\}$-valued phase label functions. 
In Section \ref{Examples of mesoscopic phase labels}, we will indicate how to 
construct phase labels in the case of Kac, percolation and nearest neighbor 
Ising models.\\

\noindent
{\bf Definition :}  A $\{\pm 1,0\}$-valued function $u$ on $\uTor$ is called
 a mesoscopic phase label, if there exists $k\in\bbN$, such that
 $u$ is an $\cF_k$-measurable function.

\medskip

\subsection{Tightness theorem for mesoscopic phase labels}

We fix now a sequence of non-negative numbers $\{\gr_k\}$ such that
\begin{equation}
\label{rhok}
\lim_{k\to\infty} \gr_k ~=~0 .
\end{equation}
The following compactness result holds uniformly in the microscopic
scales $N=2^n$.
\begin{thm}[Tightness of Mesoscopic Phase Labels]
\label{thm Compactness}
Let $N=2^n$ and assume that $\{u_k (\go ,x )\}$ is a  sequence 
of random mesoscopic phase label functions defined on the common probability
space $(\gO_N ,\cA_N ,\bbP_N )$, such that the realizations of  
$u_k\in\cF_{n-k},\ k=1,...,n$, 
and for every $k$ the following two conditions hold:

\noindent
{\bf A.} The distribution of the 
family of random variables $\{| u_k (\go ,x)|\}_{x\in\sTor{n-k}}$
is stochastically dominated by the Bernoulli site percolation  measure 
$\bbP_{\text{\rm perc}}^{\gr_k}$ on $\sTor{n-k}$. In particular,
\begin{equation}
\label{A}
\bbP_N\left( u_k (x_1 )=0,...,u_k (x_\ell )=0\right)~\leq 
(\rho_k)^\ell.
\end{equation}
\noindent
{\bf B.} If for two different points $x,y\in \sTor{n-k}$ the corresponding
$u_k$-phase labels have opposite signs, that is if $u_k (x)u_k (y)=-1$, then
on any finer scale $k^{\prime} \leq k$ any $*$-connected chain of 
$\sBox{n-k^{\prime}}$ blocks joining $\sBox{n-k}(x)$ to $\sBox{n-k}(y)$ 
contains at least one block with zero $k^{\prime}$-label.\\

\noindent
Then for every $a>0$ and $\gd >0$ there exists a finite scale 
$k_0 =k_0 (\gd )$,
such that
\begin{equation}
\label{tightbound}
\frac{1}{N^{d-1}}\log\bbP_N\left( u_k\in\cV ( K_a ,2\gd)^{\text{c}}\right)
\leq ~ -c_1 (d)\min\left\{ \gd 2^{n-dk}~ ,~\frac{a}{2^{(d-1)k_{0}}}~,
\frac{\gd  2^{n-dk_0}}{n^d}
\right\} ~,
\end{equation}
for all $k \geq k_0$ .
\end{thm}

\noindent
{\bf Remark~.} The proof of this general theorem is given
in Appendix~A. Notice that for $N$ sufficiently large we obtain a simpler 
surface order estimate which, for
every $\nu <1/d$ fixed, holds
uniformly in all mesoscopic scales $k_0 (\gd) \leq k\ \leq \nu\log N$,
\begin{equation}
\label{surface}
\frac{1}{N^{d-1}}\log\bbP_N
\left( u_k\in\cV ( K_a ,2\gd)^{\text{c}}\right)~\leq~
-c_1 (d)\frac{a}{2^{(d-1)k_{0}}} .
\end{equation}
Also an inspection of the proof shows that the tightness of the phase
labels on a certain scale $k$ does not depend on the validity of 
Assumptions~A and B on the successive scales $k^{\prime}>k$. In particular,
the estimate \eqref{surface} is valid on fixed (large) finite scales 
$k=k_0$, once the Assumption~A is satisfied, and once any $*$-connected 
sign changing chain 
of $k_0$-blocks necessarily contains a $0$-block. This simplified
version of Theorem~\ref{thm Compactness} is used in the case of Kac
potentials which we discuss in Subsection~\ref{subsection Kac potentials}.\qed \\

\subsection{Relation to magnetization profiles}

The original Gibbs measure is related to the above abstract setting in the 
following way: For every $N=2^n$, one constructs a (possibly enlarged) 
probability space $(\gO_N ,\cA_N ,\bbP_N )$, on which both the spin variables
$\gs\in \{ -1,+1  \}^{\Tor{N}}$ and various indexed families 
$\{ u_k^{\gz}\}$ of mesoscopic 
phase labels are defined. Such construction should enjoy the following 
set of properties:
\medskip

\noindent
{\bf C1.}  The marginal distribution of $\gs$ under $\bbP_N$ is precisely
$\Is_N$.

\noindent
{\bf C2.} For every $\gz >0$ the family $\{ u_k^{\gz}\}$ of mesoscopic
 phase labels satisfies Assumption~A of Theorem~\ref{thm Compactness} with
the corresponding sequence $\{ \gr_{k,\gz}\}$ of site percolation 
probabilities obeying \eqref{rhok}.

\noindent
{\bf C3.} For every $k\in \{0,...,n\}$ and $\gz > 0$ the local magnetization 
profile $\cM_k$ and the phase label $u_{k}^{\gz}$  are related as follows:
$\bbP_N$-a.s.,
\begin{equation}
\label{Mkuk}
\left| \cM_k (x) -m^{*}u_{k}^{\gz} (x)\right|~\leq~\gz\qquad\text{whenever}\ 
| u_{k}^{\gz}(x) |=1 .
\end{equation}

Notice that both functions above are $\cF_{n-k}$-measurable, that is 
\eqref{Mkuk} should be verified over the mesoscopic boxes indexed by 
the points $x\in\sTor{n-k}$.\\

Under conditions C1-C3, given any $\gd > 0$ one can choose the accuracy
$\gz$ of the coarse graining, a finite
scale $k_0 =k_0 (\gd,\gb )$ and a sequence of mesoscopic phase labels
 $\{ u_k^\gz \}$, such that for every $\nu <1/d$ fixed,
\begin{equation}
\label{Mkukbound}
\frac{1}{N^{d-1}}
\log\bbP_N \left( 
\max_{k_0\leq k\leq \nu n}\| \cM_k  -m^{*} u_{k}^\gz \|_1 >\gd
\right)~\leq ~
-c_2 \; 2^{(1 -d \nu) n} .
\end{equation}
Notice that \eqref{Mkukbound} holds uniformly in the size of the 
system $N=2^n$, once Assumptions C1-C3 do so.

Let us check \eqref{Mkukbound}. By the very construction,
$$
\| \cM_k  -m^{*}u_{k}^{\gz} \|_1~\leq ~ \gz +\frac{2}{|\sTor{n-k} |}
\sum_{x\in\sTor{n-k}} 1_{u_k^{\gz} (x)=0} .
$$
Consequently, using the domination by the Bernoulli site percolation 
(Assumption~A),
\begin{equation*}
\begin{split}
&\bbP_N\left( \| \cM_k  -m^{*}u_{k}^{\gz} \|_1 >\gd\right)~\leq ~
\bbP_N\left(\frac{1}{|\sTor{n-k} |}
\sum_{x\in\sTor{n-k}} 1_{u_k^{\gz} (x)=0} >\frac{\gd -\gz}{2}\right)\\
&\ \ \leq
\bbP_{\text{perc}}^{\rho_{k,\gz}}\left(\frac{1}{|\sTor{n-k} |}
\sum_{x\in\sTor{n-k}} 1_{u_k^{\gz} (x)=0} >\frac{\gd -\gz}{2}\right)~\leq~
\text{exp}\left\{ -c_1 2^{d(n-k)}\log\frac{\gd -\gz}{2\rho_{k,\gz}}\right\} .
\end{split}
\end{equation*}
The latter estimate is of the super-surface order once 
$ \rho_{k,\gz} \ll (\gd -\gz )/2$ and $k<n/d$.

\section{Examples of mesoscopic phase labels}
\label{Examples of mesoscopic phase labels}
\setcounter{equation}{0}

We show that mesoscopic phase labels can be constructed
in the case of Kac, percolation and Ising models.

\subsection{Kac potentials}
\label{subsection Kac potentials}

For this model mesoscopic phase labels are defined on the original space 
of spins $\gs\in\{-1, +1\}^{\Tor{N}}$ :
the coarse graining is obtained by averaging locally the magnetization.
Recall that we are using dyadic length scales $N=2^n$.

Phase labels are constructed in three steps. First, 
for any integer $k$ and $\gz>0$, we introduce the block spin variables 
$\bar u_k^\gz$ which label the boxes $\sBox{n-k}$ according to the averaged
magnetization over the boxes of the linear size $M=2^k$. These $\bar u_k^\gz$
 are constant on each of the blocks $\sBox{n-k}(x)$ with $x \in \sTor{n-k}$
\begin{eqnarray*}
\bar u_k^\gz(\gs, x)~
\df ~ \left\{
\begin{array}{l}
\pm 1 \qquad {\rm if} \ \quad
| \frac{1}{M^d} \sum_{ i \in \dBox{M} (2^n x)} \gs_i \mp m^* | < \gz,\\
0  \qquad {\rm otherwise}.
\end{array}
\right.
\end{eqnarray*}

In the Kac case we do not use Theorem~\ref{thm Compactness} in its full
generality, the object of the coarse graining is to choose a finite scale
$k_0$, such that the family of mesoscopic phase labels is exponentially
tight in $\bbL_1$. 
Recall that the scaling parameter is chosen such that $\gep = 2^{-m}$
with $m$ large but fixed.
Eventually finite renormalization scales $k_0$ are going to satisfy 
$k_0 = m+ a_0$, where $a_0$ depends on $\gb$ and $\gz$, but not on $m$.
The sign of the $k_0$-label over a box $\sBox{n-k_0}(x)$ depends on a 
more refined information on the fluctuations of the magnetization inside the 
box : we choose another scale $\ell_0;\ \ell_0 =m-b_0$, where, as in the case
of $a_0$, the scale $b_0$ will eventually depend only on $\gb$ and $\gz$,
and define the family of modified block spins $\{\tilde u_{k_0}^\gz\}$
on the $k_0$-scale as 
\begin{eqnarray*}
\tilde u_{k_0}^\gz (\gs, x)
\df \left\{
\begin{array}{l}
\pm 1 \qquad {\rm if} \quad
\qquad  \bar u_{\ell_0}^\gz (\gs, y) = \pm 1,
\qquad \forall ~y\in \sTor{n - \ell_0} \cap \sBox{n-k_0}(x) \\
0 \qquad {\rm otherwise}.
\end{array}
\right.
\end{eqnarray*}
Finally, we define the mesoscopic phase label functions 
$\{u^\gz_{k_0}(\gs,x)\}$.
If $\tilde u^\gz_{k_0} (\gs,x) = 0$, we set $u^\gz_{k_0} (\gs,x) = 0$.
If $x,y\in\sTor{n-k_0}$ are $*$-neighbors, but the corresponding modified
blocks spins satisfy
  $\tilde u^\gz_{k_0}(\gs,x) \, \tilde u^\gz_{k_0}(\gs,y) < 0$ then 
$u^\gz_{k_0} (\gs,x) = u^\gz_{k_0} (\gs,y) = 0$.
Otherwise, we set $u^\gz_{k_0}(\gs,x) = \tilde u^\gz_{k_0} (\gs,x)$.\\

A consequence of the Peierls estimate proven in \cite{CP} and 
\cite{BZ} is that assumption A is satisfied, namely 
\begin{thm}
For any $\gb >1$, there exists 
$\gz_0=\gz_0 (\gb) >0$, such that the following holds:
For any $\gz <\gz_0$ one can choose $\gep_0 = \gep_0(\gz)$, $a_0 =a_0 (\gz )$ 
and $b_0 =b_0 (\gz )$,
such that uniformly in the interaction parameters $\gep =2^{-m}<\gep_0$,
\begin{eqnarray*}
\Is_{\gep,N} \left( u^\gz_{k_0} (x_1) =0, \dots , 
u^\gz_{k_0} (x_r) =0  \right) 
\leq \exp \left( - \frac{c_0}{\gep^d} r \right),
\end{eqnarray*}
where, for every fixed $\gep =2^{-m}<\gep_0$, the mesoscopic phase labels
$u^\gz_{k_0}$ are constructed on the scales $k_0 = m+ a_0(\gz )$ and 
$l_0 =m-b_0 (\gz )$.
\end{thm}

\noindent
{\bf Remark.}
A more refined statement implying exponential decay of correlations was 
proven in \cite{BMP}. 
Notice that conditions C1-C3 of the previous Section are
 satisfied by definition of the mesoscopic phase 
label functions.  Notice also that assumption B of Theorem~\ref{thm Compactness}
is automatically satisfied on the $k_0$-scale. Thus, the family 
$\{ u^\gz_{k_0}\}$ is exponentially tight in $\bbL_1$.\qed \\

A similar renormalization procedure was carried out by Lebowitz, Mazel 
and Presutti \cite{LMP} for a system of point particles in $\bbR^d$
interacting with Kac potentials.
In this case the study of phase transition in the continuum is much more involved.
Beyond a proof of the liquid-vapor phase transition, their results 
 provide an accurate description of the system in terms of
mesoscopic phase labels  which represent the liquid and the gaseous
phases.
Such a coarse graining should be helpful to obtain further results on 
phase coexistence in the continuum.\\

\subsection{Bernoulli bond percolation}

Bernoulli bond percolation exhibits features similar to the Ising model
as phase transition and surface order behavior in a regime of phases
coexistence.
Nevertheless, as the setting is different from the Ising model,
we briefly recall some notation.
The set of edges is $\bbE = \big\{ \{x,y\} \; | \; x \sim y \big\}$,
where $x \sim y$ means that the vertices are nearest neighbors.
An edge $b$ in $\bbE$ is open if $\go_b =1$ and closed otherwise.
To any subset $\gL\Subset\bbZ^d$, we associate $[\gL]_e$ the set 
of edges in $\gL$.
The space of bonds configurations in $\gL$ is $\gO_\gL = \{ 0, 1\}^{[\gL]_e}$.
For a given $p$ in $[0,1]$, we define the Bernoulli bond percolation measure 
on $\gO_\gL$ by
\begin{eqnarray*}
\Perc^{p}_{\gL} (\go) =    
\prod_{b \in [\gL]_e} ( 1 -  p)^{1 - \go_b} p^{\go_b}  \, .
\end{eqnarray*}
For simplicity $\Perc^{p}_N$ denotes the measure on 
$\gO_N = \gO_{\Tor{N}}$.

Let $\go$ be a configuration in $\gO$,
an open path $(x_1, \dots ,x_n)$ is a finite sequence of distinct
nearest neighbors $x_1, \dots ,x_n$ such that on each edge
$\go_{\{x_i , x_{i+1} \}} = 1$.
We write $\{ A \lra B \}$ for the event such  that there exists an open
path joining a site of $A$ to one of $B$.
The connected components of the set of open edges of $\go$ are called
$\go$-clusters.

A phase transition is characterized by the occurrence of an infinite cluster. 
Define $\Theta_p$ by
\begin{eqnarray}
\label{perc Theta}
\Theta_p = 
\lim_{N \to \infty} \Perc_N^{p}(\{ 0 \lra  \partial \Tor{N} \}) \, ,
\end{eqnarray}
then there is a critical value $p_c$ in $]0,1[$ such that for any
$p$ below $p_c$ there is no percolation and $\Theta_p=0$, instead 
for any $p$ above $p_c$ the occurrence of an infinite cluster starting
from 0 has positive probability $\Theta_p$.
In the thermodynamic limit, there exists only one limiting Gibbs measure
and almost surely a unique infinite cluster with local density $\Theta_p$.
In order to mimic the coexistence of 2 phases in the finite domains
$\Tor{N}$, we say that one phase is formed by the largest cluster
and the other phase by the other clusters.

For this model, Pisztora introduced a renormalization procedure 
\cite{Pisztora1}, \cite{DP}, \cite{Pisztora2} which holds as soon as $p>p_c$
and $d \geq 3$.
The mesoscopic phase labels $\{ u^\gz_k \}$ will be defined for any 
mesoscopic scale $M = 2^k$, where $k$ is an integer which eventually depends 
on $N$.
This construction requires 2 steps.
The first step is to retain only the main features of the typical 
configurations on finite size boxes $\dBox{M}$.
Then we attribute a sign to the blocks $\sBox{n-k}$ according to the
phase they represent.
Set $M ' =2M$.
For any $x$ in $\sTor{n-k}$, the following events depend only on 
configurations in the box $\dBox{M'}(2^n x)$.
\begin{eqnarray*}
U_x = \left\{ \go \in \gO_{N} \;  \big| \;  
\text{there is a unique crossing 
cluster $C^*$ in  $\dBox{M'}(2^n x)$} \right\}.
\end{eqnarray*}
A crossing cluster is a cluster which intersects all the faces of the box.
Let $\ell$ be an integer smaller than $k$ which will be fixed later
\begin{eqnarray*}
R_x & = & U_x \bigcap \left\{ \go \in \gO_{N}  
\;  \big| \; 
\text{every open path in $\dBox{M'}(2^n x)$  with diameter 
 larger than $2^\ell$ } \right.\\
& & 
\text{is contained in $C^*$ }  \Big\},
\end{eqnarray*}
where the diameter of a subset $A$ of $\bbZ^d$ is 
$\sup_{x,y \in A} \|x - y\|_1$.
Finally, we consider an event which imposes that the density of the crossing
cluster in $\dBox{M}(2^n x)$ is close to $\Theta_p$ with accuracy
$\gz >0$
\begin{eqnarray*}
V_x^\gz = U_x \bigcap \big\{ \go \in \gO_{N}   \;  \big| \quad
| C^* \cap \dBox{M}(2^n x)| \in [\Theta_p -\gz, \Theta_p + \gz] \, 
M^d  \big\},
\end{eqnarray*}
where $| \cdot |$ denotes the number of vertices in a set.

Each box $\sBox{n-k}(x)$ is labeled by the variable 
${\tilde u}^\gz_{k} (\go,x)$ 
\begin{eqnarray*}
\forall x \in \sTor{n-k}, \qquad
{\tilde u}^\gz_{k} (\go,x)  
\df \left\{
\begin{array}{l}
1 \qquad     \text{if} \qquad  \go \in R_x \cap V_x^\gz,\\
 0 \qquad     \text{otherwise}.
\end{array}
\right.
\end{eqnarray*}
Let $\{x_1,\dots,x_r \}$ be vertices in $\sTor{n-k}$ not $*$-neighbors of 
$x$, then \cite{Pisztora1} implies that for every $p>p_c$, there exists $k_0(p, \gz)$, 
and $\ell_0 (p)$ such that for all $k \geq k_0$ and 
$k \geq \ell \geq \ell_0$
\begin{eqnarray*} 
\Perc^{p}_{N} \left({\tilde u}^\gz_{k} (x) = 0  \ 
\big| \ {\tilde u}^\gz_{k} (x_1 ), \dots,{\tilde u}^\gz_{k} (x_r) \right)
\leq  \exp( - c_1 \, 2^\ell) + \exp (- c_2 (\gz)  2^k),
\end{eqnarray*}
 From \cite{LSS} (Theorem 1.3), we deduce that for $k$ and $\ell$ 
large enough, 
the random variables $\{ {\tilde u}^\gz_{k} (x)\}$ are dominated by a 
Bernoulli site percolation measure $\bbP^{\rho_k}_{\rm perc}$
\begin{eqnarray}
\label{domination Y}
\rho_k \leq   \exp( - c(\gz) \, 2^{\ell}).
\end{eqnarray}
A straightforward way to recover the previous statement is to partition
$\sTor{n-k}$ into $c(d)$ sub-lattices $\big( \sTor{n-k-1, i} \big)_{i \le c(d)}$
which are translates of $\sTor{n-k-1}$. 
Any collection of vertices $\{x_1,\dots,x_r \}$ in $\sTor{n-k}$ can be
rearrange into $c(d)$ subsets $\{x_1^{(i)},\dots,x_{r_i}^{(i)} \}$
such that each $\{x_1^{(i)},\dots,x_{r_i}^{(i)} \}$ belongs to $\sTor{n-k-1, i}$.
Applying H\"older inequality, we get
\begin{eqnarray*} 
\Perc^{p}_N \left(
{\tilde u}^\gz_{k} (x_1 ) = 0, \dots,{\tilde u}^\gz_{k} (x_r) = 0 \right)
\leq 
\prod_{i=1}^{c(d)}
\Perc^{p}_N 
\left( {\tilde u}^\gz_{k} (x^{(i)}_1 ) = 0, \dots,
{\tilde u}^\gz_{k} (x^{(i)}_{r_i}) = 0 \right)^{\frac{1}{c(d)}}  \, .
\end{eqnarray*}
As the vertices in $\sTor{n-k-1, i}$ are not $*$-neighbors in $\sTor{n-k}$,
the domination by a Bernoulli product measure follows.

We say that a block $\sBox{n-k}(x)$ is regular if $\tilde u^\gz_k(x) =1$.
Finally we define the mesoscopic phase labels $u^\gz_k$ to be equal to 1 
on the regular blocks connected to the largest cluster
and to $-1$ on the regular blocks disjoint from the largest cluster.
Otherwise, we set 
${u}^\gz_{k} (\go,x) = {\tilde u}^\gz_{k} (\go,x) = 0$.
 From (\ref{domination Y}), the mesoscopic phase labels satisfy
assumption A.
Notice that if $x$ and $y$ are $*$-neighbors in $\sTor{n-k}$
the boxes $\dBox{M'}(2^n x)$ and $\dBox{M'}(2^n y)$ overlap. 
Choosing the parameter $\ell\leq k-3$ we insure that if the boxes $\sBox{n-k}(x)$ and 
$\sBox{n-k}(y)$ are both regular, then 
the crossing clusters in these boxes are connected.
This implies that assumption B is satisfied : two blocks with $k$-labels
of different signs cannot be $*$-connected.

The Bernoulli bond percolation model is precisely described by 
Pisztora's coarse graining, namely on a sufficiently large scale $2^k$,
the typical configurations have a unique crossing 
cluster surrounded by small islands of size smaller than $2^\ell$.
According to Theorem \ref{thm Compactness}, the family $\{ u^\gz_k \}$
is exponentially tight in $\bbL^1$.\\

\subsection{Ising nearest neighbor.}
\label{subsection : Ising nearest neighbor}

An extension of the preceding renormalization procedure applicable
to the Ising model has been also introduced in \cite{Pisztora1}.
Unlike Ising model with Kac potentials, this coarse graining is 
defined on an enlarged phase space via the FK representation.
For a  review of FK measures, we refer the reader 
to \cite{Pisztora1}, \cite{ACCN} and \cite{Grimmett}.

Let us recall the definition of the random cluster measures (or
FK measures) which are a generalization of the Bernoulli bond
percolation measures with correlated bond distribution.
To any subset $\gL$ of $\bbZ^d$  and $\pi$ included in $\partial \gL$, we 
associate a set of edges 
\begin{eqnarray*}
[\gL]_e^\pi = \big\{ \{x,y\} \; | \; x \sim y, \ x \in \gL, \ y \in \gL \cup \pi 
 \big\},
\end{eqnarray*}
and the space of configurations in $\gL$ is $\gO_\gL^\pi = \{ 0, 1\}^{[\gL]^\pi_e}$. 
The first step is to introduce a measure on $\gO_\gL^\pi$.
A vertex $x$ of $\gL$ is called $\pi$-wired if it
is connected by an open path to $\pi$.
We call $\pi$-clusters the clusters defined with respect to the boundary 
condition $\pi$ : a $\pi$-cluster is a connected set of open edges in 
$\gO_\gL^\pi$ and we identify to be the same cluster all the clusters
which are $\pi$-wired, i.e. connected to $\pi$.
For a given $p$ in $[0,1]$, we define the FK measure on $\gO_\gL^\pi$ with
boundary conditions $\pi$ by
\begin{eqnarray*}
\Perc^{\pi,p}_{\gL} (\go) = {1 \over {\bf Z}_{\gL}^{\pi,p}}   
\left( \prod_{b \in [\gL]_e^\pi} ( 1 -  p)^{1 - \go_b} p^{\go_b} \right)
2^{c^\pi (\go)},
\end{eqnarray*}
where $Z_{\gL}^{\pi,p}$ is a normalization factor and $c^\pi (\go)$ 
is the number of clusters which are not $\pi$-wired.
If $\pi = \partial \gL$ then the boundary conditions are said to be
wired and the corresponding FK measure on $\gO^{\rm w}_{\gL}$
is denoted by $\Perc^{\rm w,p}_{\gL}$.
Finally, the periodic measure on the torus $\Tor{N}$ is denoted
by $\Perc^{\rm per,p}_N$ and the phase space by $\gO_N^{\rm per}$.

In order to recover the Gibbs measure $\Is_{\gL}$, we fix the
percolation parameter $p_\gb = 1 - \exp(-2 \gb)$ and
generate the edges configuration $\go$ in $\gO_N^{\rm per}$ according to
the measure $\Perc^{\rm per ,p_\gb}_N$.
Given $\go$, we equip
randomly each $\go$-cluster with a color $\pm 1$ with probability 
${1 \over 2}$ independently from the others.
This amounts to introducing the measure $P_N^\go$ on $\{-1,1\}^{\Tor{N}}$
such that the spin $\gs_i$ has the color of the cluster attached to $i$.
The Gibbs measure $\Is_{N}$ can be viewed as the first marginal
of the coupled measure 
$\Joint_N (\gs,\go) =  P_N^\go (\gs) \Perc_N^{\rm per,p_\gb}(\go)$
on the space $\{-1,1\}^{\Tor{N}} \otimes \gO_N^{\rm per}$.
In the case of $\pi$-wired boundary conditions, the spins attached
to the $\pi$-wired cluster are equal to 1.\\

As a consequence of this representation, one has for any increasing
sequence of sets $\gL_N$
\begin{eqnarray*}
m^* = \lim_{N \to \infty} \Is^+_{\gL_N}(\gs_0)  =
\lim_{N \to \infty} \Perc_{\gL_N}^{\rm w,p_\gb}(\{ 0 \lra  \partial \gL_N \} ) 
= \Theta_{p_\gb} .
\end{eqnarray*}
In the following, we use $m^*$ or $\Theta_{p_\gb}$ depending on the context.
Furthermore, we suppose that
\begin{eqnarray}
\label{Theta}
\lim_{N \to \infty} 
\Perc_{\gL_N}^{\rm f,p_\gb}(\{ 0 \lra  \partial \gL_N \} )
= \lim_{N \to \infty} 
\Perc_{\gL_N}^{\rm w,p_\gb}(\{ 0 \lra  \partial  \gL_N \} ) 
= \Theta_{p_\gb}.
\end{eqnarray}
This property is satisfied for all $\gb$ outside a subset of $\bbR$  which 
is at most countable (see Lebowitz \cite{Lebowitz} and Pfister \cite{Pf1}).

On the scale $M= 2^k$, we define, in the same way as for Bernoulli
bond percolation, the variables $\tilde u_k^\gz (\go,x)$ which are
piecewise constant on each box $\sBox{n-k} (x)$ with $x$ in
$\sTor{n-k}$.
The mesoscopic phase labels depend on the averaged magnetization in 
regular blocks.
Define the label of $\sBox{n-k} (x)$ by 
\begin{eqnarray*}
u^\gz_{k} (\gs,\go,x)   \df
\left\{ \begin{array}{l}
{\rm{sign}}(C^*)  \qquad     {\text{if}}
 \qquad  {\tilde u}^\gz_{k} (\go,x)  =1
\ \ {\rm{and}} \ \  | \cM_k (\gs, x) - {\rm{sign}}(C^*) \,  m^* | 
< 2 \gz,\\
0 \qquad  \qquad \quad    \text{otherwise},
\end{array}
\right.
\end{eqnarray*}
where $C^*$ is the crossing cluster in $\dBox{M} (2^n x)$.\\

In a regular box $\sBox{n-k} (x)$ (i.e. $\tilde u^\gz_k (x) = 1$),
the averaged magnetization is controlled by the random coloring of the small 
clusters included in $\dBox{M} (2^n x)$.
So that the averaged magnetization in a regular box is independent
of the configurations in the neighboring boxes.
In the case of Ising model, the additional parameter $\ell=\ell(k)$ 
is tuned in order to control the fluctuations of the magnetization 
over the small clusters.
As a consequence of this, assumptions A, B and C1-C3 are satisfied for $p_{\gb}$
above a certain non-trivial slab percolation threshold 
$p_{\tilde \gb_c }$, which is conjectured to coincide with 
$p_{\gb_c}$
(see \cite{Pisztora1} for details), and Theorem \ref{thm Compactness}
holds.

\noindent
{\bf Remark~.}
Using the notations of this Subsection, the set $\frB_p$ introduced in
Subsection \ref{subsection Main results} could be defined as
\begin{eqnarray*}
\frB_p = \{ \gb \; : \; \gb > \tilde \gb_c \ \text{and \eqref{Theta} holds}   \}.
\end{eqnarray*}


\section{Surface tension}
\label{section_Surface tension}
\setcounter{equation}{0}

We are going to derive Propositions \ref{prop 2} and \ref{prop 3}
for Ising model with nearest neighbor interaction.
As explained before, the philosophy of the proof is to start
from the macroscopic level and to localize successively on
finer scales with the help of a coarse graining.
The approach itself is quite general. 
Nevertheless the coarse graining is model dependent, therefore we will 
need first to state an alternative representation of the surface tension in 
terms of the FK representation in order to use the estimates which will be
obtained from Pisztora's coarse graining.
The idea of such definitions has been introduced in \cite{Cerf}.\\

\subsection{FK representation}

We fix $\vec{n}$ a vector in $\bbS^{d-1}$ and study $\tau_\gb (\vec{n})$.
Following notation of Subsection~\ref{subsection Surface tension},
we consider, for any $\gep$ positive, the parallelepiped
$\widehat \gL(N,\gep N)$  of $\bbR^d$ oriented according to $\vec{n}$.
Namely, the basis of $\widehat \gL(N,\gep N)$ with side lengths 
equal to $N$ is orthogonal to $\vec{n}$ and 
the other sides have lengths equal to $\gep N$.
For simplicity its microscopic counterpart
$\widehat \gL(N,\gep N) \cap \bbZ^d$ will be denoted by $\gL_N (\gep)$.\\

By using the correspondence between the Ising model and the FK representation,
one can rewrite $\tau_\gb$ in terms of the bond model.
Let $\{ \partial^+ \gL_N (\gep) \nlra \partial^- \gL_N (\gep)  \}$ be 
the event such that there is no open path inside $\gL_N (\gep)$  
joining $\partial^+ \gL_N (\gep) $ to $\partial^- \gL_N (\gep) $. 
Then,
\begin{eqnarray}
\label{ST 1}
\tau_\gb (\vec{n}) = \lim_{N \to \infty} \,  - {1 \over N^{d-1}} \log 
\Perc^{\rm{w}, p_\gb}_{\gL_N(\gep)} 
\big( \{ \partial^+ \gL_N (\gep) \nlra \partial^- \gL_N (\gep) \} \big).
\end{eqnarray}
Notice that the event $\{ \partial^+ \gL_N  (\gep) \nlra \partial^- \gL_N 
(\gep) \}$ takes only into account the paths inside $\gL_N (\gep) $ and not the
identification produced by wired boundary conditions. The relation above will
be useful only in the proof of Proposition  \ref{prop 2}.\\

We are now going to state an approximate expression of the surface
tension which is weakly dependent on the boundary conditions.
It will be used in the derivation of Proposition \ref{prop 3}.
Let $\gL_N ' (\gep)$ be the the parallelepiped
\begin{eqnarray}
\label{ST slab}
\gL_N ' (\gep)  = \left\{ i \in \gL_N(\gep) \quad \big| \quad 
\vec{i} \cdot  \vec{n} \in [- \frac{\gep}{4} N,\frac{\gep}{4} N] \right\},
\end{eqnarray}
and denote by $\partial^\top \gL_N ' (\gep) $ 
(resp $\partial^\bot \gL_N ' (\gep) $) the face of
$\partial^+ \gL_N ' (\gep) $ (resp $\partial^- \gL_N ' (\gep) $) orthogonal 
to $\vec{n}$. 
Let $\{ \partial^\top \gL_N ' (\gep)  \nlra \partial^\bot  \gL_N '  (\gep) \}$ 
be the event such that there is no open path inside $\gL_N ' (\gep)$  connecting 
$\partial^\top \gL_N ' (\gep) $ to $\partial^\bot \gL_N ' (\gep) $. 
One has
\begin{lem}{\rm \bf [\cite{Bo} $\gb \gg 1$, \cite{CePi} $\gb \in \frB_p$]}
\label{lem surface tension}
For any $\gb \in \frB_p$
\begin{eqnarray}
\label{ST 2}
\tau_\gb (\vec{n}) = - {1 \over N^{d-1}}
 \log \Perc^{\pi,\rm p_\gb}_{\gL_N (\gep)  } 
\left( 
\{ \partial^\top \gL_N ' (\gep) \nlra \partial^\bot \gL_N ' (\gep)\} 
\right)  + c_{\gep,N} (\pi),
\end{eqnarray}
where the function $c_{\gep,N}$ goes to 0 as $N$ tends to infinity and 
$\gep$ goes to 0,  uniformly over the boundary conditions $\pi$ and 
$\vec{n} \in \bbS^{d-1}$.
\end{lem}

As it will be explained in Part~\ref{part_boundary}  on the wetting phenomenon,
the system is  in fact extremely sensitive to boundary conditions. Nevertheless
in the above Lemma, the interface is constrained to be in  $\gL_N '(\gep)$, so
that it does not feel the influence of the boundary : the boundary conditions
are screened because the system  relaxes to equilibrium in the region $\gL_N
(\gep) \setminus \gL_N '(\gep)$.

Let us first examine the influence of the boundary conditions $\pi$ on 
the faces of $\gL_N (\gep)$ orthogonal to $\vec{n}$.
As $\{ \partial^\top \gL_N ' (\gep) \nlra \partial^\bot \gL_N ' (\gep)\}$
is a decreasing event, FKG inequality imply that it is enough to
check that
\begin{eqnarray}
\label{ST 3}
\tau_\gb(\vec{n}) = \lim_{N \to \infty} \,  - {1 \over N^{d-1}} \log 
\Perc^{\rm f, w, p_\gb}_{\gL_N (\gep)} 
\big( \{ \partial^+ \gL_N ' (\gep)  \nlra \partial^- \gL_N ' (\gep)  \} \big),
\end{eqnarray}
where $\Perc^{\rm f, w, p_\gb}_{\gL_N (\gep) }$ is the FK measure
with free boundary conditions on the faces orthogonal to 
$\vec{n}$ and wired on the others.
This can be proved by means of a Peierls argument for $\gb$ large enough
\cite{Bo} or by an analysis of the relaxation of the clusters density
for $\gb$ in $\frB_p$ \cite{CePi}.

As already noticed in \cite{Cerf} in the context of percolation,
the influence of the boundary conditions on the sides of 
$\gL_N (\gep)$ parallel to $\vec{n}$ is negligible as $\gep$ goes to 0.
This explains that the factor $c_{\gep,N} ( \cdot)$ vanishes uniformly over
the boundary conditions.\\

\subsection{Extended representation}

We would like to stress that the previous treatment of the
surface tension is not satisfactory and a more coherent approach would 
be to consider a more general definition independent of the model in terms 
only of mesoscopic phase labels.
In fact, a definition of surface tension valid in an abstract setting 
would be difficult to use because  the surgical procedure of the minimal
section argument requires a precise knowledge of how
the microscopic system is related to the mesoscopic phase labels.

\section{Lower bound : Proposition \ref{prop 2}}
\setcounter{equation}{0}

The proof is divided into 3 steps. We first start by approximating
the surface $\partial^* \cK_m$ by a regular surface $\partial{\widehat K} $ 
and imposing the condition that a mesoscopic interface exists close to
$\partial {\widehat K}$.
Then, using the definition of surface tension (\ref{ST 1}), 
we derive Proposition \ref{prop 2}.\\

\subsection{Step 1 : Approximation procedure.}

\begin{figure}[t]
\centerline{
\psfig{file=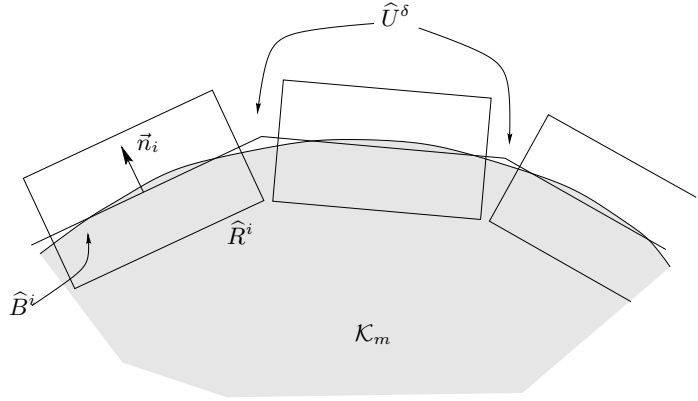,height=5cm}}
\caption{Polyhedral approximation.} 
\figtext{ 
\writefig        0.35   6.50    {\footnotesize $\widehat U^\gd$} 
\writefig       -2.90   4.85    {\footnotesize $\vec n_i$}
\writefig       -4.60   2.65    {\footnotesize $\widehat B^i$}
\writefig       -1.70   3.60    {\footnotesize $\widehat R^i$} 
\writefig        0.00   2.30    {\footnotesize $\cK_m$} 
}
\label{fig_approx_lowerbd}
\end{figure}

A polyhedral set has a  boundary included in the union of 
a finite number of hyper-planes. 
The surface $\partial^* \cK_m$ can be approximated as follows
(see Fig. \ref{fig_approx_lowerbd})
\begin{thm}
\label{thm approx}
For any $\gd$ positive, there exists a polyhedral set ${\widehat K}$ such that 
$$
\1_{\widehat K} \in \cV( \1_{\cK_m}, \gd)
\qquad {\rm  and} \qquad
\big| \cW_\gb ({\widehat K})-\cW_\gb (\cK_m) \big| \leq \gd.
$$
For any $h$ small enough there are $\ell$ disjoint parallelepipeds 
$\widehat R^1, \dots, \widehat R^{\ell}$   with basis 
$\widehat B^1, \dots, \widehat B^{\ell}$ included 
in $\partial {\widehat K}$ of side length $h$ and height $\gd h$.
Furthermore, the sets $\widehat B^1, \dots, \widehat B^{\ell}$ cover 
$\partial {\widehat K}$ up to a set 
of measure less than $\gd$ denoted by 
$\widehat U^\gd =\partial {\widehat K} \setminus \bigcup_{i =1}^\ell \widehat B^i$ 
 and they satisfy
\begin{eqnarray*}
\Big| \sum_{i = 1}^{\ell} \int_{\widehat B^i} \tau_\gb (\vec{n}_i) \, d \cH^{(d-1)}_x 
- \cW_\gb (\cK_m)  \Big| \leq \gd,
\end{eqnarray*}
where the normal to $\widehat B^i$ is denoted by $\vec{n}_i$.
\end{thm}
\noindent
The proof is a direct application of Reshtnyak's Theorem and can be found
in the paper of Alberti, Bellettini \cite{AlBe}.

Using Theorem \ref{thm approx}, we can reduce the proof of Proposition 
\ref{prop 2} to the computation of the probability of 
$\{ \frac{\cM_k}{m^*} \in \cV( \1_{\widehat K}, \gd) \}$.
According to (\ref{Mkukbound}) the estimates can be restated in terms
of the mesoscopic phase labels.
For any $\gd >0$, there exists $\gz = \gz(\gd)$ and $k_0 = k_0 (\gd)$ 
such that Proposition \ref{prop 2} will be implied by
\begin{eqnarray}
\label{lower 1}
\liminf_{N \to \infty} \;  { 1 \over N^{d-1}}
\min_{k_0 (\gd) \leq k \leq \nu n} \, \log  \Joint_N 
\left(  u^\gz_k \in \cV (\1_{\widehat K},\gd) \right) \geq
- \cW_\gb ({\widehat K}) - o(\gd).
\end{eqnarray}

\subsection{Step 2 : Localization of the interface.}

The images of ${\widehat K}$, $\widehat R^i$ and $\widehat U^\gd$
in $\Tor{N}$ will be denoted by $K_N$, $R^i_N$ and $U^\gd_N$.
In order to enforce a mesoscopic interface which crosses each $R^i_N$, 
we define the event 
$$
\cA = \inter_{i = 1}^{\ell} 
\{ \partial^+ R^i_N  \nlra \partial^- R^i_N \} \;  .
$$
We consider also $\cB$ the set of configurations such that
the bonds at distance less than 10 of $U^\gd_N$ are closed.
Notice that these events depend only on bonds variables.
One has 
\begin{eqnarray}
\label{lower 2}
\Joint_N  \left(  u^\gz_k \in \cV (\1_{\widehat K}, \gd) \right) \geq
\Joint_N  \left(  \left\{ u^\gz_k \in \cV (\1_{\widehat K}, \gd ) \right\} 
\cap \cA \cap \cB \right).
\end{eqnarray}
The interface imposed by the event $\cA \cap \cB$ decouples 
$K_N$ from its complement, therefore
the system is in equilibrium in $K_N$ and $K_N^c$ : a proof similar 
to the one of Theorem \ref{thm Compactness} implies that one can
choose $\gz' = \gz' (\gd)$ and $k_0 ' = k_0 '(\gd)$ such that
\begin{eqnarray*}
\lim_{N \to \infty} \; \max_{k_0 ' (\gd) \leq k \leq \nu n} \, 
\Joint_N \left(
\int_\gL | u_k^{\gz '} (x) - 1| \, dx \geq \frac{\gd}{2} \  
{\rm or} \ 
\int_\gL | u_k^{\gz '} (x) + 1| \, dx \geq \frac{\gd}{2} \ \Big| \  
\cA \cap \cB \right) = 0 \; , 
\end{eqnarray*}
where $\gL$ stands for $\widehat K$ or $\widehat K^c$.
So that (\ref{lower 2}) can be rewritten for $N$ large 
enough as
\begin{eqnarray}
\label{lower 3}
\min_{k_0 (\gd) \leq k \leq \nu n} \,
\Joint_N  \left(  u^\gz_k \in \cV (\1_{\widehat K},\gd) \right) \geq
\frac{1}{8} \, 
\Perc^{\rm per, p_\gb}_N  \left( \cA \cap \cB \right).
\end{eqnarray}

\subsection{Step 3 : Surface tension.}

Combining  the definition of surface tension (\ref{ST 1}),
inequality  (\ref{lower 3}) and Theorem \ref{thm approx}, we get
\begin{eqnarray*}
\liminf_{N \to \infty} \, \frac{1}{N^{d-1}} \, 
\min_{k_0 (\gd) \leq k \leq \nu n} \, \log
\Joint_N  \left(  u^\gz_k \in \cV (\1_{\widehat K},\gd) \right) \geq
- \sum_{i =1}^\ell \int_{{\widehat B}^i} \, \tau_\gb (\vec{n}_i) \, 
d \cH_x^{d-1} - o(\gd).
\end{eqnarray*}
We have also used the fact that the event $\cB$ is supported by at most 
$c(d,\gd) N^{d-1}$ edges where $c(d,\gd)$ vanishes as $\gd$ goes to 0.
Therefore the probability of $\cB$ is negligible with respect to a surface 
order.

\section{Upper bound  : Proposition \ref{prop 3}}
\setcounter{equation}{0}

The proof is divided into 3 steps. First we decompose $\partial^* v$ 
in order to reduce the proof to local computations in small regions.
Then in each region we localize the interface on the mesoscopic level 
via the minimal section argument.
Finally the last step is devoted to the computation of the surface tension
factor.

\subsection{Step 1 : Approximation procedure.}

\begin{figure}[t]
\centerline{
\psfig{file=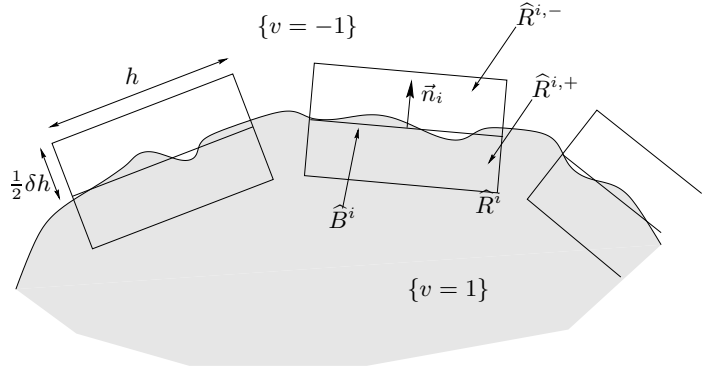,height=5cm}}
\figtext{ 
\writefig       -3.05   4.22    {\footnotesize $h$} 
\writefig       -4.6    2.80    {\footnotesize $\tfrac12\gd h$} 
\writefig       -1.3    4.9     {\footnotesize $\{v=-1\}$} 
\writefig       0.70    1.40    {\footnotesize $\{v=1\}$}
\writefig       0.88    4.00    {\footnotesize $\vec n_i$}
\writefig       -0.35   2.30    {\footnotesize ${\widehat B}^i$}
\writefig       2.35    4.05    {\footnotesize ${\widehat R}^{i,+}$} 
\writefig       2.14    5.00    {\footnotesize ${\widehat R}^{i,-}$} 
\writefig       1.60    2.50    {\footnotesize ${\widehat R}^i$} 
}
\caption{Approximation by parallelepipeds.} 
\label{fig_approx_upperbd}
\end{figure}

We approximate $\partial^* v$ with a finite number of parallelepipeds
(see Fig. \ref{fig_approx_upperbd}).
\begin{thm}
\label{theo ABCP}
For any $\gd$ positive, there exists $h$ positive such that there are 
$\ell$ disjoint parallelepipeds ${\widehat R^1}, \dots, \widehat R^{\ell}$ 
included in $\uTor$ 
with basis $\widehat B^1, \dots, \widehat B^\ell$ of 
size $h$ and  height $\gd h$.
The basis $\widehat B^i$ divides $\widehat R^i$ in 2 parallelepipeds 
$\widehat R^{i,+}$ and $\widehat R^{i,-}$
and we denote by $\vec{n}_i$ the normal to $\widehat B^i$.
Furthermore, the parallelepipeds satisfy the following properties
\begin{eqnarray*}
\int_{\widehat R^i} | \cX_{\widehat R^i}(x)  - v (x)| \, dx \leq \gd \, 
\vol(\widehat R^i) \quad
{\rm{and}} \quad
\Big| \sum_{i = 1}^{\ell} \int_{\widehat B^i} \tau_\gb (\vec{n}_i) \,  d \cH^{(d-1)}_x - 
\cW_\gb (v) \Big| \leq \gd,
\end{eqnarray*}
where $\cX_{\widehat R^i} = 1_{\widehat R^{i,+}} - 1_{\widehat R^{i,-}}$ 
and the volume of $\widehat R^i$ is ${\vol(\widehat R^i)} = \gd h^d$.
\end{thm}
\noindent
This Theorem is a rather standard assertion  of the  geometric measure Theory.
 A variation of it has been formulated and applied in the context of the $\bbL_1$-theory
of phase segregation in 
 \cite{ABCP} along with a sketch of the proof, which, however, contained a gap 
 (see \cite{Bo} for a detailed proof along the lines of \cite{ABCP}). 
A very clean alternative derivation of a similar result  has been given 
by Cerf \cite{Cerf} using the Vitali covering Theorem.

Theorem~\ref{theo ABCP} enables us to decompose the boundary into regular sets
(see Fig. \ref{fig_approx_upperbd})
so that it will be enough to consider events of the type
\begin{eqnarray*}
\left\{ \frac{\cM_k}{m^*} \in 
\bigcap_{i =1}^{\ell} \,  \cV(\widehat R^i , \gd \vol(\widehat R^i))
\right\} \, ,
\end{eqnarray*}
where $\cV( \widehat R^i , \gep)$ is the $\gep$-neighborhood of 
$\cX_{\widehat R^i}$ 
\begin{eqnarray*}
\cV (\widehat R^i, \gep) = \left\{ v^\prime \in \bbL^1 \big( \uTor  \big)
\ \big| \quad  \int_{\widehat R^i} | v^\prime(x) - \cX_{\widehat R^i}(x) |  
\, dx 
\leq \gep  \right\}.
\end{eqnarray*}
Using (\ref{Mkukbound}), we see that to derive Proposition \ref{prop 3},
it is equivalent to prove the following statement for any $\gd$
positive and $k_0 = k_0 (\gd)$, $\gz = \gz (\gd)$
\begin{eqnarray*}
\limsup_{N \to \infty}  \frac{1}{N^{d-1}}
\max_{k_0 (\gd) \leq k \leq \nu n} \, \log \Joint_N \big( u^\gz_k \in 
\bigcap_{i =1}^{\ell} \,  \cV(\widehat R^i , \gd \vol(\widehat R^i)) \big)
\leq - \cW_\gb (v) + C(\gb,v) \gd.
\end{eqnarray*}

\subsection{Step 2 : Minimal section argument.}

\begin{figure}[t]
\centerline{
\psfig{file=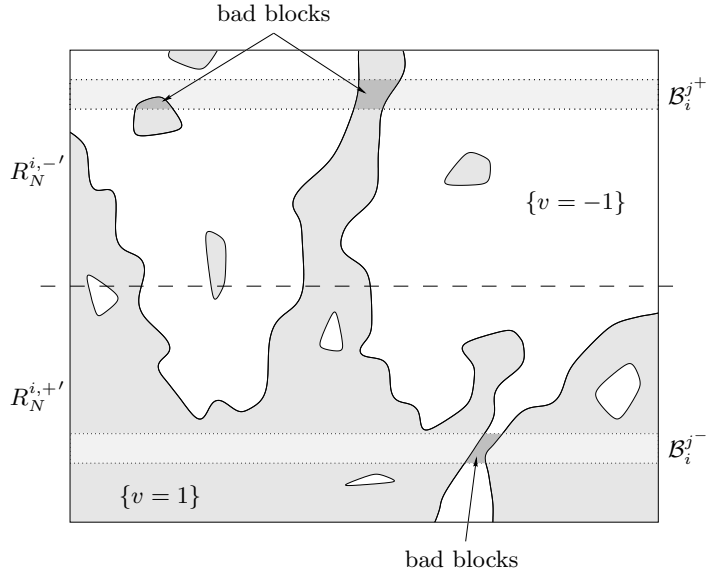,height=7cm}}
\figtext{ 
\writefig       0.60    .14     {\footnotesize bad blocks} 
\writefig       -1.9    7.40    {\footnotesize bad blocks} 
\writefig       2.20    4.9     {\footnotesize $\{v=-1\}$} 
\writefig       -3.20   1.00    {\footnotesize $\{v=1\}$}
\writefig       4.10    1.60    {\footnotesize $\cB_i^{j^-}$}
\writefig       4.10    6.30    {\footnotesize $\cB_i^{j^+}$}
\writefig       -4.65   2.30    {\footnotesize ${R_N^{i,+}}'$} 
\writefig       -4.65   5.30    {\footnotesize ${R_N^{i,-}}'$} 
}
\caption{Minimal sections.} 
\label{fig_minsec}
\end{figure}

The microscopic domain associated to $\widehat R^i$ is 
$R^i_N = N \widehat  R^i \cap \Tor{N}$.
We also set $R^{i,+}_N = N \widehat R^{i,+} \cap \Tor{N}$ and 
$R^{i,-}_N = R^i_N  \setminus R^{i,+}_N$.
At the scale $M =2^k$,
we associate to any configuration $(\gs,\go)$ the set of {\it bad} 
boxes which are the boxes $\dBox{M}$ intersecting $R_N^i$ labeled by $0$ 
and the ones 
intersecting $R^{i,+}_N$ (resp $R^{i,-}_N$) labeled by $-1$ (resp $1$).
For any integer $j$, we set $\widehat B^{i,j} = \widehat  B^i +
j \, c(d) 2^{n-k} \, \vec{n}_i$ and define
\begin{eqnarray*}
B^{i,j}_N = \big\{ j '  \in {R^i_N}   \ | \ 
\exists x \in \widehat B^{i,j}, \qquad  \|j' - Nx\|_1 \leq 10 \big\}.
\end{eqnarray*}
Let $\cB_i^j$ be the smallest connected set of boxes $\dBox{M}$ 
intersecting $B^{i,j}_N$.
By construction the $\cB_i^j$ are disjoint surfaces of boxes.
For $j$ positive, let $n_i^+(j)$ be the number of {\it bad} boxes in
$\cB_i^j$ and define
\begin{eqnarray*}
n^+_i = \min \big\{ n_i^+(j): \qquad 0< j < \frac{\gd h}{2c(d)} 2^{n-k} \big\}.
\end{eqnarray*}
Call $j^+$ the smallest location  where the minimum is achieved and define
the minimal section in $R^{i,+}_N$as $\cB_i^{j^+}$.
For $j$ negative, we denote by $\cB_i^{j^-}$ the minimal section in 
${R^{i,-}_N}$ and $n_i^-$ the number of {\it bad} boxes in $\cB_i^{j^-}$
(see Fig.~\ref{fig_minsec}).

For any configuration $(\gs,\go)$ such that 
$u^\gz_k (\gs,\go)$ belongs to 
$\bigcap_{i =1}^{\ell} \,  \cV(\widehat R^i , \gd \vol(\widehat R^i))$,
one can bound the number of {\it bad} boxes in the minimal sections by
\begin{eqnarray}
\label{upper 3}
\sum_{i =1}^\ell n^+_i + n_i^- \leq \gd C_1(v) 2^{(d-1)(n-k)} \, .
\end{eqnarray}
Such an estimate implies that a mesoscopic interface is mainly located
between the 2 minimal sections and that only some mesoscopic fingers
attached to the interface may percolate.
As these fingers will cross the minimal sections through {\it bad} boxes,
the strategy is therefore to modify the configuration $\go$ on the {\it bad} 
boxes so that no fingers can percolate in the new configuration.
More precisely, we introduce the set 
\begin{eqnarray*}
\cA = \big\{ \go \in \gO^{\rm per}_{N} \ \big| \quad
\exists \gs \ {\rm such \  that}  \  u^\gz_k (\gs,\go) \in 
\bigcap_{i =1}^{\ell} \,  \cV(\widehat R^i , \gd \vol(\widehat R^i)) 
 \big\} \, ,
\end{eqnarray*}
and for any $\go$ in  $\cA$ define $\bar \go$ the configuration with closed edges on 
the boundary of the {\it bad} blocks in the minimal sections
and equal to $\go$ otherwise.
Inequality (\ref{upper 3}) implies that $\go$ and $\bar \go$
differ only on at most $\gd C_2(v) N^{d-1}$ edges, so that we
can control precisely the cost of the surgical procedure which 
consists in isolating the {\it bad} blocks in the minimal sections
by closing the edges around them.
\begin{eqnarray}
\label{upper 4}
\Joint_N \left( u^\gz_k (\gs,\go) \in 
\bigcap_{i =1}^{\ell} \,  \cV(\widehat R^i , \gd \vol(\widehat R^i))  
\right) & \leq &
\Perc_N^{\rm per, p_\gb} \big( \cA \big)\\
& \leq &
\exp \big( \gd \, C_3 (v,\gb) N^{d-1} \big) \;
\Perc_N^{\rm per, p_\gb} \big( \bar \cA \big) \, , \nonumber
\end{eqnarray}
where $\bar \cA = \{ \bar \go \; |\; \go \in \cA\}$.

\subsection{Step 3 : Surface tension estimates.}

Let ${\widehat R^i} \, '$ be the parallelepiped included in $\widehat R^i$ 
with basis $\widehat B^i$ and height $\frac{\gd}{2} h$. 
Its microscopic counterpart is ${R^i_N} \, '$.
We are going to check now that $\bar \cA$ is included in 
$\bigcap_{i =1}^\ell \{ \partial^\top {R^i_N} ' \nlra \partial^\bot 
{R^i_N} '\}$.
This amounts to say that not only the minimal section argument 
enables us to find a mesoscopic interface in $R^i_N$ but that in fact
this interface exists on the microscopic level.
To see this, choose any configuration $\go$ in $\cA$ which contains an 
open path ${\bf C}$ joining 
$\partial^\top {R^i_N} '$ to $\partial^\bot {R^i_N} '$ and suppose
that  ${\bf C}$ crosses the minimal sections without intersecting a 
{\it bad} box.
Then ${\bf C}$ intersects 2 regular boxes $\dBox{M}(2^n x^+)$ and 
$\dBox{M}(2^n x^-)$ in $\cB_i^{j^+}$  and $\cB_i^{j^-}$.
According to the definition of the coarse graining,
this would imply that the crossing clusters of $\dBox{M}(2^n x^+)$ and 
$\dBox{M}(2^n x^-)$ are connected to ${\bf C}$, so that
${\tilde u}^\gz_k (x^+) = {\tilde u}^\gz_k (x^-)$.
Therefore one of these boxes has to be a {\it bad} box.

From (\ref{upper 4}), we get 
\begin{eqnarray*}
\label{upper 5}
\Joint_{N} \left( u^\gz_k \in 
\bigcap_{i =1}^{\ell} \,  \cV(\widehat R^i , \gd \vol(\widehat R^i))  
\right) && \leq \exp \big( \gd \, C_3 (v,\gb) N^{d-1} \big) \\
&& \qquad \Perc^{\rm per,p_\gb}_{N} \big(  
\bigcap_{i =1}^\ell \{ \partial^\top {R^i_N} ' \nlra 
\partial^\bot {R^i_N} ' \} \big).
\end{eqnarray*}
Conditioning outside each domain $R^i_N$ and using (\ref{ST 2}), 
we derive
\begin{eqnarray*}
\limsup_{N \to \infty} \; {1 \over N^{d-1}}
\max_{k_0(\gd) \leq k \leq \nu n} \, 
& \log & \Joint_N \left( u^\gz_k  \in 
\bigcap_{i =1}^{\ell} \,  \cV(\widehat R^i , \gd \vol(\widehat R^i))  
\right)
\leq\\
&& \qquad 
- \sum_{i =1}^\ell \int_{\widehat B_i} \tau_\gb (\vec{n}_i) \, d \cH_x  
+ C_4 (\gb,v)  \gd.
\end{eqnarray*}
This concludes the Proposition.
\section{Open problems}

We would like mention some open questions related to the 
$\bbL_1$-theory
\begin{enumerate}
\item
Extention of the $\bbL_1$-theory to general finite range
models and to the context of Pirogov-Sinai Theory.
\item
Proof of the Wulff construction for continuum models
in an $\bbL_1$-setting.
\item
Upgrade of the concentration properties to the Hausdorff distance,
based on more delicate versions of the minimal section argument;
some results of this sort should appear in 
\cite{BodineauIoffeVelenik99}.
\item
A more challenging problem would be to provide an accurate 
description of phase segregation \`a la DKS. In particular
one should understand how to control phase boundaries and 
prove local limit results with boundary conditions which 
are only statistically pure.
\end{enumerate}

\part{Dobrushin-Koteck\'{y}-Shlosman (DKS) theory in 2D}
\label{part_strongWulff}
\setcounter{section}{0}

In this part we review and explain the results on phase 
separation in the two-dimensional  nearest neighbor Ising model as enforced
by the canonical constraint on the magnetization \cite{DKS}, \cite{IS}. 
The theory is built upon sharp local estimates over finite volume vessels
 $\Lambda_N$ and on the probabilistic analysis of the random microscopic
phase separation line. We focus here on the ``free'' spatial geometry of
 the phase
segregation, that is disregarding the boundary effects. These effects
 could enter
the picture in two different ways: in terms of the boundary conditions on 
$\partial\gL_N$ and in terms of the geometry of $\gL_N$. 
In the former case the minority phase could be absorbed by part of the boundary
$\partial\gL_N$. This and related phenomena are discussed in Part~4.
In the second case  the finite vessel $\gL_N$ might not be able to accommodate
 the corresponding optimal crystal shape. Such a geometric constraint is, from
 the point of view of the  microscopic theory, merely a technical nuisance, 
though,
 on the macroscopic level, it might lead to formidable variational problems.
We go around this domain geometry issue by choosing $\gL_N$ to be of the Wulff
shape itself
$$
\gL_N~=~N{\mathcal K}_1\cap\bbZ^2 ,
$$
where ${\mathcal K}_1$ is the unit area Wulff shape. Thus, 
 $\gL_N$ accommodates any optimal shape of area smaller than $N^2$. 

The corresponding finite volume canonical Gibbs measure is then defined by
\begin{equation}
\label{3.1.muN}
\Is_{\gL_N,-}^{\gb}\lb~\cdot~\big|~M_N (\gs ) =-N^2m^* +a_N\rb ,
\end{equation}
where $M_N \df\sum_{i\in\gL_N}\gs_i$
is the total spin, $m^* =m^* (\gb )$ is the
spontaneous magnetization, and  
 $a_N$ points inside the phase transition region, $a_N\in (0,2N^2m^* )$.
 In the sequel we shall use the shortcut $\Is_{N,-}^{\gb}$ for the 
finite volume measure $\Is_{\gL_N,-}^{\gb}$.

\noindent
{\bf Notation.}  The values of positive constants $c_1, c_2, ...$ are 
updated with each subsection.
\section{Main Result}
\label{dima_main}
\setcounter{equation}{0}
DKS theory gives a comprehensive solution to the following problem 
of phase separation:
\vskip 0.1cm
\noindent
{\bf Problem~1.} For $\gb >\gb_c$ and $a_N\in (0,2N^2 m^* )$ characterize typical 
spin configurations $\gs$ under the canonical measure \eqref{3.1.muN} .
\vskip 0.1cm
An ostensibly simpler problem is
\vskip 0.1cm
\noindent
{\bf Problem~2.} For $\gb >\gb_c$ and $a_N\in (0,2N^2 m^* )$ find sharp local asymptotics 
of 
$$
\Is_{N,-}^{\gb} \left( M_N~=~-m^*N^2 +a_N\rb  .
$$
In fact both problems are equivalent. 
 In particular, the phenomenon behind the shift of the magnetization  is
inside the phase transition region
 not a bulk one (and hence is not in the realm of the usual theory of large deviations),
 and the crucial role is played by the spatial geometry
of symmetry breaking.

\noindent
\subsection{Heuristics} 
\label{dima_main_heuristics}
Under the finite volume pure state $\Is_{N,-}^{\gb}$ the typical maximal size of
$\pm$~contours is of order $\log N$. One could then visualize a typical microscopic 
configuration $\gs$ on $\gL_N$  in terms of an archipelago of small~(that is of the maximal
size $\sim\log N$) ``$+$'' islands which could contain still smaller ``$-$'' lakes
 etc.  This archipelago spreads out uniformly over $\gL_N$, and 
the density of the plus ``soil'', which spells out in terms of the magnetization 
$M_N (\gs )$ as $\lb |\gL_N |+M_N (\gs )\rb/2|\gL_N |$, is close to its equilibrium
value
$$
\frac{|\gL_N |+\la M_N\ran{\gb}{N, -}}{2|\gL_N |}\ \sim\ \frac{1-m^*}{2} .
$$
Thus, one could think
 of two different competing patterns behind the $a_N$-shifts, 
$a_N\geq 0$, of
the magnetization $M_N$ from its equilibrium value 
$\la M_N\ran{\gb}{N,- }\sim -m^* |\gL_N |$:
\vskip 0.cm
\noindent
1) The density of the archipelago increases in a spatially homogeneous fashion without, 
however, altering the typical sizes of the islands.

\noindent
2) Spatial symmetry is broken, and an abnormally huge island of the ``$+$''~phase 
of excess area  
$\sim a_N/2m^*$
appears.
\vskip 0.1cm
\noindent
Heuristically, the first scenario corresponds to Gaussian 
fluctuations, and its price, in terms
of probability, should be of order
$$
\exp\lb -c_1 (\gb )a_N^2/N^2\rb .
$$

Phase segregation manifests itself in the second scenario, and the probabilistic price 
for creating such a huge island is proportional to the length of its boundary
$$
\exp\lb -c_2 (\gb )\sqrt{a_N}\rb .
$$
A comparison between the two expressions above suggests that the first scenario should be
preferred whenever $a_N\ll N^{4/3}$, whereas large shifts $a_N\gg N^{4/3}$ should result in
the phase segregation picture described in the second scenario. This indeed happens to
 be the case,
and we refer to \cite{DS} and \cite{IS} for a complete rigorous 
treatment\footnote{The critical
case of $a_N\sim N^{4/3}$ is still an open problem.}. 

\begin{figure}[t]
\centerline{
\psfig{file=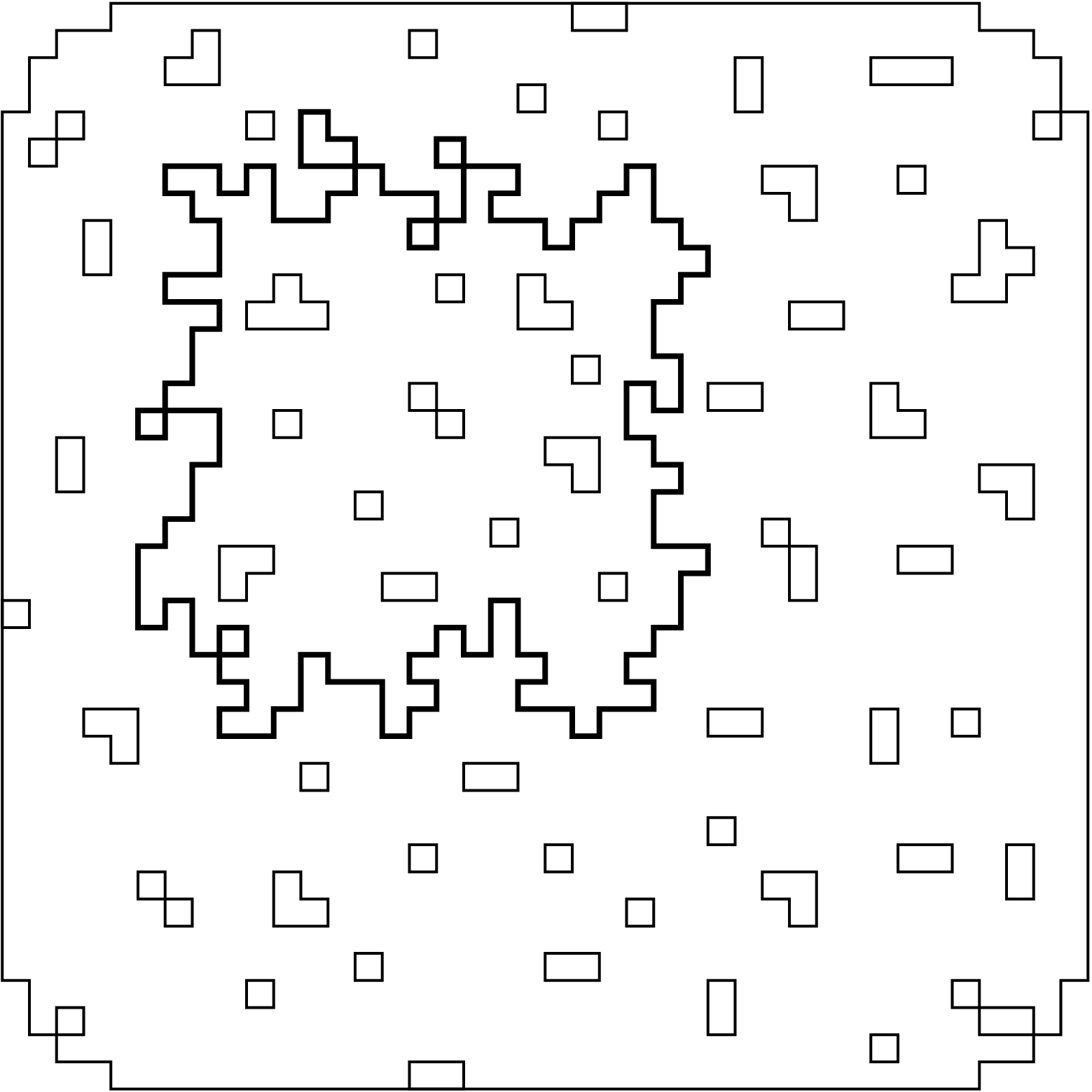,height=6cm}\hspace{3cm}
\raise 1.5 truecm \hbox{\psfig{file=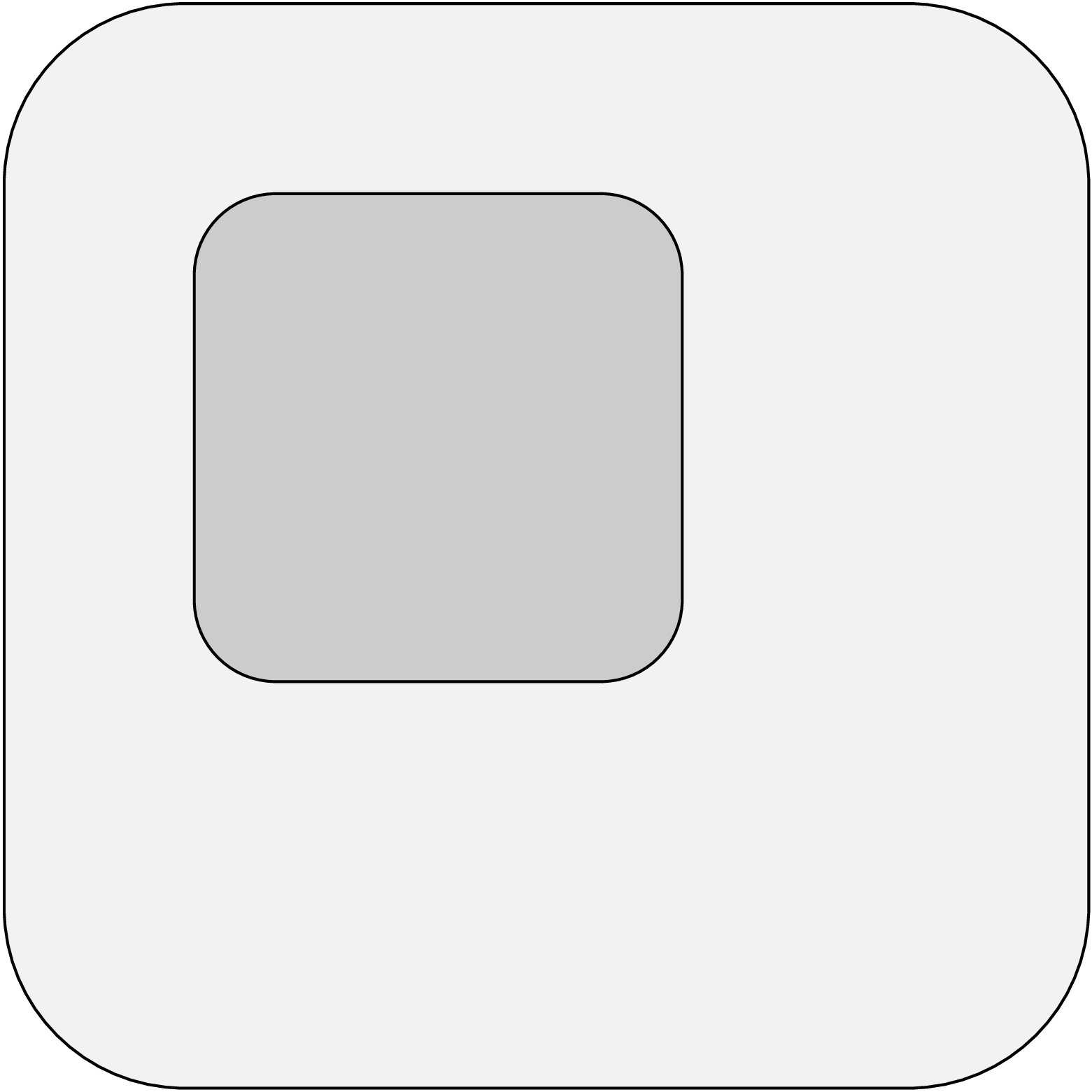,height=3cm}}
}
\figtext{ 
\writefig        1.10    3.30     {\Large $\stackrel{1/N}{\longrightarrow}$} 
\writefig        4.22    2.50     {\small $-$ phase} 
\writefig        3.63    3.70     {\small $+$ phase} 
\writefig       -2.27    2.90     {\small $\gga$} 
}
\caption{DKS picture under the  $1/N$ scaling: On the left the microscopic $\gL_N$
box with the unique $K\log N$-large contour $\gamma$.  On the right the 
continuous box $\cK_1$ with the scaled image of $\gamma$.} 
\label{fig_rescaling}
\end{figure}

For the sake of the exposition, we shall
stick here to the possibly most interesting case of $a_N\sim N^2$, which corresponds also
 to the macroscopic type of scaling discussed in Part~2.  
The DKS theory
 gives then
the following sharp characterization of the phase segregation in the canonical ensemble:
under $\Is_{N,-}^{\gb}\lb~\cdot~\big| M_N =-m^* N^2 +a_N\rb$ a typical spin configuration
 $\gs$ contains exactly one abnormally large contour $\gamma$ which decouples between
 the ``$+$''~phase (inside $\gamma$) and the ``$-$''~phase (outside $\gamma$). In 
particular, the average magnetization inside (respectively outside) $\gamma$ is
 close to $m^*$ (respectively $-m^*$), and the area encircled by $\gamma$ can be
 thus recovered from the canonical constraint,
$$
m^*\left|\text{int}\lb\gamma\rb\right| -m^*\lb N^2 -\left|\text{int}\lb\gamma\rb\right|\rb
~\approx~-m^*N^2+a_N\ \ \Longrightarrow\ \  \left|\text{int}\lb\gamma\rb\right|~\approx~
\frac{a_N}{2m^*} .
$$ 
Under the scaling of $\gL_N$ by $1/N$, that is into the normalized continuous 
shape $\cK\subset\bbR^2$,
the microscopic phase boundary $\gamma$ sharply concentrates around a shift of the Wulff
shape of the corresponding scaled area $a_N/2m^*N^2$ (Fig.~\ref{fig_rescaling}). 
\subsection{DKS theorem} 
\label{dima_main_thm}
More precisely, for any $r\in\bbR_{+}$ let ${\mathcal K}_{r}$ to denote the
Wulff shape of  the area $r$. Also given a number $s\in\bbR_+$, let us say that
a microscopic contour $\gamma$ is $s$-large, if ${\rm diam}_{\infty}(\gamma
)>s$.
\begin{thm}[\cite{DKS}\footnote{},\cite{IS}]\footnotetext{In the original
monograph [DKS] the corresponding  results has been derived in the context of
the Ising model with periodic boundary condition.}
\label{thm DKS}
Let the inverse temperature $\beta >\beta_c$ be fixed, and let  
the sequence $\{a_N\}$, $-m*N^2 +a_N\in\text{\rm Range}(M_N )$, be such  
 that the limit
$$
a~=~\lim_{N\to\infty}\frac{a_N}{N^2}~\in~(0,2m^* (\beta ))
$$
exists. Then,
$$
\log\Is^{\gb}_{N,-}\big(~M_N~=~-m^*N^2~+a_N~\big)\ =\ 
-\cW_{\gb}\lb{\partial\mathcal K}_{\frac{a_N}{2m^*}} \rb\big( 1~+~\mbox{O}\big(
 N^{-1/2}\log N\big)\big) .
$$
Moreover, if $K=K(\gb )$ is large enough, 
 with $\IsN\left(~\cdot~|M_N=-N^2m^* + a_N\right)$-probability converging to $1$ as
 $N\to\infty$:
\begin{enumerate}
\item There is exactly one 
$K(\gb )\log N$-large contour $\gamma$. 
\item This $\gamma$ satisfies
\begin{equation}
\begin{align}
\label{thm3.1.1}
&\min_{x}\frac1{N}d_{{\Bbb H}}\big(~\gamma,x+\partial {\mathcal K}_{\frac{a_N}{2m^*}}
\big)
\ \leq\ c_1(\beta )N^{-1/4}\sqrt{\log N}\\
\intertext{and}
\label{thm3.1.2}
&\min_{x}\frac1{N^2}\text{\rm Area}\Big(~\text{\rm int}\lb\gamma \rb
\Delta \lb x+{\mathcal K}_{\frac{a_N}{2m^*}}\rb~\Big)\ \leq\ c_2 (\beta )N^{-3/4}\sqrt{\log N}.
\end{align}
\end{equation}
\end{enumerate}
\end{thm}

\subsection{DKS theory}
\label{dima_main_theory}
The DKS theory views the production of the event $\{ M_N -m^*N^2 +a_N\}$ in terms of
a two-step procedure: On the first stage a length scale $s =s(N)$ is chosen, and {\bf all}
 the microscopic $s$-large contours $(\gamma_1 ,...,\gamma_n )$ are
fixed. If the total area inside these $s(N)$-large contours is smaller than $a_N/2m^*$,
 then the total magnetization $M_N$ still has to be steered towards the imposed value 
$M_N = -m^*N^2 +a_N$, but already under the constraint that all the $\pm$~contours different
from $(\gamma_1 ,...,\gamma_n )$ are $s(N)$-small. The probability 
$\IsN\lb M_N = -m^*N^2 +a_N\rb$ reflects the price of the optimal strategy along these lines.

We record the two steps of the DKS theory as follows:

\noindent
1) Study the statistics of $s(N)$-large contours under $\IsN$.

\noindent
2) Give local limit estimates on the magnetization in the $s(N)$-restricted phases.

\begin{figure}[t]
\centerline{
\psfig{file=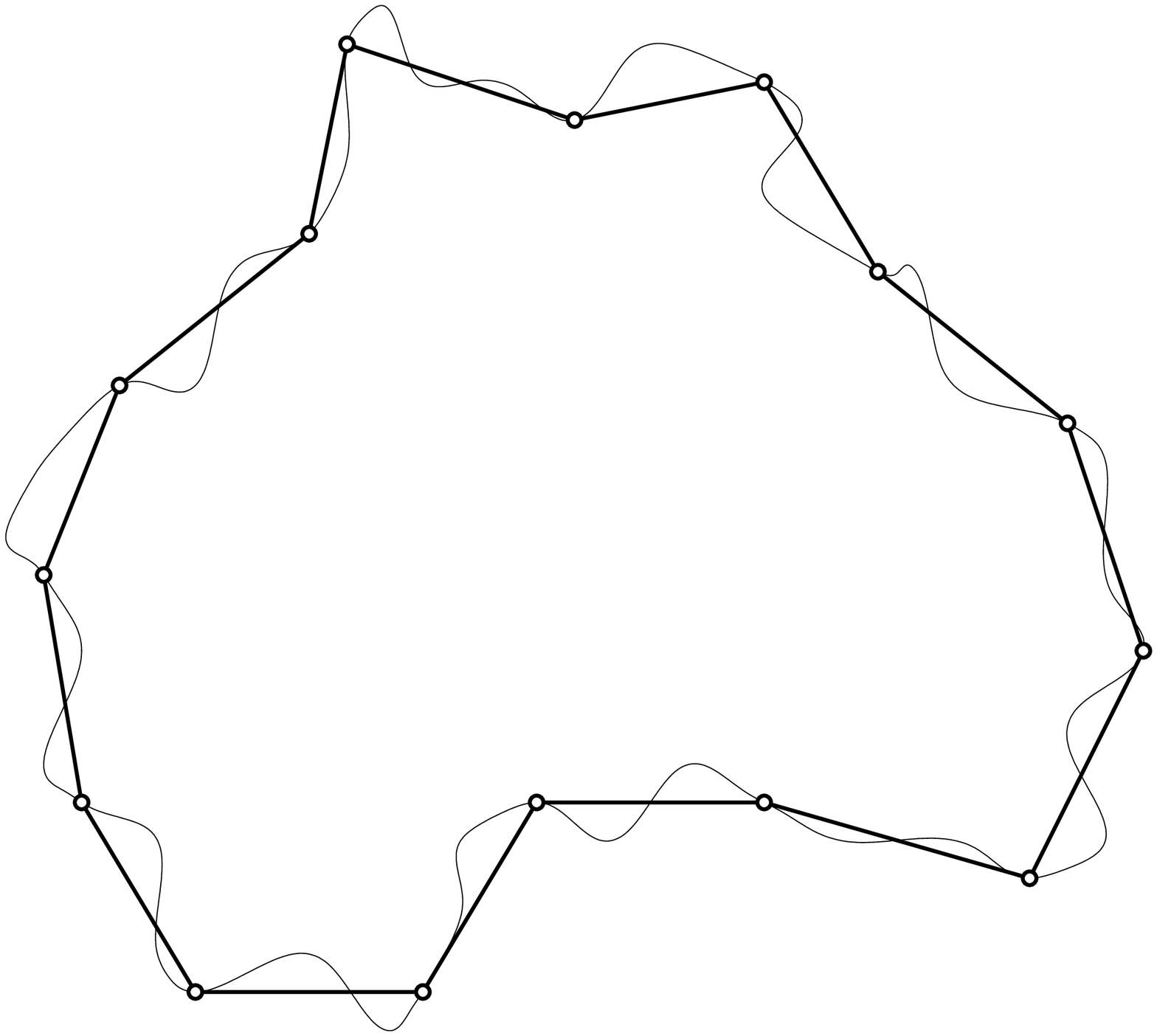,height=5cm}\hspace{1cm}
\psfig{file=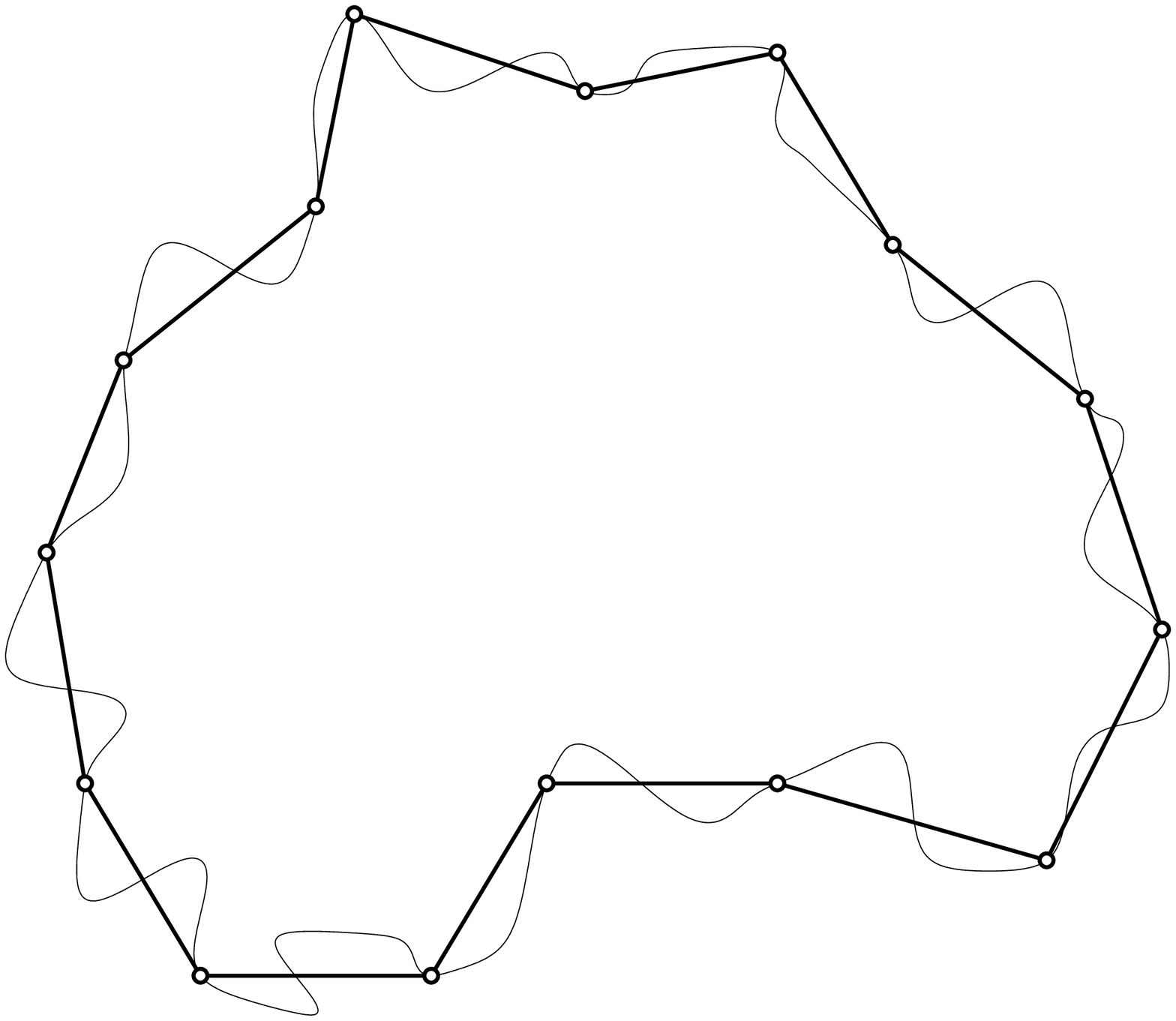,height=5cm}
}
\figtext{ 
\writefig       -6.05   3.70    {\footnotesize $u_1$} 
\writefig       -6.40   2.50    {\footnotesize $u_2$} 
\writefig       -6.20   1.40    {\footnotesize $u_3$} 
\writefig       -5.30   5.30    {\footnotesize $u_{n-1}$}
\writefig       -5.20   4.40    {\footnotesize $u_n$}
\writefig       -3.58   5.20    {\footnotesize $\gga_1$}
\writefig        2.50   4.80    {\footnotesize $\gga_2$}
}
\caption{Two microscopic contours $\gamma_1$  and $\gamma_2$ are compatible 
with the same skeleton $S= (u_1,...u_n)$.} 
\label{fig_skeletonS}
\end{figure}

\noindent 
The introduction of $s(N)$-cutoffs leads to the separation of the length scales which 
has a double impact on the problem: it sets up the stage for the renormalization 
analysis of microscopic phase boundaries, and it improves the control over the 
bulk magnetization inside the corresponding microscopic phase regions. Let us try to 
explain this in more details: As far as the statistics of the $s(N)$-large contours is
considered, we are interested in giving sharp estimates on the $\IsN$-probability of
the events of the type
$$
\lbr~\text{$s(N)$-large contours of $\gs$ encircle a certain prescribed area~}\rbr .
$$
The point is that the contribution of any particular microscopic contour to the 
probability of such an event is negligible. In other words, one also has to take into account
the entropy (number) of all the contributing contours. The required entropy cancelation
(and hence the production of the relevant limiting thermodynamic quantity - surface
 tension) is achieved by means of a certain coarse graining procedure, the so called
skeleton calculus, which we describe in Section~\ref{dima_skeletons}.
 Roughly, instead of studying the
probabilities of
individual  microscopic contours one considers the packets of all contours passing through
the vertices of a given ``$s(N)$-skeleton'' $S= (u_1, u_2,...,u_n)$ and staying within a
distance of the order $s(N)$ from the closed polygonal 
line $\text{Pol}(S)$ (Fig.~\ref{fig_skeletonS}). 
The distance between successive vertices of $S$ complies with the length scale
$s(N)$,  $\normsup{ u_{i+1} -u_{i}}\sim s(N)$. Surface tension is produced on
the level of skeletons. In fact, the probability of observing a $\pm$~contour
compatible with a given skeleton $S$ admits an asymptotic with
$s(N)\nearrow\infty$ description
\begin{equation}
\label{3.1.skeleton}
\IsN\lb S\rb\ \asymp\ \exp\lbr -\cW_{\gb}\lb \text{Pol} (S)\rb\rbr .
\end{equation}
We quote the precise result  in 
Section~\ref{dima_skeletons}, which we  devote to
 a general exposition of the skeleton calculus.

Since the vertices of $S$ are $s(N)$-apart, and the surface tension $\st$ is strictly
positive for all $\gb >\gb_c$, the energy $\cW_{\gb}\lb \text{Pol} (S)\rb$
controls the number $\# (S)$ of vertices of $S$ as
\begin{equation}
\label{3.1.snumber}
\# (S)~\leq ~c_3 (\gb )\frac{\cW_{\gb}\lb \text{Pol} (S)\rb}{s(N)} .
\end{equation}
When combined with \eqref{3.1.skeleton} this leads to the reduction of the 
combinatorial complexity of the problem: the number of different skeletons of a fixed
energy $\widehat{\cW}_N$ does not compete with the approximate probability 
$\exp\{-\widehat{\cW}_N\}$ to observe any such skeletons. Thus, the study of
 $\{ M_N =-m^*N^2 +a_N\}$ reduces, in terms of skeletons, to the maximal term estimation.
It should be stressed, however, that unlike the coarse graining procedures of the $\bbL_1$
theory, the mesoscopic objects (skeletons) of the DKS theory closely follow the
microscopic structure of phase boundaries.

The local limit estimates in the $s(N)$-restricted phases are, therefore,
required  uniformly over finite lattice domains whose boundaries are carved
with $s(N)$-large contours compatible with not too costly skeletons. This
imposes a natural restriction on the length of these boundaries, and we shall
describe the appropriate family of domains  in Section~\ref{dima_estimates}
along with the exposition of the corresponding uniform local limit  results.
Intuitively, long contours are responsible for long range dependencies between
spins, and, therefore,  the $s(N)$-cutoff constraint improves the mixing
properties of the  system and helps to extend the validity of classical
(Gaussian) behavior of moderate deviations. In Section~\ref{dima_bulk} we
quote the corresponding relaxation and decay properties which lie in the heart
of the local limit estimates. In Section~\ref{dima_structure} we give an
outline of the proof of the DKS theorem.

Finally, the (long) list of open problems is briefly addressed in
Section~\ref{dima_problems}.

\section{Estimates in the phases of small contours}
\label{dima_estimates}
\setcounter{equation}{0}
 As it has been mentioned, the 
estimates in the phase of small contours should be derived uniformly over a family
 of lattice domains whose boundaries are composed of not too costly $s(N)$-large contours.

\noindent
{\bf Definition} Basic family $\cD_N$ of subsets $A\subseteq\gL_N$: We fix two numbers
$a$ (small) and $R$ (big).
$$
A\in\cD_N~\Longleftrightarrow~aN^2\leq |A| \quad\text{and}\quad |
\partial A|\leq RN\log N.
$$
\qed

We fix a basic scale $s(N) =K\log N$ of large contours, where $K= K(\gb )$ is a
sufficiently large number, so that $K\log N$-contours are highly improbable under 
the pure state $\Is_{N,-}^{\gb}$. Of course, exactly the same number $K$ appears in 
the statement of Theorem~\ref{thm DKS}. The upper bound on $\partial A$ in the definition
of the family $\cD_N$ states that the configurations with total length of $K\log N$
 large contour exceeding $RN\log N$ are ruled out. This conclusion is explained in
 more detail in Section~\ref{dima_skeletons} (see the remark following Lemma~\ref{energy}).

\medskip
\noindent
\subsection{Structure of local limit estimates} 
\label{dima_estimates_structure}
Let us turn now to the structure of local limit estimates in the $s(N)$-restricted 
phases. First of all, given any $A\subset\bbZ^2$, the $s$-restricted phase
 on $A$ is 
defined via
$$
\Is_{A,-}^{\gb,s}\left(~\cdot ~\right)~\df ~
\Is_{A,-}^{\gb}\left(~\cdot ~\Big|\text{All $\pm$ contours are $s$-small}\right).
$$
We would like to study the probabilities of deviations $a_N\geq 0$ of the total
 magnetization $M_A$ from the 
corresponding averaged value $\langle M_A\rangle_{A,-}^{\gb ,s}$. Let us define
 the set of feasible values of such deviations as
$$
{\bf M}_A^{+}~= ~\left\{ a_N\geq 0:\ \langle 
 M_A\rangle_{A,-}^{\gb ,s} +a_N\in\text{\rm Range}(
 M_A)\right\} .
$$
Roughly, the cutoff $s$ extends the validity of Gaussian moderate deviations for
 the following reason: The price of shifting the magnetization by $a_N$ on the 
expense of $s(N)$-small contours is of the order $(a_N /s^2)s\sim a_N/s(N)$. This should
be tested against the Gaussian moderate deviation exponent of the order $a_N^2 /N^2$.
 Thus the Gaussian behavior should prevail once $a_N\ll N^2/s(N)$. Of course, the latter
constraint on $a_N$ becomes less stringent as $s(N)$ decreases. 
On the rigorous mathematical part the classical approach to estimating
$$
\Is_{A,-}^{\gb,s}\left( M_A =\langle M_A\rangle_{A,-}^{\gb ,s} +a_N\right),
$$
amounts to first finding the value of magnetic
field
$$
g = g(A, s(N),a_N ),
$$ 
such that the expected magnetization under the $g$-tilted state is precisely
what we want, 
\begin{equation}
\label{3.2.sgexp}
\langle M_A\rangle_{A,-,g}^{\gb ,s} \ =\ \langle M_A\rangle_{A,-}^{\gb,s} ~+~a_N,
\end{equation}

\noindent
and, then, to rewrite the $\Is_{A,-}^{\gb ,s}$-probability in terms of the 
$\Is_{A,-,g}^{\gb ,s}$ one:
\begin{equation}
\label{3.2.transform}
\begin{split}
&\Is_{A,-}^{s}\left( M_A=\langle M_A\rangle_{A,-}^{\gb ,s}+a_N\right)~\\
&\qquad =~\exp\big\{-(\langle M_A\rangle_{A,-}^{s} +a_N)g~+~
\log\big<\text{e}^{g M_A}\big>_{A,-}^{\gb ,s}\big\}~
\Is_{A,-,g}^{\gb ,s}\big(~ M_A~=~\big< M_A\big>_{A,-,g}^{\gb ,s}~\big)\\
&\qquad=~
\exp\left\{-\int\limits_0^g\int\limits_r^g\big< M_A ; M_A\big>_{A,-,h}^{\gb ,s}
 \text{d}h\text{d}r\right\}
~\Is_{A,-,g}^{\gb ,s}\left(~ M_A~=~\big< M_A\big>_{A,-,g}^{\gb ,s}\right) .
\end{split}
\end{equation}
One then tries to derive sufficiently precise estimates on the semi-invariants of 
$\Is_{A,-,h}^{\gb ,s}$ and to prove a local CLT under $\Is_{A,-,g}^{\gb ,s}$.
 Thus, it is extremely important to 
understand how the magnetization $\langle M_A\rangle_{A,-,g}^{\gb,s}$ and 
other semi-invariants of $\Is_{A,-,g}^{\gb,s}$ change with the magnetic field
$g$ in the phase of $s(N)$-small contours. 

 Breaking of the classical limit behavior 
in the $s(N)$-restricted phase 
manifests itself by the jump of the magnetization which is 
related to the appearance of
abnormally large $\pm$-contours. Without cutoffs this jump occurs for $g\sim 1/N$,
 and imposing the $s(N)$ constraint would delay such a jump \cite{ScS2}.
 It is easy to imagine 
what
should be the critical order of the magnetic field  $g$, at which those large
contours should start to be favored in the $s$-restricted phase: for a $\pm$~contour 
of the linear
size $s(N)$ one wins $\sim s^2g$ on the level of magnetization and loses 
$\sim s$ on the level of surface energy. These two terms start to
be comparable when $sg\sim 1$.  Therefore no particular deviation 
from the classical behavior should be expected as far as  $gs(N)\ \ll\ 1$.
 We refer to \cite{IS}, where all these heuristic considerations 
 have been made precise. 

\noindent
\subsection{Basic local estimate on the $K\log N$ scale}
\label{dima_estimates_basic}
Actually~\cite{IS}  it is enough to consider
 only the basic $K\log N$-scale:

\begin{lem}[\cite{IS}]
\label{KlogN}
Assume that a sequence of numbers $\{ b_N\}$ satisfies
$$
\lim_{N\to\infty}\frac{b_N\log N}{N^2} ~= ~0 .
$$
Then, on the basic scale $s(N)= K\log N$, the estimate 
\begin{equation}
\label{3.2.KlogN}
\begin{split}
&\Is_{A,-}^{\gb ,s}\left( M_A=\langle M_A\rangle_{A,-}^{\gb ,s}+a_N\right)\\
&\qquad =\  \frac1{\sqrt{2\pi\chi_{\gb} |A|}}
\exp\big\{-\frac{a_N^2}{2\chi_{\gb} |A|}+\mbox{O}
\bigl(\frac{a_N^2}{N^3}(\log N\vee\frac{a_N}{N})\bigr)\big\}
\big(1~+~\text{\small{o}}(1) \big),
\end{split}
\end{equation}
holds uniformly in domains $A\in\cD_N$ and in $a_N\in{\bf M}_A^{+}\cap [0,b_N]$,
where $\chi_{\gb}$ is the susceptibility under the pure state $\Is_{-}^{\gb}$.
\end{lem}

\noindent
\subsection{Super-surface estimates in the restricted phases}
\label{dima_estimates_supersurface} 
Moderate deviations on the intermediate scales $s(N)\gg\log N$ are, for the 
purposes of the theory, controlled by the following super-surface order estimate
in the phase of small contours (c.f. Lemma~2.5.1 in \cite{IS})

\begin{lem}
\label{supersmall}
Let the large contour parameter $s(N)\gg \log N$ be fixed.
There exists a constant $c_1 =c_1 (\beta )>0$, such that for 
 all $N>0$,  $A\in\cD_N$ and all $a_N\in{\bf M}_A^{+}$, 
\begin{equation}
\label{supper}
\Is_{A,-}^{\gb ,s}
\big(~M_A = \langle M_A\rangle_{A,-}^{\gb ,s} +a_N~\big)\ \leq\ 
\exp\big(~-c_1\frac{a_N^2}{N^2}\wedge
\frac{a_N}{s(N)}~\big) .
\end{equation}
\end{lem}
The idea of the proof is simple: either an area of order $a_N /2m^*$ is exhausted by
 the $K\log N$ large contours, which, in the $\Is_{N,-}^{\gb s}$-restricted phase,
 should have a surface tension price with the exponent of the order $a_N /s(N)$, or
$K\log N$ large contours cover an area much less than $a_N /2m^*$, which means that
the remaining deficit of the magnetization should be compensated in the basic
 $K\log N$ restricted phase, where we can use Lemma~\ref{KlogN}.

\section{Bulk Relaxation in Pure Phases}
\label{dima_bulk}
\setcounter{equation}{0}
  The term relaxation is used here in the 
equilibrium setting in order to describe the approximation of local 
finite volume statistics by the infinite
volume ones. We successively describe the relaxation properties of 
pure ``$-$'' states with non-positive
 and small positive magnetic fields and in the restricted phases of small contours.

\noindent
\subsection{Non-positive magnetic fields $h\le 0$.}
\label{dima_bulk_nonpositive}
  The crucial property of 
low temperature pure phases could
be stated as follows: Let us say that the sites $i$ and $j$ are $*$-neighbors
if $\| i-j\|_1 =1$. Given a spin configuration $\sigma$ on 
$\{-1,+1\}^{\bbZ^2}$,
 let us say that the sites $i$ and $j$ are $+*$-connected, if there exists a
$*$-connected chain of sites $i_1,...,i_n$, $i_1=i$ and $i_n =j$, such that
$\sigma (i_k)=1$ for every $k=1,...,n$. 
\begin{thm}[ \cite{CCSc}]
\label{ccs}
For every $\gb >\gb_c$ there exists $c_{1} =c_{1} (\gb ) >0$, such that 
uniformly in subsets $A\subseteq \bbZ^2$, 
$i,j\in A$  and in 
magnetic fields $h\le 0$,
\begin{equation}
\label{3.2.decay}
\Is_{A,-,h}^{\gb}\left(~i\stackrel{+*}{\longleftrightarrow}j~\right)~
\le~{\rm e}^{-c_1 (\gb )\normsup{i-j}}.
\end{equation}
\end{thm}  
\noindent
{\it Remark.} Of course, since
$\left\{i\stackrel{+*}{\longleftrightarrow}j\right\}$ is a non-decreasing
event, the uniformity follows from the FKG ordering, once \eqref{3.2.decay} is
verified for the infinite volume zero-field measure $\Is_{-}^{\gb}$.
\begin{cor}[Relaxation of local observables]
Fix $k\in\bbZ$.  
Uniformly in $A\subseteq \bbZ^2$, magnetic fields $h\le 0$ and local
observables $f$ with $|{\rm supp}(f)|=k$,
\begin{equation}
\label{3.2.spin}
\left|\langle f\rangle_{A,-,h}^{\gb}~-~\langle f\rangle_{-,h}^{\gb}\right|~\le~
c_2 (k){\rm e}^{-c_3(\gb ) {\rm dist}_{\infty}\big({\rm supp}(f),\partial
A\big)}
\end{equation}
\end{cor}
\noindent
Furthermore,
\begin{cor}[Relaxation and decay of semi-invariants]
Fix $n\in\bbZ$.  Uniformly in $A\subseteq \bbZ^2$, magnetic fields $h\le 0$ and
sites $i_1 ,...,i_n\in A$,
\begin{equation}
\begin{align}
\label{3.2.semiA}
&\left|\langle\gs (i_1);...;\gs(i_n )\rangle_{A,-,h}^{\gb}~-~\langle\gs (i_1);...;\gs(i_n )
\rangle_{-,h}^{\gb}\right|~\le~c_4 (n)
{\rm e}^{-c_5 (\gb ){\rm dist}_{\infty}\big( \{i_1 ,...,i_n\},\partial A\big)}\\
\intertext{and}
\label{3.2.semidecay}
&\left|\langle\gs (i_1);...;\gs(i_n )\rangle_{A,-,h}^{\gb}\right|~\le ~c_6 (n){\rm exp}\left\{-c_7(\gb )
\frac{{\rm diam}_{\infty}\big( i_1, ...,i_n\big)}{n}\right\} .
\end{align}
\end{equation}
\end{cor}
\noindent
Finally,
\begin{cor}[Asymptotic expansions]
Fix $n\in\bbZ$. 
Uniformly in $A\subseteq\bbZ^2$ and in $i\in A$,
\begin{equation}
\label{3.2.expansion1}
\left|\langle\gs (i)\rangle_{A,-,h}^{\gb}~-~\big( -m^*(\gb )+
\sum_{k=1}^{n}\frs_k \frac{h^k}{k!}\big)
\right|~\le ~c_8 (n) |h|^{n+1} +c_9(n){\rm e}^{-c_{10} (\gb )
{\rm dist}_{\infty}\big( i,\partial A\big)} ,
\end{equation}
where $\frs_k$ is the $k$-th semi-invariant of the zero-field infinite volume
measure $\Is_{-}^{\gb}$,
$$
\frs_k ~\df~\sum_{i_1 ,...,i_k\in\bbZ^2}\langle \gs (0) ;
\gs (i_1);...;\gs(i_n )\rangle_{-}^{\gb} .
$$
\end{cor}
\noindent
{\bf  Remark}  It is possible (and straightforward) to 
formulate \eqref{3.2.semiA}, \eqref{3.2.semidecay} and 
\eqref{3.2.expansion1} 
 in the general case of $n$ local  observables $f_1,...,f_n$.\qed

\noindent
\subsection{Positive magnetic fields $h > 0$.}
\label{dima_bulk_positive}
  Modifying ``$-$'' 
states by negative magnetic fields $h<0$ 
amounts to moving
away from  the phase transition region.  Relaxation properties of $\Is_{A,-,h}^{\gb}$ with $h>0$ are 
radically different - uniformity is lost, and the size of the domain $A$ starts to play a crucial role. Indeed, the 
unique infinite
volume measure $\Is_{-,h}^{\gb}=\Is_{h}^{\gb}$ stochastically dominates $\Is_{+}^{\gb}$ whatever
small $h>0$ is.  Thus, for large domains $A$, the configuration in the bulk is flipped under 
$\Is_{A,-,h}^{\gb}$ into
the ``$+$'' dominated state.  It is easy to understand on the heuristic grounds what should be
 the order of the critical
size of $A$ for such a ``flip'' to occur:  given $h>0$, the surface energy of a $\pm$-contour $\gamma$ is of 
the order $|\gamma |$ and it competes with the bulk gain inside the contour which, in its turn, is proportional
to $h{\rm Area}(\gamma )$.  The latter factor wins (loses), once the linear size of 
$\gamma$ is much larger (respectively much smaller) than $1/h$.  Thus the sign of the
dominant spin under  $\Is_{A,-,h}^{\gb}$ should depend on whether $A$ can accommodate large enough 
contours, or, in other words, on how the linear size of $A$ relates to $1/h$. 

The important and remarkable fact is  that exponential  relaxation properties of finite volume ``$-$'' states are
uniformly preserved for domains of the sub-critical size.

\begin{thm}[\cite{ScS2},~\cite{IS}]
\label{thmh}
There exists a constant $a=a(\gb )>0$ such that for any $h>0$ fixed,  
\begin{equation}
\label{3.2.decayh}
\Is_{A,-,h}^{\gb}\left(~i\stackrel{+*}{\longleftrightarrow}j~\right)~
\le~{\rm e}^{-c_1 (\gb )\normsup{i-j}}.
\end{equation}
uniformly in domains $A\subset\bbZ$ such that
any connected component of $A$ has diameter bounded above by $a/h$. As a consequence exponential decay of
semi-invariants \eqref{3.2.semidecay} and the asymptotic expansion estimate \eqref{3.2.expansion1} hold 
uniformly in such domains as well. 
\end{thm}
 
\noindent
\subsection{Phases of small contours}
\label{dima_bulk_phases}
 Theorem~\ref{thmh} explains how the cutoff parameter
$s(N)$ upgrades the regular behavior of ``$-$''-states with positive magnetic fields
 $h$: By the definition of the restricted phase $\Is_{A,-}^{\gb ,s}$ the diameter of 
any relevant microscopic domain is at most of the order $s(N)$. 
\begin{thm}[\cite{ScS2},~\cite{IS}]
\label{thmhs}
There exists a constant $a=a(\gb )>0$ such that for any $h>0$ and $s$ satisfying $hs\leq a(\gb )$,  
\begin{equation}
\label{3.2.decayhs}
\Is_{A,-,h}^{\gb ,s}\left(~i\stackrel{+*}{\longleftrightarrow}j~\right)~
\le~{\rm e}^{-c_1 (\gb )\normsup{i-j}}\,,
\end{equation}
uniformly in domains $A\subseteq\bbZ$ . \\
Furthermore, the expectations in restricted phase are
controlled as follows: for every $k\in\bbZ$,
\begin{equation}
\label{3.2.spins}
\left|\langle f\rangle_{A,-,h}^{\gb, s}~-~\langle f\rangle_{A\cap\gL_s
(f),-,h}^{\gb}\right|~\le~
c_{2} (k){\rm e}^{-c_{3}(\gb )s},
\end{equation}
uniformly in $A\subseteq \bbZ^2$  and in local functions $f$, 
$\left|\big({\rm supp(f)}\big)\right| =k$, 
where we have used the following notation: 
$\gL_s (f)\df \left\{i:{\rm d}_{\infty}\left( i,{\rm supp}(f)\right)\leq s\right\}$.
Finally, the decay of the semi-invariants is controlled in the restricted phases as
\begin{equation}
\label{3.2.semidecays}
\left|\langle\gs (i_1);...;\gs(i_n )\rangle_{A,-,h}^{\gb,s}\right|~\le 
~c_4 (n){\rm exp}\left\{-c_5(\gb )
\frac{{\rm diam}_{\infty}\big( i_1, ...,i_n\big)}{n}\wedge s\right\} .
\end{equation}
\end{thm} 

\section{Calculus of Skeletons}
\label{dima_skeletons}
\setcounter{equation}{0}
The renormalization analysis of large $\pm$~contours is performed on various cutoff
scales $s$, the appropriate choice of $s$ typically depending on the linear size $N$ of the
system $s = s(N)$. We shall state coarse graining estimates uniformly in finite domains 
$A\subset\bbZ^2$ and in the cutoff scales $s$.

\noindent
\subsection{Definition}
\label{dima_skeletons_definition}
 A $\pm$~contour $\gamma$ is said to be $s$-large if $\text{diam}_\infty(\gamma
)\geq s$. Given a cutoff scale  $s\in\bbN$ and an $s$-large  $\pm$~contour
$\gamma$ we say that $S =(u_1,...,u_n)$ is an $s$-skeleton of $\gamma$,
$\gamma\sim S$ if
\begin{enumerate}
\item All vertices of $S$ lie on $\gamma$.
\item $s(N)/2 \leq \normsup{u_i - u_{i+1}}\leq 2s,\ \forall~i=1,...,n$, where we
have identified $u_{n+1}\equiv u_1$. 
\item The Hausdorff distance $d_{{\mathbb H}}$ between $\gamma$ and the polygonal
line ${\rm Pol} (S)$ through the vertices of $S$ satisfies
$$
d_{{\Bbb H}}\big(\gamma ,{\rm Pol} (S)\big)\ \leq\ s(N) .
$$
\end{enumerate}
Similarly, given the collection $\lb \gamma_1 ,...,\gamma_n\rb$ of all $s$-large contours
of a configuration $\gs\in\gO_{A,-}$, let us say that a collection  $\frS =(S_1 ,...,S_n) $
of $s$-large skeletons is compatible with $\gs$, $\gs\sim\frS$, if $\gamma_i\sim S_i$ for all
 $i=1,...,n$.

Of course, a configuration $\gs\in\gO_{A,-}$ has, in general, many different 
compatible collections of $s$-skeletons. Nonetheless, for each particular $\frS$
the probability
\begin{equation}
\label{3.4.frS}
\Is_{A,-}^{\gb}\lb\frS\rb ~\df ~\Is_{A,-}^{\gb}\lb\gs :~\gs\sim \frS\rb
\end{equation}
is well defined. 

\noindent
\subsection{Energy estimate}
\label{dima_skeleton_energy} 
As the renormalization scale $s$ grows, the probabilities \eqref{3.4.frS} start to admit
 a sharp characterization in terms of the energies $\cW_{\gb}(\frS)$,
$$
\cW_{\gb }\left(\frS \right)~\df ~\sum_1^n\cW_{\gb }\left( {\rm Pol}(S_i) \right) ,
$$
for a collection $\frS  =\lb S_1 ,...,S_n \rb$. Below we a give precise version of this
crucial statement in terms of the upper and lower bounds on the corresponding probabilities. 
 The first important renormalization energy
estimates  could be \cite{Pfister} formulated as follows
\begin{lem}[\cite{Pfister}]
\label{energy}
On every skeleton scale $s$ and independently of $A\subset\bbZ^2$,
\begin{equation}
\label{3.4.energy}
\Is_{A,-}^{\gb}\big(~\frS~\big)\ \leq\ \exp\big\{~-\cW_{\gb}(\frS)~\big\} .
 \end{equation}
Furthermore, uniformly in $A\subset\bbZ$ , $r >0$ and cutoff parameters $s$,
\begin{equation}
\label{3.4.energyge}
\Is_{A,-}^{\gb}\left( \cW_{\gb}(\frS)\geq r\right)~\leq~{\rm exp}\left\{ -r\big(
1-\frac{c_1\log |A |}{s}\big)\right\} .
\end{equation}
\end{lem}
Energy estimate \eqref{3.4.energy} provides an upper bound on 
the probability of observing
$\pm$~contours in the vicinity of a skeleton. Before going 
to a complementary lower bound let us dwell on the sample path
structure of the contours which is hidden behind these renormalization
estimates.

\noindent
\subsection{Calculus of skeletons} 
\label{dima_skeleton_calculus}
By definition a contour is a self-avoiding closed path of nearest neighbor
 bonds of $\bbZ^2$. For every set $A\subseteq \bbZ^2$ the Ising measure
 $\Is_{A,-}^{\gb}$
 induces a weight function $q_{A^*}^{\gb^*}$ on the space of such self-avoiding
polygons (see Subsection~\ref{ssec_2dIsing}),
$$
q_{A^*}^{\gb^*}\left( \gamma\right)~=~\Is_{A,-}^{\gb}\lb \gs\in\gO :~
\gamma\ \text{is a $\pm$~contour 
of }\gs\rb .
$$
In terms of these weights the probability of observing a certain skeleton
$S =\lbr u_1 ,...,u_n\rbr$ could be written as

$$
\Is_{A,-}^{\gb}\lb S\rb ~=~\sum_{\gamma\sim S} q_{A^*}^{\gb^*}\left( \gamma\right ) .
$$

Each microscopic contour $\gamma$ compatible with $S$, $\gamma\sim S$, splits
into the  union of disjoint open self-avoiding lattice paths 
$\gamma_k: u_k\to u_{k+1}, \ k=1,...,n$. The analysis of limit 
properties of $\Is_{A,-}^{\gb}$ comprises two main steps which 
could be loosely described as follows:

\noindent
1) As the renormalization scale $s$ grows, the statistical behavior of
 different pieces $\gamma_k$ decouple under $q_{A^*}^{\gb^*}$, that is
\begin{equation}
\label{3.4.gasplit}
\sum_{\gamma\sim S}  q_{A^*}^{\gb^*}\left( \gamma\right)\ 
\approx\ \prod_{k=1}^{n}\lb \sum_{\gamma_k : u_k\to u_{k+1}}
q_{A^*}^{\gb^*}\left( \gamma_k\right)\rb .
\end{equation}

\noindent
2) The $k-th$ term ($k=1,...,n$) in the above product corresponds
to a $\pm$~interface stretched in the direction of the vector
 $u_{k+1}-u_{k}\in\bbR^2$, in other words
\begin{equation}
\label{3.4.taubk}
q_{A^*}^{\gb^*}\left( \gamma_k\right) \ \approx\ \text{e}^{-\tau_{\gb}
(u_{k+1}-u_{k})} .
\end{equation}

Thus, the skeleton calculus resembles a refined version of
the sample path large deviation principle for genuinely two-dimensional
random curves. At very low temperatures, a very precise local analysis of the
phase separation line 
 has been developed in \cite{DKS},\cite{DS} using the method of cluster
expansions. Our approach here pertains
to the whole of the phase transition region $\gb >\gb_c$, but is strongly
linked to the very specific self-duality properties of the two-dimensional
nearest neighbor Ising model. 
We refer to Subsection~\ref{ssec_2dIsing} and, eventually, to 
\cite{PfisterVelenik97,PfisterVelenik98} for comprehensive description and
study of the relevant properties of the duality transformation. The output of
these techniques could be recorded in the following form

\begin{lem}[Probabilistic Structure of the Phase Separation Line
\cite{PfisterVelenik97}]
\label{psline}
Given any $A\subset\bbZ^2$ and any two compatible self-avoiding paths $\gl_1$ and 
$\gl_2$,
\begin{equation}
\label{3.4.qasuper}
q_{A^*}^{\gb^*}\lb\gl_1\cup\gl_2 \rb~\geq ~q_{A^*}^{\gb^*}\lb\gl_1\rb q_{A^*}^{\gb^*}\lb\gl_2\rb .
\end{equation}
Furthermore, 
\begin{equation}
\label{3.4.expbouns}
\text{e}^{-c_1 (\gb )|\gl_2 |}~\leq ~\frac{q_{A^*}^{\gb^*}\lb\gl_1\cup\gl_2 \rb}{q_{A^*}^{\gb^*}\lb\gl_1\rb}~
\leq~\text{e}^{-c_2 (\gb )|\gl_2 |} 
\end{equation}
On the other hand, given any $A\subseteq \bbZ^2$ and any three points $u,v,w\in A^*$,
 the $q^{\gb^*}_{A^*}$ weight of the paths going from $u$ to $v$ through $w$ is
 bounded above as \cite{PfisterVelenik97}
\begin{equation}
\label{3.4.qasub}
\sumtwo{\gl :u\to v}{w\in\gl} q_{A^*}^{\gb^*}\lb \gl\rb ~\leq~
\lb\sum_{\gl_1 :u\to w}q_{A^*}^{\gb^*}\lb \gl_1\rb\rb
\lb\sum_{\gl_2 :w\to v}q_{A^*}^{\gb^*}\lb \gl_2\rb\rb .
\end{equation}
Finally, the weights $q_{A^*}^{\gb^*}$ are non-increasing in $A$, and are related to the dual connectivities
as
\begin{equation}
\label{3.4.qarepr}
\sum_{\gl :~ u\to v}q_{A^*}^{\gb^*}\lb\gl\rb ~=~\la \gs(u) \gs (v)\ran{\gb^*}{A^*,f} .
\end{equation}
\end{lem}
Relation \eqref{3.4.qarepr} is the link to the surface tension: first of all the impact of
a particular set $A$ exponentially diminishes with the distance to $\partial A$
\cite{I1},
\begin{equation}
\label{3.4.decay}
\la\gs (u) \gs (v)\ran{\gb^*}{f}-\exp\lbr - c_2 (\gb )\text{d}\lb\{ u,v\},
\partial A\rb\rbr~\leq~
 \la\gs (u) \gs (v)\ran{\gb^*}{A^*,f}~\leq ~
\la\gs (u) \gs (v)\ran{\gb^*}{f} .
\end{equation}
uniformly in  $A^*\subseteq\bbZ^2$ and 
any $u,v\in A^*$. 
Moreover the following Ornstein-Zernike type correction formula \cite{Al} holds 
uniformly in $u,v\in\bbZ^2$:
\begin{equation}
\label{3.4.OZ}
\exp\left\{-\tau_{\gb}\lb u-v\rb -c_3 (\gb )\log \normsup{u-v}\right\}~\leq ~
\la\gs (u) \gs (v)\ran{\gb^*}{f}~\leq~
\exp\left\{-\tau_{\gb}\lb u-v\rb\right\},
\end{equation}

\noindent
\subsection{Skeleton lower bound}
\label{dima_skeletons_lb}
The energy estimate \eqref{3.4.energy} is an immediate consequence of the
(iterated) sub-multiplicative property \eqref{3.4.qasub}, the representation 
formula \eqref{3.4.qarepr} and the right-most inequalities in \eqref{3.4.decay}
and \eqref{3.4.OZ}. In order to prove a lower bound one essentially needs to
reverse the inequality in \eqref{3.4.qasub}. An indirect way to do so is to use
the FK representation (see \cite{ScS1} and \cite{IS}). We shall briefly present
here a more direct approach which has been developed in \cite{I1} and
\cite{PfisterVelenik97}. Qualitatively it gives the same order of corrections
as the FK one, but has a clear advantage of being explicitly related to the
statistics of the microscopic phase boundaries at different length scales.  The
basic idea is that the phase separation line has rather strong mixing
properties, in particular  paths $\gl_1$ and $\gl_2$ on the right hand side of
\eqref{3.4.qasub} should interfere, in the  case of $(u,v,w)$ being in a
general position,  only in a vicinity of $w$. Thus, at a price of lower order
corrections (as we shall see these corrections are logarithmic with the
skeleton scale $s$) the inequality \eqref{3.4.qasub} could  be reversed using
the  super-multiplicativity property \eqref{3.4.qasuper}. The notion of
``general position'' simply means that $u,w$ and $v$ do not form too small an
angle and live on the same length scale, and  it is quantified by the following

\noindent
{\bf Definition.} Given a skeleton scale $s\in\bbN$ and a number $\gep >0$, let us say that
 that a triple $(u,w,v)$ of $\bbZ^2$-lattice points is $(s,\gep )$-compatible, if
$$
\frac{s}2 ~\leq ~\min\lbr \normsup{ w-u},\normsup{ v-w}\rbr ~\leq
\max\lbr \normsup{ w-u},\normsup{ v-w}\rbr ~\leq ~ 2s ,
$$
whereas $\cos\lb w-u ,v-w\rb \geq -1+\gep $.\qed

We shall state the lower bound in terms of the limiting weights 
$q^{\gb^*}\lb\cdot\rb\df\lim_{A^*\nearrow\bbZ^2_\star}q^{\gb^*}_{A^*}$ 
(which exist by Lemma~\ref{psline}).
\begin{lem}
\label{triple}
Fix $\gep >0$. Then there exists a scale $s= s(\gep )$, such that
\begin{equation}
\label{3.4.triple}
\sumtwo{\gl :u\to v}{w\in\gl} q^{\gb^*}\lb \gl\rb~\geq~\exp\lbr -\lb
\st (w-u)+\st (v-w )\rb -c_1 (\gb )\log s\rbr ,
\end{equation}

uniformly in all skeleton scales $s\geq s(\gep )$ and in all $(s,\gep )$-compatible
 triples $(u,w,v)$.
\end{lem}
       
We sketch the proof of this lemma in Appendix~B. Iterating \eqref{3.4.triple} 
we arrive to the following lower bound on the probability of observing
a certain regular skeleton:

\noindent
{\bf Definition.} A skeleton $S =(u_1 ,...,u_n)$ is said to be $(s,\gep )$-regular, if
any triple $(u_{i-1},u_i ,u_{i+1})$ of  successive points of $S$ is $(s,\gep )$-compatible,
and the distance between any two non-neighboring intervals $[u_i ,u_{i+1}]$ and 
$[u_j ,u_{j+1}]$ exceeds $\gep s$.\qed
\begin{lem}
\label{skeletonlb}
For every $\gep >0$, there exists a number $c_2 =c_2 (\gep )<\infty$, such that 
uniformly in the skeleton scales $s$
and in all $(s, \gep )$-regular skeletons $S$, 
\begin{equation}
\label{3.4.lb}
\begin{split}
\Is_{N, -}^{\gb}&\left( \exists~\text{a}~\pm~\text{contour}~\gamma :\ 
d_{\bbH}(\gamma , {\rm Pol}(S) )\le K(\gb )\sqrt{s}\log s \right) \\
&\ge~
{\rm exp}\left\{ -W_{\gb} \left({\rm Pol}(S) \right) ~-~c_2 (\gep )\# (S) \log s\right\}\\
&\ge ~
{\rm exp}\left\{ -W_{\gb} \left({\rm Pol}(S) \right)\left(1 ~-~c_3 (\gep,\gb )\frac{\log s}{s}
\right)\right\} ,
\end{split}
\end{equation}
where $\#  (S)$ denotes the number of vertices in $S$, and the last inequality follows from 
\eqref{3.1.snumber}. 
\end{lem} 

In fact we need lower bounds only for a very specific set of $s$-skeletons, namely on those
approximating the Wulff shape $\cK_{a_N /2m^*}$. These skeletons always satisfy the conditions
of the above theorem. An academic attempt to prove a lower bound for all possible shapes will
lead to annoying, though solvable, technicalities, but will fail to contribute much to the 
microscopic theory of phase separation, as we see it.

\section{Structure of The Proof}
\label{dima_structure}
\setcounter{equation}{0}
In order to give a probabilistic characterization of the microscopic canonical state
$\Is_{N,-}^{\gb}\left( ~\cdot ~\big| M_N =-m^*N^2 +a_N\right)$ one first derives a sharpest
possible lower bound on the probability $\Is_{N,-}^{\gb}\left(M_N =-m^*N^2 +a_N\right)$, 
and then rules 
out those geometric events (in terms of skeletons, but with an eventual translation to
the language of microscopic spin configurations), which happen to qualify as improbable 
when compared with this lower bound.

\noindent
\subsection{Lower bound}
\label{dima_structure_lb}
 The best lower bound comes as an outcome of the optimal 
combination of the basic local limit Lemma~\ref{KlogN} and the skeleton lower bound
 \eqref{3.4.lb}. We choose a skeleton approximation of the corresponding Wulff shape
$\cK_{a_N/2m^*}$, and using local limit estimates steer the magnetization towards the
desirable value $-m^*N^2 +a_N$.
 Optimality reflects the choice of the best possible skeleton scale:
Notice that the estimate \eqref{3.4.lb} becomes sharper with the growth of the cutoff 
parameter $s(N)$. On the other hand, the area of the microscopic phase region is 
controlled, with respect to the area inside ${\rm Pol}(S)\sim a_N/2m^*$, up to a 
$N\sqrt{s(N)}\log s(N)$ correction (see Appendix~B or \cite{IS}), which, of course,
makes the local limit step more expensive for large values of $s(N)$. It happens that
the bounds are balanced on the skeleton scale $s(N)\sim \sqrt[4]{a_N}$.

\begin{thm}[\cite{IS}]
\label{lb}
Uniformly in $a_N\in {\bf M}_N^{+}$, that is for all $a_N\geq 0$, such that
$-m^*N^2+a_N\in{\rm Range}(M_N)$,
\begin{equation}
\label{3.5.lb}
\Is_{N,-}^{\gb}\left(M_N =-m^*N^2 +a_N\right)~\geq ~\exp\left\{ 
-\sqrt{\frac{a_N}{2m^*}}\cW_{\gb}\left(\partial\cK_1\right) -
c_1 (\gb )\sqrt[4]{a_N}\log N
\right\} .
\end{equation}
\end{thm}

\noindent
\subsection{Upper bounds}
\label{dima_structure_ub}
 First of all, one derives an upper bound on the shift
of the magnetization. On any skeleton scale,
\begin{equation}
\label{3.5.sdecomp} 
\Is_{N,-}^{\gb}\left(M_N =-m^*N^2 +a_N\right)~\leq ~\sum_{\frS}
\Is_{N,-}^{\gb}\left(M_N =-m^*N^2 +a_N~;~\frS \right) .
\end{equation}
Due to the intrinsic entropy cancelation under the skeleton coarse graining,
and in view of the lower bound \eqref{3.5.lb} and the energy estimate
\eqref{3.4.energy} one could, for example, shoot for the maximal term in the
above sum. If the phase volume (see \cite{DKS}  for the precise definition ) of
$\frS$ is much less than $a_N/2m^*$, then the deficit of the magnetization
should be compensated in the phase of $s(N)$-small contours, which, by 
Lemma~\ref{supersmall} exerts a super-surface price in the exponent. On the
other hand, if the phase volume of $\frS$ is close to $a_N/2m^*$, then by the
isoperimetric inequality  and by the energy estimate \eqref{3.4.energy}, the
best possible price one should be prepared  to pay is already close to ${\rm
exp}\left\{ -\cW_{\gb}\left(\cK_{a_N/2m^*}\right)\right\}$.  Again the
resulting estimate is subject to an optimization via a careful choice of the
skeleton scale $s(N)$.

\begin{thm}[\cite{IS}]
\label{ub}
Uniformly in $a_N\sim N^2$,
\begin{equation}
\label{3.5.ub}
\Is_{N,-}^{\gb}\left(M_N =-m^*N^2 +a_N\right)~\leq ~\exp\left\{ 
-\sqrt{\frac{a_N}{2m^*}}\cW_{\gb}\left(\partial\cK_1\right) +c_1 (\gb )\sqrt[4]{a_N}\log N
\right\} .
\end{equation}
\end{thm}

A more delicate study \cite{DKS},\cite{IS} of the typical sample properties of
the microscopic configuration $\gs$ under  $\Is_{N,-}^{\gb}\left( ~\cdot ~\big|
M_N =-m^*N^2 +a_N\right)$ is again based on the analysis of
\eqref{3.5.sdecomp}. At this point the stability Bonnesen-type estimates (see
Subsection~1.3 of the Introduction)  for the Wulff variational problem become
important - they enable to quantify the  conclusion that only those collections
$\frS$, which are close to the shifts of  the Wulff shape $\cK_{a_N/2m^*}$,
have a chance to survive a comparison with the lower bound \eqref{3.5.lb}. A
step further, involving local limit estimates of Lemma~\ref{KlogN}, is to
conclude that all these collections actually contain exactly one large
skeleton, which corresponds to the unique large contour as asserted by the DKS
theorem.

\section{Open Problems}
\label{dima_problems}
\setcounter{equation}{0}

There are still important open problems even in the nearest neighbor Ising
case. Notably, one knows how to control precise fluctuations of the phase
separation line only at very low temperatures, that is using the method of
cluster expansions \cite{DH}. This is a serious gap in the theory, since large
scale statistics of  microscopic phase boundaries are ultimately responsible
for exact (up to  zero order terms) expansions of canonical partition functions
\cite{H}.  So far qualitative probabilistic results have been obtained either
for very low temperature models \cite{H}, or in the simplified setting of
self-avoiding polygons \cite{I3}, \cite{HI} or Bernoulli bond percolation
\cite{CI}. Another interesting and apparently important problem is to
understand sample path properties of spin configurations in a situation when a
canonical constraint is imposed in the  restricted phase. Apart from giving
rise to a potentially fascinating probabilistic structure, this question is
closely related to the issue of the dynamical spinodal decomposition.

There is absolutely no matching probabilistic study of the phase separation in 
multiphase two-dimensional models, for example $q$-states Potts models. Some
results in this direction are reported in \cite{Velenik97}, but this issue is
almost entirely open even in the context of the $\bbL_1$-theory. In particular,
the corresponding phenomena is still not worked out on the level of macroscopic
variational problems, see, however \cite{ABFH}, \cite{MoS} and the references
 therein. 

The key issue, however, which we feel is largely misunderstood is that 
at moderately low temperatures the DKS
theory of two-dimensional phase segregation, say in the general context of
finite range ferromagnetic models with pair interactions is far from being complete.
What currently exists is an
example of how these ideas could be implemented in the nearest
neighbor case. At least from the mathematical point of view, the nearest
neighbor case is a degenerate one, in a sense that it enables a reduction
 to pure boundary conditions over decoupled microscopic regions even at 
temperatures only moderately below critical. This should not be the case for more
general range of interactions. 
In this respect the assertion that low temperature
expansions should go through for general interactions much along the same lines
 as they do for the nearest neighbor model, seems to be rather irrelevant - the
real issue is not to kill mixed boundary conditions, but to understand how they
 should be incorporated into the DKS theory.

\part{Boundary effects}\label{part_boundary}
\setcounter{section}{0}
In the previous parts, we explained how the thermodynamical variational
problem describing the macroscopic geometry of coexisting phases can be derived
in various lattice models of statistical physics. To simplify the analysis, we
restricted our attention to periodic boundary conditions or to systems
contained in a Wulff-shaped box, avoiding thus a discussion of the effect of a
confining geometry on the behavior of the system. In this part, we would
like to explain what happens when we take such effects into account. Boundary
conditions play a particularly important role in the kind of problems presented
in this review, since they concern the asymptotic behavior of large but finite
systems and therefore the boundary cannot be simply ``sent to infinity'' as
usually done. We will see that taking care of boundary effects not only
provides a complete description of the geometry of these constrained systems
thus allowing a rigorous description of the interaction between an equilibrium
crystal and a substrate, but also allows to study the effect of so-called {\em
boundary phase transitions}.

For simplicity, we only discuss the case of the Ising model with nearest
neighbors interaction.
\section{Wall free energy}
\setcounter{equation}{0}
The vessel containing the system has not only the property of confining it, but
can also act in an asymmetric way on the various phases inside, favoring
some of them; indeed this is what happens typically in real systems. In fact,
this is precisely the reason one introduces boundary conditions in the first
place: To impose the equilibrium phase the system realizes. It appears
to be convenient to have a parameter allowing a fine-tuning of the asymmetry,
interpolating between pure $+$ or $-$ boundary conditions. Let us now describe
how this is done.

\bigskip
Let $\theplane = \setof{i\in\bbZ^d}{i(d) = 0}$ and $\halfspace =
\setof{i\in\bbZ^d}{i(d)\geq 0}$. The vessel of our system is the box
\begin{equation*}
\BoxxNM = \setof{i\in\halfspace}{-N\leq i(n) \leq N,\, n=1,\dots,d-1,\,
0\leq i(d) \leq M}\,,
\end{equation*}
and the {\em wall} is $\thewall = \BoxxNM\cap\theplane$. 

Let $\bdf\in\bbR$; we consider the following Hamiltonian,
\begin{equation*}
\Ham_\BoxxNM^\bdf(\gs) = - \sumtwo {\nnb ij\subset\halfspace}
{\nnb ij\cap\BoxxNM\neq\eset} \gs_i\gs_j - \bdf\sum_{i\in\thewall} \gs_i\,.
\end{equation*}
Let $\overline\gs\in\{-1,1\}^\halfspace$; the Gibbs measure in $\BoxxNM$ with
boundary condition $\overline\gs$ is the following probability measure on
$\{-1,1\}^\halfspace$~\footnote{\label{footnote}Note that we could equivalently consider
$\Isbd{\overline\gs}{\BoxxNM}$ as a probability measure on $\{-1,1\}^{\bbZ^d}$
by extending the b.c. $\overline\gs$ by $\overline\gs_i=1$ for all
$i\in\bbZ^d\setminus\halfspace$; it is then possible to replace the boundary
magnetic field $\bdf$ by a coupling constant: $\bdf\sum_{i\in\thewall} \gs_i =
\bdf\sum_{\nnb ij:\,i\in\thewall,\,j\not\in\halfspace} \gs_i\gs_j$. This will be
used when dealing with negative boundary field, see Subsection~\ref{ssec_2D}.},
\begin{equation*}
\Isbd{\overline\gs}{\BoxxNM}(\gs) =
\begin{cases}
(\PFbd{\overline\gs}{\BoxxNM})^{-1}\exp[-\gb\,\Ham_\BoxxNM^\bdf(\gs)] & \text{if
$\gs_i=\overline\gs_i$, $\forall i\not\in\BoxxNM$,}\\
0 & \text{otherwise.}
\end{cases}
\end{equation*}
We'll usually use the short-hand notations $\Isbd{\overline\gs}{N,M}$,
$\PFbd{\overline\gs}{N,M}$, ....
As usual, we write $+$ for $\overline\gs\equiv 1$ and $-$ for
$\overline\gs\equiv -1$. We therefore distinguish one of the sides of the box
$\BoxxNM$, $\thewall$, which we call the ``wall''. Notice that instead
of usual boundary conditions, a {\em boundary magnetic field} $\bdf$ is
acting on $\thewall$; since setting $\bdf=1$ produces $+$ b.c. on the
wall, while setting $\bdf=-1$ results in $-$ b.c., this provides the promised
interpolation parameter. Of course, we could also consider more complicated
situations, where (possibly inhomogeneous) boundary magnetic fields act on the
whole boundary of the box. However, for simplicity, we restrict our attention
to this particular case, which will turn out to be general enough that the
basic phenomena induced by the use of boundary fields can already be analyzed.

\medskip
To quantify the preference of the wall toward one of the phases, it is
convenient to introduce a new thermodynamic quantity, the {\em wall free
energy},
\begin{equation}\label{eq_taubd}
\taubd(\gb,\bdf) \df \limtwo{N\ra\infty}{M\ra\infty} \frac 1{\abs{\thewall}}\log
\frac{\PFbd{+}{N,M}}{\PFbd{-}{N,M}}\,.
\end{equation}
The existence of this quantity, and the remarkable fact that the two limits can
be taken in any order, has been established in \cite{FroehlichPfister87a}; the
proof relies on the simple identity
\begin{equation}\label{eq_taubd2}
\taubd(\gb,\bdf) = \limtwo{N\ra\infty}{M\ra\infty} \gb \int_{-\bdf}^\bdf \frac
1{\abs{\thewall}} \sum_{i\in\thewall} \Ebdf{+}{N,M}{\bdf'}{\gs_i}\,\dd \bdf'\,.
\end{equation}
We'll return to this formula in the next section. The heuristics behind the
definition of $\taubd(\gb,\bdf)$ is that the free energy $\Fbd{+(-)}{N,M} = -\log
\PFbd{+(-)}{N,M}$ of the $+$ ($-$) phase can be decomposed in the
following way:
\begin{align*}
\Fbd{+}{N,M} &= f_{\scriptscriptstyle\rm
b}(\beta) \,\abs{\BoxxNM} + f^+_{\scriptscriptstyle\rm s}(\beta)
\,\abs{\bnd\BoxxNM\setminus\thewall} + f^+_{\scriptscriptstyle\rm
w}(\beta,\bdf)\,\abs{\thewall} +
\smallo(\abs{\bnd\BoxxNM},\abs{\thewall})\,,\nonumber\\
\Fbd{-}{N,M} &= f_{\scriptscriptstyle\rm
b}(\beta) \,\abs{\BoxxNM} + f^-_{\scriptscriptstyle\rm s}(\beta)
\,\abs{\bnd\BoxxNM\setminus\thewall} + f^-_{\scriptscriptstyle\rm
w}(\beta,\bdf)\,\abs{\thewall} + \smallo(\abs{\bnd\BoxxNM},\abs{\thewall})\,,
\end{align*}
where 
\begin{align*}
f_{\scriptscriptstyle\rm b}(\beta) &\df -\lim_{N,M\ra\infty} \abs{\BoxxNM}^{-1}\,
\log\PFbd{\overline\gs}{N,M}\,, \\
f^+_{\scriptscriptstyle\rm s}(\beta)&\df -\lim_{N,M\ra\infty}
\abs{\bnd\BoxxNM}^{-1}\, \bigl(\log \PFbdf{+}{N,M}{1} - f_{\scriptscriptstyle\rm
b}(\beta) \abs{\BoxxNM}\bigr)\,, \\
f^+_{\scriptscriptstyle\rm w}(\beta,\bdf) &\df -\lim_{N,M\ra\infty}
\abs{\thewall}^{-1}\, \bigl(\log\PFbd{\overline\gs}{\BoxxNM} - 
f_{\scriptscriptstyle\rm b}(\beta) \abs{\BoxxNM}-f^+_{\scriptscriptstyle\rm s}
(\beta) \abs{\bnd\BoxxNM\setminus\thewall}\bigr)\,,
\end{align*}
(and similarly for $f^-_{\scriptscriptstyle\rm s}(\beta)$ and
$f^-_{\scriptscriptstyle\rm~w}~ (\beta,\bdf)$). As the notations suggest,
$f_{\scriptscriptstyle\rm b}(\beta)$ is independent of $\bdf$ and
$\overline\gs$, $f^+_{\scriptscriptstyle\rm s}(\beta)$ is independent of $\bdf$
and by symmetry $f^+_{\scriptscriptstyle\rm s}(\beta) =
f^-_{\scriptscriptstyle\rm s}(\beta)$. Therefore, we see that
$\taubd(\beta,\bdf) = \limtwo{N\ra\infty}{M\ra\infty} \frac 1{\abs{\thewall}}\,
(\Fbd{-}{N,M}-\Fbd{+}{N,M}) = f^-_{\scriptscriptstyle\rm w}(\beta,\bdf) -
f^+_{\scriptscriptstyle\rm w}(\beta,\bdf)$ is nothing else than the leading
order term of the difference in free energy between the two phases in the
presence of the wall.

The ultimate justification of \eqref{eq_taubd} however is that this quantity
plays exactly the role of its thermodynamical analogue in the variational
problem describing the macroscopic geometry of phase coexistence, see Theorems
\ref{thm_sessile2D} and \ref{thm_sessile3D} below.

The following Theorem states basic properties of $\taubd(\gb,\bdf)$; since
$\taubd(\gb,\bdf)$ is obviously odd in $\bdf$, we just state them for $\bdf\geq
0$ (also $\taubd(\gb,0)=0$). 
\begin{thm}\label{thm_taubdprop} \cite{FroehlichPfister87b} Let
$\tau^*_\gb=\tau_\gb(\vec{e}_d)$ and suppose $\bdf\geq 0$. Then
\begin{itemize}
\item $\taubd(\gb,\bdf)$ is a non-negative, increasing function of $\gb$ and
$\bdf$, concave in $\bdf$; moreover, if $\bdf>0$,
$$
\taubd(\gb,\bdf) >0 \Leftrightarrow \gb > \gbc\,.
$$ 
\item For all $\gb$ and $\bdf$, $\taubd(\gb,\bdf)\leq \tau^*_\gb$.
\item For all $\gb>\gbc$, there exists $1\geq\hw(\gb)>0$ such that
$$
\taubd(\gb,\bdf)<\tau^*_\gb \Leftrightarrow \bdf<\hw(\gb)\,.
$$
\end{itemize}
\end{thm}
In the case of the 2D Ising model, $\hw(\gb)$ can be computed explicitly, see
\cite{Abraham80, McCoyWu73} and Fig.~\ref{fig_droplets}.

The following terminology is standard\footnote{This terminology only makes
sense once we have chosen one of the equilibrium phase as reference; here it is
the $-$ phase.}: when $\bdf\geq\hw(\gb)$, we say that the system is in the
{\em complete drying} regime; when $\abs{\bdf}<\hw(\gb)$, it is in the {\em
partial wetting} regime; and when $\bdf\leq-\hw(\gb)$, it is in the {\em
complete wetting} regime. The reason for this terminology should become clear
later.
\section{Surface phase transition}
\setcounter{equation}{0}
In this section, we will see that the boundary magnetic field can trigger
{\em surface phase transitions}: The behavior of the system in the vicinity of
the wall depends dramatically on $\abs\bdf$ being greater or smaller than
$\hw(\gb)$. A more detailed discussion of these issues can be found in
\cite{PfisterVelenik96}.

\bigskip
The state of the system in the middle of a big box $\BoxxNM$ is entirely
determined by the boundary conditions, and is independent of the value of the
boundary field, so that the usual (infinite volume) Gibbs state simply doesn't
provide any information on the behavior of the system close to the wall. To
analyze the behavior of the system ``in the vicinity'' of the wall, it is
therefore useful to introduce the notion of {\em surface Gibbs states}; these
differ from the Gibbs states usually considered in these models by the fact
that one does not work with a sequence of boxes converging to $\bbZ^d$, but
instead converging only to the half-space $\halfspace$. More precisely, the
surface Gibbs states are the weak limits of the measures
$\Isbd{\overline\gs}{N,M}$ when $N,M\ra\infty$ (observe that
$\BoxxNM\nearrow\halfspace$). Two of them are of particular importance for our
discussion, $\sgp$ and $\sgm$, obtained respectively by taking weak limits of
the measures with $+$ and $-$ boundary conditions. It is not difficult to show
\cite{FroehlichPfister87a} that these two measures exist, are extremal, and are
invariant under translations parallel to the wall; moreover, there is
uniqueness of the surface Gibbs state if and only if $\sgp=\sgm$.

There is a close relation between $\taubd(\gb,\bdf)$ and the behavior of the
system near the wall; this can be most easily seen from the following identity,
consequence of \eqref{eq_taubd2} and symmetry \cite{FroehlichPfister87a},
\begin{equation}\label{eq_taubd3}
\taubd(\beta,\bdf) = \int_{-\bdf}^\bdf \Ebdf{+}{\halfspace}{\bdf'}{\gs_0}\;\dd
\bdf' = \int_0^\bdf \bigl( \Ebdf{+}{\halfspace}{\bdf'}{\gs_0} -
\Ebdf{-}{\halfspace}{\bdf'}{\gs_0} \bigr)\,\dd \bdf'\,.
\end{equation}
Using \eqref{eq_taubd3}, it is possible to prove the following Theorem showing
that a surface phase transition occurs at $\bdf=\hw(\gb)$; this is the so-called
{\em wetting transition}.
\begin{thm}\label{thm_wettingtransition} \cite{FroehlichPfister87b}
There is a unique surface Gibbs state if and only if $\abs\bdf\geq\hw(\beta)$.
\end{thm}
Let us briefly discuss the heuristics behind this result. The $+$ and $-$
boundary conditions fix the phase present in the bulk (i.e. in the middle of a
big box $\BoxxNM$). However, Theorem \ref{thm_wettingtransition} shows that
when $\bdf\geq\hw(\gb)$, the surface Gibbs state is unique, and therefore the
state of the system near the wall is {\em independent} of the boundary
conditions, i.e. of the phase present in the bulk. The mechanism responsible
for this is the following. Suppose that $\bdf<0$ and consider $+$-boundary
conditions; then it is natural to regard the boundary field as a negative b.c.,
and therefore to introduce an open contour with boundary $\bnd\thewall$
separating the $-$ phase favored by the wall from the $+$ phase present in the
bulk (see Section \ref{sec_tools} for more details). As long as $\bdf>-1$, there
is a competition between two effects: On the one hand it is energetically
favorable for the open contour to follow the wall, on the other hand this
would lead to a loss in entropy, since there is less room for fluctuations.
When $\bdf\leq-\hw(\gb)$, the entropy wins: The contour is repelled away from
the wall, at a distance diverging with the size of the box; this is the
phenomenon of {\em entropic repulsion}. The surface Gibbs state then describes
the behavior of the system below this surface, i.e. a mesoscopic film of $-$
phase along the bottom wall. The fact that the contour is sent away from the
wall explains why we recover the surface tension,
$\taubd(\gb,\bdf)=\tau^*_\gb$. When $\bdf>-\hw(\gb)$ energy wins, and this
modifies completely the behavior of the microscopic surface: it sticks to the
wall, making only small excursions away from it; in this case, the phase in the
bulk can reach the wall and the surface Gibbs state depends on the choice of
boundary conditions.

Part of these heuristics can be made quite precise in the 2D case. Consider $+$
boundary conditions. When $0>\bdf>-\hw(\gb)$, one can prove that the
probability that a connected piece $I$ of the wall is not touched by the open
contour is bounded above by $K\exp[-(\tau^*_\gb-\taubd(\gb,\bdf))\,\abs{I}]$,
showing that the phase separation line really sticks to the wall
\cite{PfisterVelenik97}. The informations available when $\bdf\leq-\hw(\gb)$ are
much less precise; the magnetization profile computed in \cite{Abraham80} shows
that there is a film of width of order $\sqrt N$ along the wall. A related,
much more precise result, which holds at sufficiently low temperature and for
$\bdf=-1$ is that the phase separation line, once suitably rescaled, converges
weakly to the Brownian excursion \cite{Dobrushin93}; this should be true for
any $\bdf\leq-\hw(\gb)$.

In higher dimensions, much less is known. When $\bdf>-\hw(\gb)$, one can show
that the probability that the open contour touches the middle of the wall is
bounded away from $0$ uniformly in the size of the box
\cite{FroehlichPfister87b}. When $\bdf\leq-\hw(\gb)$, very little is
known,except in the simpler case of SOS models. Also, if it is known in
dimension 2 that $\hw(\gb)<1$ (since the exact expression for $\hw(\gb)$ has
been computed \cite{Abraham80}), this is an open problem in higher dimensions.

\medskip
Theorem \ref{thm_wettingtransition} gives a first explanation of the
terminology introduced above: when the system is in the complete drying regime,
the equilibrium phase along the wall is the $+$ phase, whatever the phase in
the bulk is; when there is complete wetting, it is the $-$ phase; only in the
regime of partial wetting can both phases be present near the wall. The fact
that the phase transition is determined by $\hw(\gb)$ (i.e. the
characterization of the partial wetting regime by
$\tau^*_\gb=\abs{\taubd(\gb,\bdf)}$) is known as {\em Cahn's criterion}.
\section{Derivation of the Winterbottom construction}
\setcounter{equation}{0}
In this section, we show how Winterbottom construction, describing the
equilibrium shape of a crystal in the presence of an attractive substrate, can
be recovered from a microscopic theory. To do this, we consider the measure
$\Isbd{+}{N,rN}$, for some $r\in\bbR$, conditioned with some canonical
constraint (exact or approximate, see below). Of course, the situation here is
more complicated than the one described in the introduction, since instead of
an infinite wall, the system is contained in a finite vessel. This, of course,
makes the problem more difficult: When the solution of the Winterbottom
variational problem does not fit inside the box $\rBox \df
\setof{x\in\bbR^d}{\abs{x(n)}\leq 1,\, n=1,\dots, d-1,\, 0\leq x(d) \leq r}$,
the solution of the constrained problem will differ from Winterbottom shape. In
fact, the general solution of the constrained problem is not known. In the way
we state them below, the derivation of this variational problem from
statistical mechanics still applies in the case when the solution is not known.

\begin{figure}[t]
\centerline{\psfig{figure=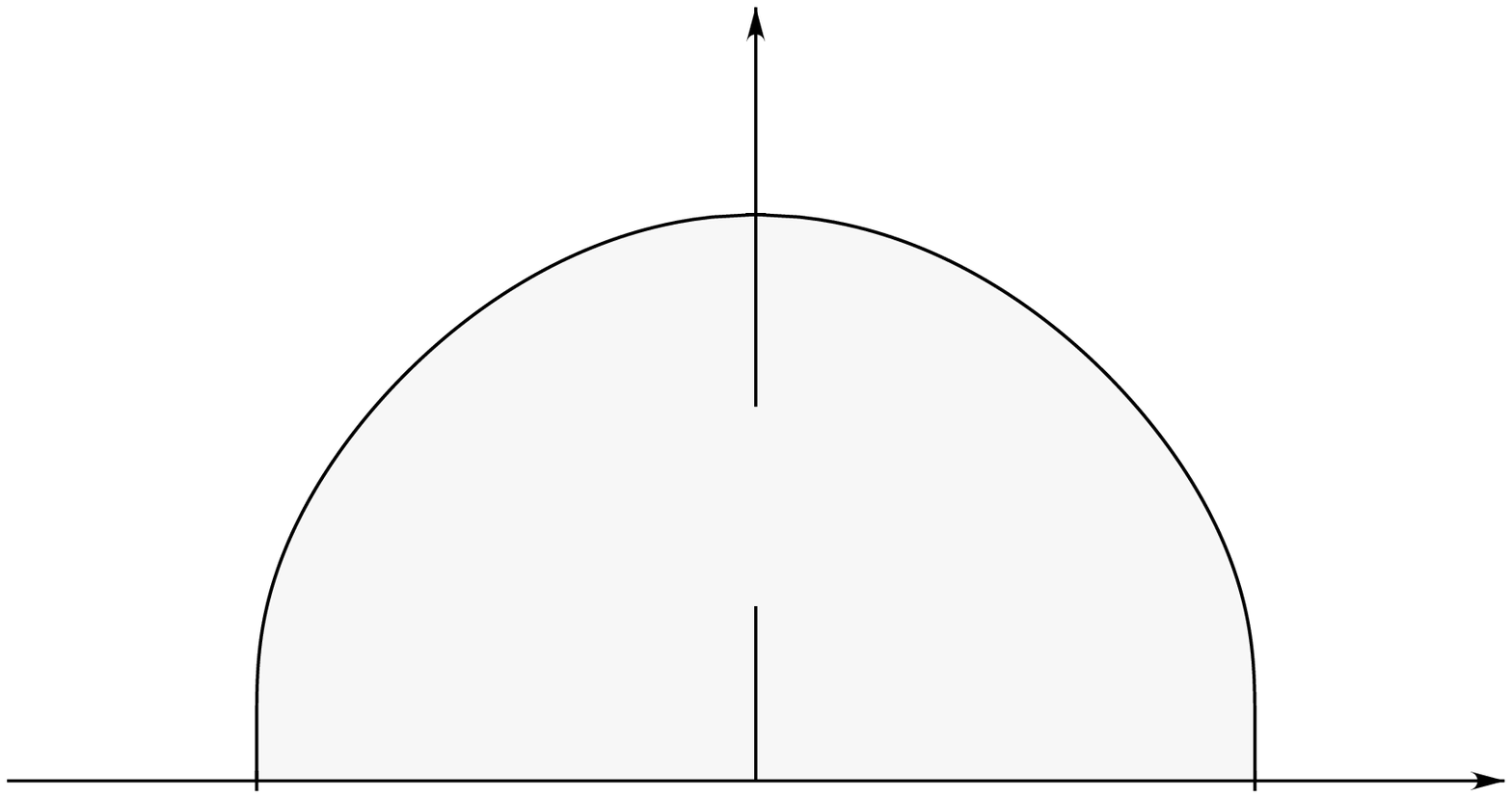,height=40mm}\hspace{.8cm}
\psfig{figure=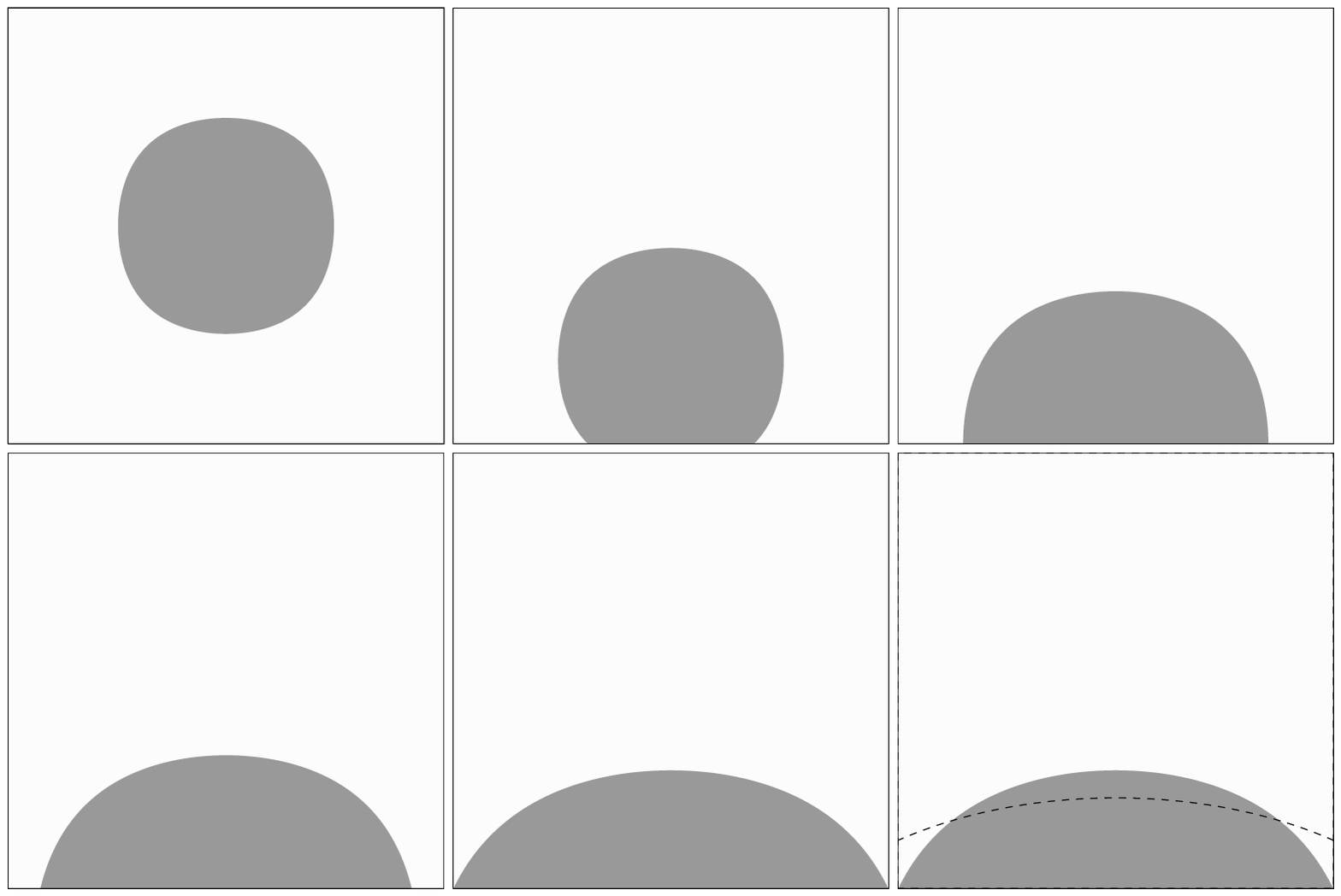,height=40mm}}
\figtext{ 
\writefig        1.40   4.17    {\footnotesize a}
\writefig        3.40   4.17    {\footnotesize b} 
\writefig        5.40   4.17    {\footnotesize c}
\writefig        1.40   2.17    {\footnotesize d} 
\writefig        3.40   2.17    {\footnotesize e}
\writefig        5.40   2.17    {\footnotesize f}
\writefig       -3.80   4.20    {\footnotesize $T$}
\writefig       -3.85   3.50    {\footnotesize $T_c$} 
\writefig        0.00   0.20    {\footnotesize $\bdf$}
\writefig       -1.04   0.20    {\footnotesize $1$} 
\writefig       -6.27   0.20    {\footnotesize $-1$}
\writefig       -4.78   1.85    {\footnotesize (non-uniqueness of}
\writefig       -4.80   1.50    {\footnotesize surface Gibbs state)}
\writefig       -4.50   2.20    {\footnotesize Partial wetting}
} 
\caption{The case of the 2D Ising model. Left: The phase diagram; the region of non-uniqueness of the surface
Gibbs state is shaded. In the other region, there is a single surface Gibbs
state. Right: A sequence of equilibrium shapes.} 
\label{fig_droplets}
\end{figure} 
Before stating the main Theorems of this Part, we briefly describe how the
wetting transition manifests itself in the macroscopic geometry of phase
separation. To do this, let $\gb>\gbc$ be fixed, and choose a value $m$ for the
canonical constraint so that the corresponding Wulff shape is small enough to
be placed inside the box $\rBox$. If $\bdf\geq\hw(\gb)$, then
$\taubd(\gb,\bdf)~=~\tau^*_\gb$, and the typical configurations will consist of
a macroscopic droplet of $-$ phase, with Wulff shape, immersed in a background
of $+$ phase; in particular, the shape of the droplet is independent of the
value of the boundary field (Fig.~\ref{fig_droplets} a). This behavior persists
up to the value $\bdf=\hw(\gb)$. Notice that as soon as $\bdf<1$, it becomes
{\em energetically} more favorable for the droplet to touch the wall. In
dimension 2, however, since $\hw(\gb)<1$, the droplet stays away from the wall,
because entropy loss is not compensated by energy gain until $\bdf$ reaches the
value $\hw(\gb)$. It is an interesting open problem to decide whether
$\hw(\gb)=1$ for $d>2$. When $\bdf<\hw(\gb)$, the typical
configurations consist of a macroscopic droplet, with Winterbottom shape, tied
to the wall. The shape of the droplet now depends on the value of $\bdf$, and
decreasing the boundary field amounts to letting the droplet spread more and
more (Fig.~\ref{fig_droplets} b--e). For some value $\widetilde\bdf$, the
droplet covers for the first time the entire wall (Fig.~\ref{fig_droplets} e).
From this point on, the shape of the droplet is left unchanged when $\bdf$ is
decreased (Fig.~\ref{fig_droplets} f; the dashed line represent part of a possible
``true'' equilibrium shape for the unconstrained problem).

From this discussion, we see that the wetting transition at $\hw(\gb)$ has a
macroscopic manifestation in the canonical ensemble. Because of the confined
geometry, however, the second transition, at $\bdf=-\hw(\gb)$ cannot be seen. To
be able to detect it, one has to consider mesoscopic droplets (in the form of
large moderate deviations, see the remark after Theorem~\ref{thm_sessile2D}).

\medskip
This also explains pretty well the terminology introduced previously: In the
complete drying regime, the droplet stays away from the wall, and so the wall is
completely dry w.r.t. the $-$ phase; in the partial wetting regime, the droplet
touches the wall, and both the $+$ and $-$ phase are in contact with it
(provided $\bdf<\widetilde\bdf$). The complete wetting regime cannot be
distinguished from the partial wetting regime in this setting, 
but see the remark
after Theorem~\ref{thm_sessile2D} for a discussion of this issue.
\subsection{2D Ising model}
Let $r\in\bbR$. The aim of this subsection is to describe the typical
configurations under the measure
\begin{equation*}
\Isbd{+}{N,rN}\bigl(\,\cdot\,\big|\,M_N = m\,\abs{\BoxxrN}\bigr)\,,
\end{equation*}
where $m\in (-m^*, m^*)$ and $M_N = \sum_{i\in\BoxxrN}\gs_i$; we will
simplify the notations further by writing simply $\Isbd{+}{N}$ ($r$ being kept
fixed).
As in Part~\ref{part_strongWulff}, it is possible to obtain precise asymptotics
for the large deviations, in the form of the following generalization of the
first part of Theorem~\ref{thm DKS}. Let $\intSTbdmin(m)$ be the infimum of the
functional $\intSTbd$ on subsets of $\rBoxTwoD$ with volume $\frac{m^* - m}{2m^*
}\abs{\rBoxTwoD}$.
\begin{thm}
\label{thm_largedevbd}
Let the inverse temperature $\beta >\beta_c$ and the boundary magnetic field
$\bdf\in\bbR$ be fixed; 
let the sequence $\{a_N\}$;  $-m*\abs{\BoxxrN}+a_N\in\text{\rm Range}(M_N )$, be
such that the limit
$$
a~=~\lim_{N\to\infty}\frac{a_N}{\abs{\BoxxrN}}~\in~(0,2m^* (\beta ))
$$
exists. Then,
$$
\log\Isbd{+}{N}\bigl(M_N = m^* \abs{\BoxxrN} - a_N \bigr) = 
-\intSTbdmin\, (1 +
\bigO(N^{-1/2}\log N)).
$$
\end{thm}
A version of this Theorem, in an approximate canonical ensemble (as
in~\eqref{eq_approx_can}), has been proven in~\cite{PfisterVelenik97}; this
stronger version can be obtained by combining the techniques
of~\cite{PfisterVelenik97} and of~\cite{IS}, see
Section~\ref{sec_tools}.

In Theorem~\ref{thm_largedevbd}, we have made no statement about the asymptotic
description of the typical configurations under the conditioned measure. The
reason is the following: These strong concentration results require the
knowledge of stability properties of the variational problem in the form, for
example, of Bonnesen inequality. However, in the present case, one does not
always have that much information about the variational problem; in fact, even
its solution is not always known. This prevents us from translating the energy
estimates on the skeletons (see~\eqref{eq_energybd}, \eqref{eq_energybd2} and
\eqref{eq_energybd3}) into strong concentration properties of
the microscopic contours. Of course, in the situations when such stability
properties are known (\cite{KoteckyPfister94} contains a simple derivation of
such a result for many situations), it is possible to obtain statements of the
same kind as those of Part~\ref{part_strongWulff}.

This illustrates the fact that although the probabilistic theory in the 2D case
is complete, in the sense that all the relevant information on the {\em
microscopic scale} is available, the sharpness of the statements one can make
on the macroscopic scale still depends on {\em macroscopic} stability
properties, which are logically separated from the probabilistic aspect of the
analysis.

However, even without information about the stability properties of the
variational problem, it is still possible to derive weak concentration
properties, in a $\bbL_1$ setting close to the one of Part~\ref{part_weakWulff}. We
present such a result in the way it is stated in~\cite{PfisterVelenik97}. In
this paper, an approximate canonical ensemble was considered, i.e. the measure
was $\Isbd{+}{N}(\,\cdot\,|\,\cA(m;c))$, where
\begin{equation}
\label{eq_approx_can}
\cA(m;c) = \bigsetof{\gs}{\bigabs{\abs{\BoxxrN}^{-1}M_N(\gs) - m}
\leq N^{-c}}\,,
\end{equation}
with $-m^* < m < m^*$, and $c$ is some real number not too large (see
Theorem \ref{thm_sessile2D} below).
We are going to prove that the phases concentrate near macroscopic droplets
which belong to the set
$\cD(m)$ 
\begin{equation*}
\cD(m)=\bigsetof{ V\subset
\rBoxTwoD}{\abs{ V}=\frac{m^* - m}{2m^* }\abs{\rBoxTwoD}\,,\,
\intSTbd(\bnd V)=\intSTbdmin(m)}\,,
\end{equation*}
Recall that to each $V\in\cD(m)$, we associate the
function $\1_V = 1_{V^c} - 1_V$.\\

To state this phase segregation Theorem, we use the mesoscopic 
notation introduced
in Part~\ref{part_weakWulff}.
Recall that $N = 2^n$. For any $a < 1$, we define a magnetization profile 
$\cM_{[an]} (\gs,x)$ at the $2^{[an]}$-scale which is piecewise constant
on boxes $\sBox{n - [an]}(x) $ with $x \in \sTor{n-[an]}$,
\begin{equation}
\cM_{[an]}(\gs,x) = 2^{- d [an]} \sum_{i \in  \dBox{[an]}(2^n x)}
\gs_i \, .
\end{equation}
We get
\begin{thm}\label{thm_sessile2D}\cite{PfisterVelenik97}
Let $\gb>\gbc$,  $\bdf\in\bbR$, $-m^*<m<m^*$ and $1/4>c>0$. Then
there exist a function $\gd(N)$ such that $\lim_{N\ra\infty}\gd(N)=0$, a real
number $\gk>0$ and a coarse-graining parameter $1>a>0$ such that for $N$ large
enough
\begin{equation*}
\Isbd{-}{N}\bigl(\frac{\cM_{[an]}}{m^*}\in\union_{V\in\cD(m)} \cV(\1_V,
\gd(N)) \, \big| \ \cA(m;c)\,\bigr)\geq
1-\exp\{-\bigO(N^{\gk})\}\,.
\end{equation*}
\end{thm}

\bigskip
\noindent
{\bf Remark:} In this case, it should also be possible to study the whole range
of moderate deviations, combining the techniques of \cite{IS} and
\cite{PfisterVelenik97}, although this has not been done explicitly. We
briefly describe the results obtained for large deviations sufficiently close to
volume order \cite{Velenik97}.

As long as $\bdf>-\hw(\gb)$, the results are similar to those obtained in the
setting of Part~\ref{part_strongWulff}: The measure concentrates on
configurations containing a single large droplet of $-$ phase, with Wulff or
Winterbottom shape depending on $\bdf$; in particular, the order of the large
moderate deviations is still $\exp\{-\bigO(\sqrt{a_N})\}$. There should not
be any problem to extend this to the whole large deviations regime ($a_N\gg
N^{4/3}$).

More interesting is the case $\bdf\leq-\hw(\gb)$. For those values of the
boundary field, the system is in the complete wetting regime
($\taubd(\gb,\bdf)=-\tau^*_\gb$), and the solution of the unconstrained
variational problem is degenerate. The solution of the constrained variational
problem in $\rBoxTwoD$ is however still well-defined for every $N$; it is
obtained by extracting the cap of a Wulff shape and rescaling it so that the
basis of the cap completely covers the wall and the rescaled cap has the
required volume. When $N$ goes to infinity, this droplet spreads out to become
a thin film in the limit (covering the entire wall, hence the terminology
complete wetting), and the corresponding value of the surface free energy
functional goes to zero. As a result of this, the scale of the large moderate
deviations {\em is not the same as when $\bdf<\hw(\gb)$}; indeed the leading
term of the asymptotics can again be computed explicitly, and is found to be
of order $\exp\{-\bigO((a_N)^2\,N^{-3})\}$. In particular, we see that the large
moderate deviations cannot extend up to $a_N\sim N^{4/3}$, since $(a_N)^2\,
N^{-3}$
is of order $1$ already when $a_N\sim N^{3/2}$. This should not be surprising
since, in the complete wetting regime, the volume under the microscopic contour
is expected to have typical fluctuations of order $N^{3/2}$ (this can be shown
when $\bdf=-1$ and $\gb$ is very large using the convergence to Brownian
excursion stated in \cite{Dobrushin93}). Therefore, typical fluctuations of
magnetization in the complete wetting regime are not governed by bulk
fluctuations anymore, but by fluctuations of the microscopic phase separation
line. To prove that this behavior is valid up to $a_N\sim N^{3/2}$ might
be a non-trivial task.\qed
\subsection{Ising model in $D\geq 3$}
Let $r\in\bbR$ and let $\cD(m)$ be the set of macroscopic droplets at
equilibrium in $\rBox$,
\begin{equation*}
\cD(m)=\bigsetof{ V\subset
\rBox}{\abs{V}=\frac{m^*-m}{2m^*}\abs{\rBox}\,,\,
\intSTbd(\bnd V)=\intSTbdmin(m)}\,.
\end{equation*}
The rest of the notations were introduced in
Part \ref{part_weakWulff}. The main result is the following
\begin{thm}\cite{BodineauIoffeVelenik99}
\label{thm_sessile3D}
For any $\gb$ in $\frB_p$, any $\bdf\in\bbR$, any $m$ in $(-m^*, m^* )$, the
following holds: For any $\gd >0$, there is $k_0 = k_0 (\gd)$ such that for 
$\nu < \frac{1}{d}$
\begin{eqnarray*}
\lim_{N \to \infty} \; 
\min_{k_0 \leq k \leq \nu n} \; 
\Isbd{+}{N} \Biggl(
\frac{\cM_k}{m^*} \in \union_{V\in\cD(m)} \cV(\1_V,\gd)
\  \Big| \quad  M_{N} \leq m\,\abs{\BoxxrN} \Biggr) = 1\,.
\end{eqnarray*}
\end{thm}
\section{The tools}
\label{sec_tools}
\setcounter{equation}{0}
In this Section, we explain how the procedures described in Parts
\ref{part_weakWulff} and \ref{part_strongWulff} have to be modified 
to take into account the effect of the boundary.
\subsection{2D Ising model}
\label{ssec_2D}
We describe the main modifications one needs to apply to the proofs of Part
\ref{part_strongWulff} in order to get the results stated in
Theorems~\ref{thm_largedevbd} and  \ref{thm_sessile2D}. We split this
Subsection into two parts, one dealing with the lower bound on
$\Isbd{-}{N}(\cA(m;c))$ or $\Isbd{-}{N}(M_N=-m^*\abs{\BoxxrN}+a_N)$, the other
one with the upper bound.

\bigskip\noindent
\underline {The lower bound.}
The constrained variational problem is more difficult than the usual one. In
fact, as noted above, the solution (and {\it a fortiori} its stability) is not
known in general, although it is in many cases. This prevents us from
proceeding as in Part \ref{part_strongWulff}, where the lower bound
follows from summing over large contours fluctuating around the Wulff shape. It
would then appear necessary to make the same kind of proof, but for any
configurations of droplets surrounding the right volume (all potential
solutions to the variational problem). This, however, would be tricky; indeed,
since we want our results to hold for large, but {\em finite} boxes, it is
compulsory to obtain estimates {\em uniform} over the droplet in the chosen set!
Fortunately, properties of the surface tension and wall free energy allow us to
restrict our analysis to a small class of well-behaved droplets: The solution of
the variational problem is necessarily taken on a {\em single convex} droplet.
This is a consequence of the convexity of $\tau_\gb$ (use Jensen inequality)
and the fact that $\taubd(\gb,\bdf)\leq\tau_\gb^*$, which imply that replacing a
droplet by its convex hull cannot increase the surface free energy; rescaling
the resulting droplet decreases the energy even more. 
It is thus enough to prove the following
\begin{pro}\cite{PfisterVelenik97}
\label{pro_lowersessile2D}
Let $\gb>\gbc$ and $\bdf\in\bbR$. There exists $N_0=N_0(\gb,\bdf,m,c,r)$ and a
constant $C$ such that, for any simple closed rectifiable curve $\cC$ which is
the boundary of a convex body of volume $\abs{\rBoxTwoD}
(m^*(\gb)+m)/2m^*(\gb)$ contained in $\rBoxTwoD$, and for all $N\geq N_0$,
\begin{equation*}
\Isbd{-}{N}(\cA(m;c)) \geq \exp\{-\intSTbd(\cC)\; N -
\gb\,C\,N^{1/2}\log N\}\,.
\end{equation*}
A completely analogous statement holds in the case of the exact canonical
ensemble.
\end{pro}
The proof of Proposition \ref{pro_lowersessile2D} is similar to the proof of
Theorem \ref{lb}. We sketch now the main changes needed to deal with
the boundary conditions. The case $\bdf\leq 0$ requires a slightly more
complicated proof than the case $\bdf>0$ so we first consider the latter.

\medskip\noindent
{\it First case: $\bdf>0$}

\begin{figure}[t]
\centerline{\psfig{file=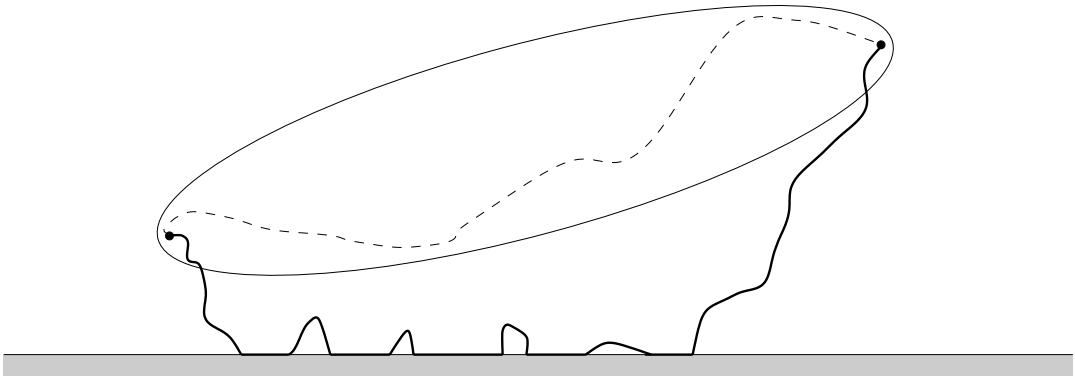,height=5cm}}
\caption{When $\taubd(\gb,\bdf)<\tau^*_\gb$, the open contour connecting two
sites close enough to the wall might not stay inside an elliptical set as in
the bulk (dashed contour), but instead might get pinned by the wall (full
contour). In such a case, the exponential decay-rate is in general not given by
$\tau_\gb$ or $\taubd(\gb,\bdf)$.} 
\label{fig_pinning}
\end{figure} 

As in the usual case, we want to approximate $\cC$ with some polygonal curve
with vertices on the dual lattice, and then sum over all contours going through
the latter; this would allow us to extract, for each piece of the contour, the
surface tension of the corresponding part of the polygonal line. Here, however,
we want to be able to extract the wall free energy when the curve $\cC$ follows
the wall. There are some complications related to this: If two vertices are
close to the wall, but don't belong to it\footnote{Consider, for example, a
family of curves $\cC$ getting closer and closer to the wall; since we need
estimates uniform in all such curves, one has to be able to deal with such a
situation.}, the sum over the corresponding piece of contour might not yield
simply $\tau_\gb$ or $\taubd(\gb,\bdf)$, but some complicated mixture, since typical such
contours might first go down to the wall, then follow it on some length, and
only then go up to the other vertex, see Fig.~\ref{fig_pinning}; this kind of
behavior has been studied in details in \cite{PfisterVelenik98}. It turns out
that it is possible to construct a polygonal approximation to the curve $\cC$
whose surface tension is not too large in comparison with that of $\cC$, while
removing these possible pathologies.

\begin{figure}[t]
\centerline{
\psfig{file=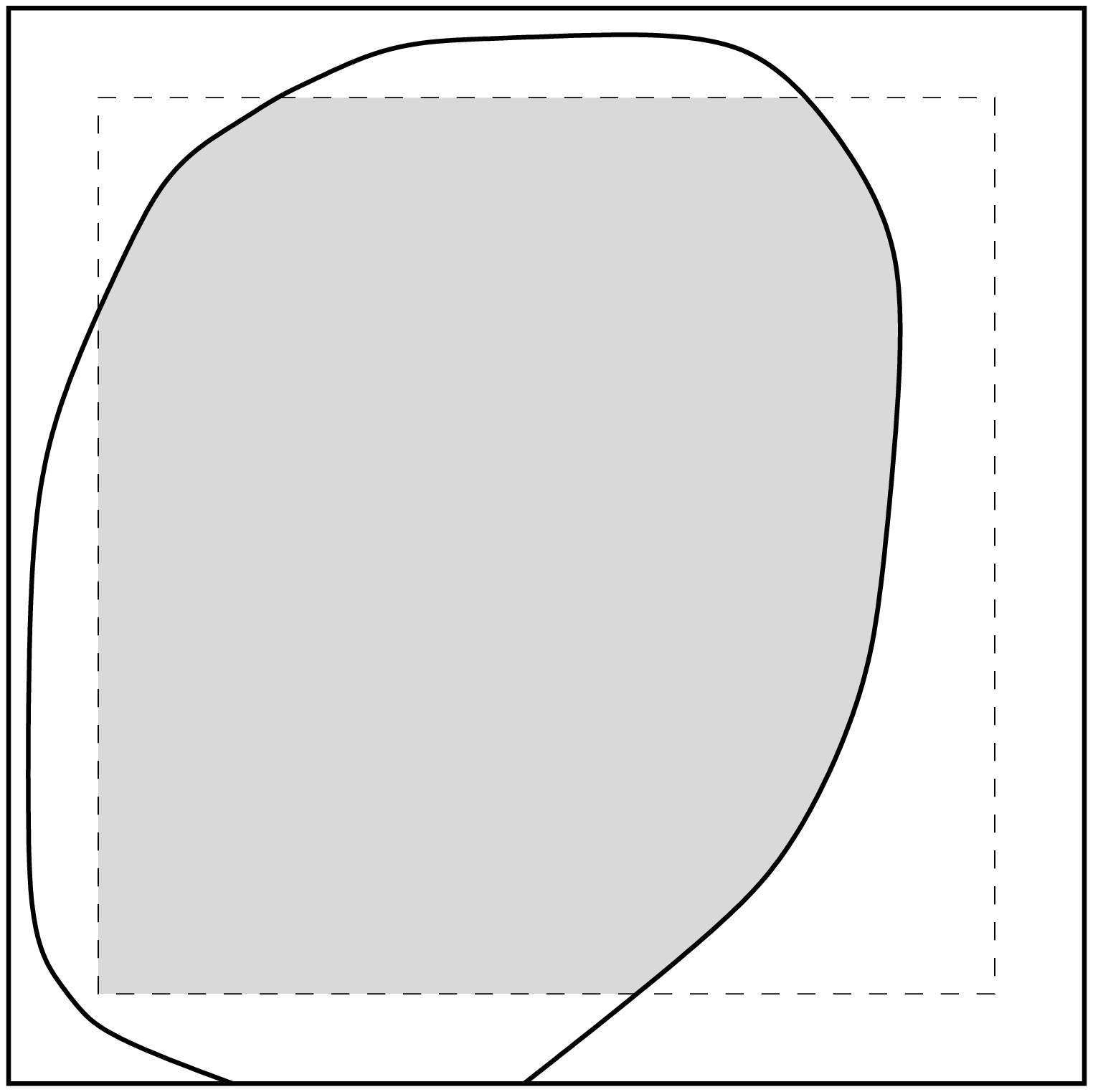,height=6cm}\hspace{1cm}
\psfig{file=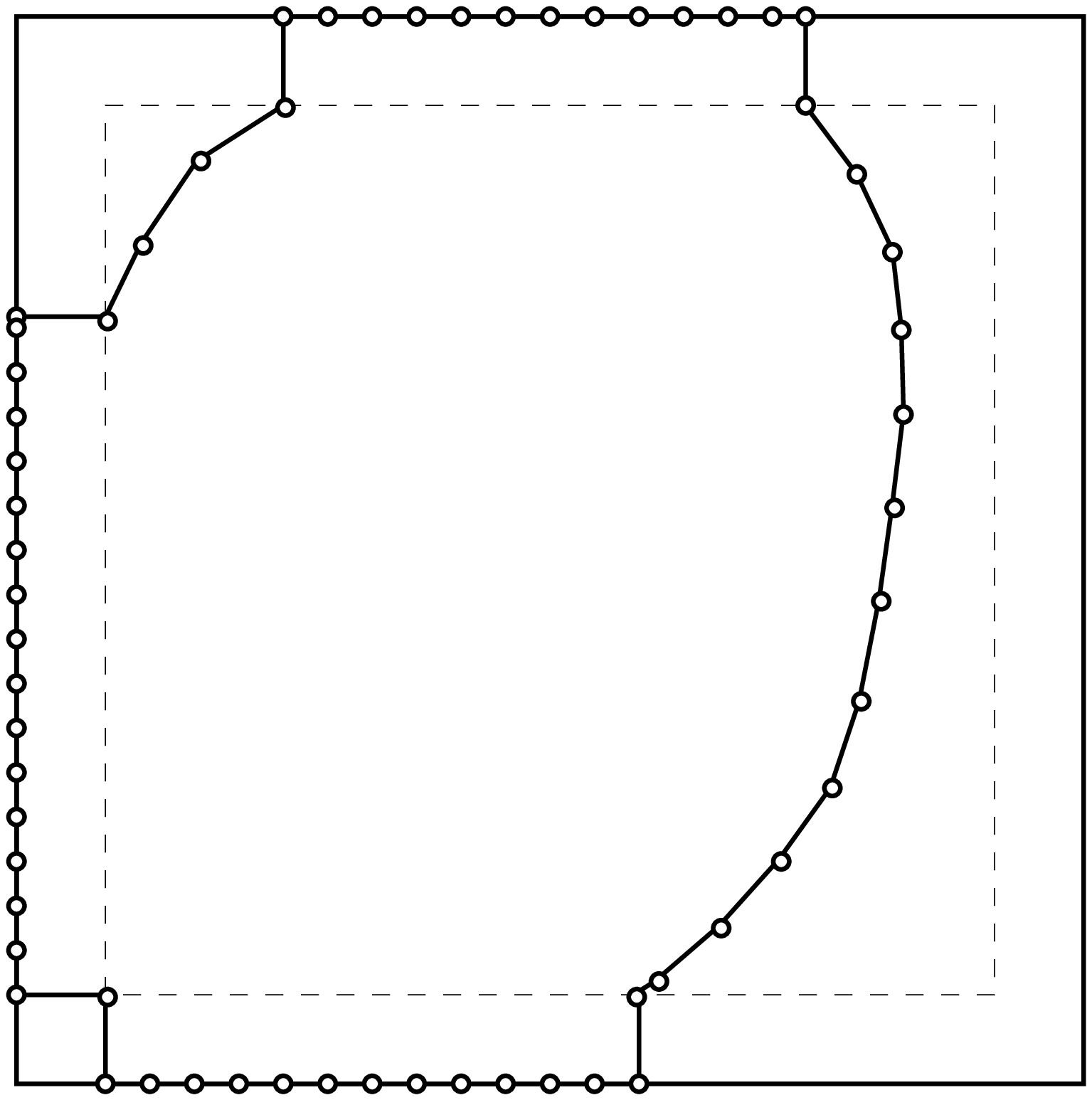,height=6cm}
}
\caption{Left: The curve $\cC$; the shaded area represents the convex body
whose boundary is $\cC_N$ and the dashed line is the boundary of
$\rBoxTwoD(N)$. Right: The polygonal approximation $\widehat\cP_N$, the dots
representing its vertices.} 
\label{fig_bd_lb_cg}
\end{figure}

The idea is the following. Let $\gd_N = N^{-1/2}\log N$, and set
\begin{equation*}
\rBoxTwoD(N) = \setof{x\in\rBoxTwoD}{\min_{y\not\in\rBoxTwoD}\,\normI{y-x} >
\gd_N}\,.
\end{equation*}
Let $V$ be the convex body with boundary $\cC$ and set $\cC_N =
\bnd(V\cap\rBoxTwoD(N))$ We first construct a polygonal approximation for each
of the components of $\cC_N\cap\rBoxTwoD(N)$ with segments of length $\gd_N$
(apart from at most 8 of them which may be shorter). Set $[x,y] =
\setof{z\in\cC_N}{z(2)=\gd_N}$. If $[x,y]\neq\eset$, we connect the two
corresponding pieces of polygonal lines by a broken line from $x$ to
$(x(1),0)$, then to $(y(1),0)$, and finally to $y$; we divide the segment
between $(x(1),0)$ and $(y(1),0)$ into segments of length $\gd_N/2$ (except
possibly for the last one which can be shorter). We repeat this construction
for the three other sides of the box. The resulting closed polygonal line is
denoted by $\widehat\cP_N$ (see Fig.~\ref{fig_bd_lb_cg}). Notice that by
construction there exists an absolute constant $C$ such that
\begin{align*}
\intSTbd(\cC) \geq \intSTbd(\widehat\cP_N)& - C\gb\gd_N\,,\\
\abs{\vol(\cC)-\vol\widehat\cP_N} \leq& C\,\abs{\rBoxTwoD}\, \gd_N\,.
\end{align*}
We then rescale the polygonal line $\widehat\cP_N$ by a factor $N$ and if necessary
move slightly the rescaled vertices so that they belong to the dual lattice;
the rescaled polygons is denoted by $\cP_N$. We then define a class $\frG$ of
closed contours going through the vertices of $\cP_N$ (in the right order),
and staying in some small boxes along its edges. For all edges of length
smaller than $N\gd_N$, as well as for the (up to 8) pieces we added above to
join $\cC_N$ to the boundary, we impose that the corresponding piece of the
contour is a fixed length-minimizing path between the vertices.

The rest of the argument proceeds in a similar way as in the standard case. The
estimates in the phase of small contours carry over without any problems since
in that case the effect of the boundary field cannot propagate far away from the
wall.

We still have to explain how one can extract the correct surface tension for
$\widehat\cP_N$ from the sum over contours in the class $\frG$ introduced
above. To do this, we use several results about the {\em random-line
representation}, proved in \cite{PfisterVelenik97,PfisterVelenik98}. To lighten
the notation, we simply write $\weightB{\bdf}$ instead of
$q_{\dBoxxrN}^{\gb^*,\bdf^*}$; $\gb^*$ and $\bdf^*$ are the dual of $\gb$ and
$\bdf$, see \eqref{eq_dualitygeneral}. The first inequality is just the
analogue of \eqref{3.4.qasuper} in our case, which turns out to be valid for
arbitrary ferromagnetic coupling constants: The weight of any high-temperature
contour $\gga\in\frG$ satisfies (\cite{PfisterVelenik97}, Lemma~5.4)
\begin{equation*}
\weightB{\bdf}(\gga) \geq \prod \weightB{\bdf}(\gga_k)
\end{equation*}
where $\gga_k$ denotes the piece of the contour $\gga$ between the $k$th and
$k+1$th vertices of $\cP_N$. The next step is to replace $\weightB{\bdf}(\gga_k)$ by
the corresponding infinite-volume quantity. First, for any $\gga_k$ joining
vertices not belonging to $\thedwall \df \setof{i\in\BoxxrN^*}{i(2)=-\tfrac12}$
(note that $\gga_k$ stays necessarily at a distance $\bigO(N\gd_N)$ from
$\thedwall$)
\begin{equation*}
\weightB{\bdf}(\gga_k) \geq (1-e^{-\bigO(N\gd_N)})\;\weight(\gga_k)\,;
\end{equation*}
second, for the pieces $\gga_k$ joining two sites of $\thedwall$, we
use
\begin{equation*}
\weightB{\bdf}(\gga_k) \geq \weightL{\bdf}(\gga_k)\,,
\end{equation*}
where $\dhalfspace \df \setof{i\in\bbZ^d_\star}{i(2)\geq -\tfrac12}$ (both
results are proved in \cite{PfisterVelenik97}, Lemma~5.3).
Finally, the remaining pieces have a length at most $8N\gd_N$, so that their
total weight is larger than $e^{-C\bigO(N\gd_N)}$.

The last step is to extract the surface free energy. The basic tool to do this
is, as in the proof of Theorem \ref{skeletonlb}, concentration properties
for open contours between 2 fixed dual sites. For the pieces $\gga_k$ not
touching the boundary, we can use the usual infinite volume results
based on \eqref{eq_concentration}, setting $s=N\gd_N$. For the pieces along the
boundary, one can use the following statement (\cite{PfisterVelenik98},
Lemma~6.10):
\begin{equation}\label{eq_concwall}
\sumtwo{\gl:\,i\to j}{\gl\subset{\bf N}_K (i,j)\cap\dhalfspace}
\weightL{\bdf}(\gl) \geq
\bk{\gs_i\gs_j}^{\gb^*,\bdf^*}_{\dhalfspace}\;\left( 1+ \smallo(1)\right) ,
\end{equation}
where ${\bf N}_K (i,j)$ is defined in Appendix~B (with $s=N\gd_N$).
(In fact, \eqref{eq_concwall} can be strengthened when $\bdf<\hw(\gb)$: in this
case, the set ${\bf N}_K (i,j)\cap\dhalfspace$ can be replaced by the set
(\cite{PfisterVelenik98}, Lemma~6.13)
$$
\setof{k\in\dhalfspace}{(i(1)\wedge j(1))-K\log \gd_N \leq k(1)\leq (i(1)\vee
j(1))+K\log \gd_N,\, k(2) \leq K\log \gd_N}\,,
$$
which is compatible with our picture of partial wetting.)

The result then follows from lower bounds on the corresponding 2-point
functions. The only new inputs are the following lower bounds on the boundary
2-point function,
\begin{align}
\bk{\gs_i\gs_j}^{\gb^*,\bdf^*}_{\halfspace^*} &\geq
C\,\frac{\exp\{-\taubd(\gb,\bdf)\norm{j-i}\}} {\norm{j-i}^{3/2}} \quad\quad
&\forall\bdf\geq\hw(\gb)\,, \label{eq_lower2ptfwall1}\\
\bk{\gs_i\gs_j}^{\gb^*,\bdf^*}_{\halfspace^*} &\geq
C\,\exp\{-\taubd(\gb,\bdf)\norm{j-i}\} &\forall\bdf<\hw(\gb)\,, \label{eq_lower2ptfwall2}
\end{align}
for any $i,j\in\thedplane\df\setof{k\in\dhalfspace}{k(2)=-\tfrac12}$.
\eqref{eq_lower2ptfwall2} is proved in \cite{PfisterVelenik97}, Prop.~7.1, while
\eqref{eq_lower2ptfwall1} follows from exact computations in the case
$\bdf^*=1$ \cite{McCoyWu73}, and \cite{PfisterVelenik97}, Prop.~7.1,
\begin{equation*}
\bk{\gs_i\gs_j}^{\gb^*,\bdf^*}_{\halfspace^*} \geq (\tanh\gb^*)^2\,
\bk{\gs_i\gs_j}^{\gb^*,1}_{\halfspace^*},\quad\quad\forall\bdf\geq 0\,.
\end{equation*}

\medskip\noindent
{\it Second case: $\bdf=0$}

This is a somewhat marginal case. The apparent difficulty is that in this case
$\bdf^*=\infty$. However, this does not create any real complications. One just
has to modify the construction of the first case as follows: We replace the
polygonal line $\widehat\cP_N$ by the (possibly open) polygonal line
$\widehat\cP_N\setminus\setof{u\in\bbR^2}{u(2)=0}$; we then sum over contours
going through the vertices of this polygonal line (contours which are open if
the polygonal line is open). This does not give any contribution for the part
of $\cC$ along the wall, which is what we want since
$\taubd(\gb,0)=0$.

\medskip\noindent
{\it Third case: $\bdf<0$}

\begin{figure}[t]
\centerline{
\psfig{file=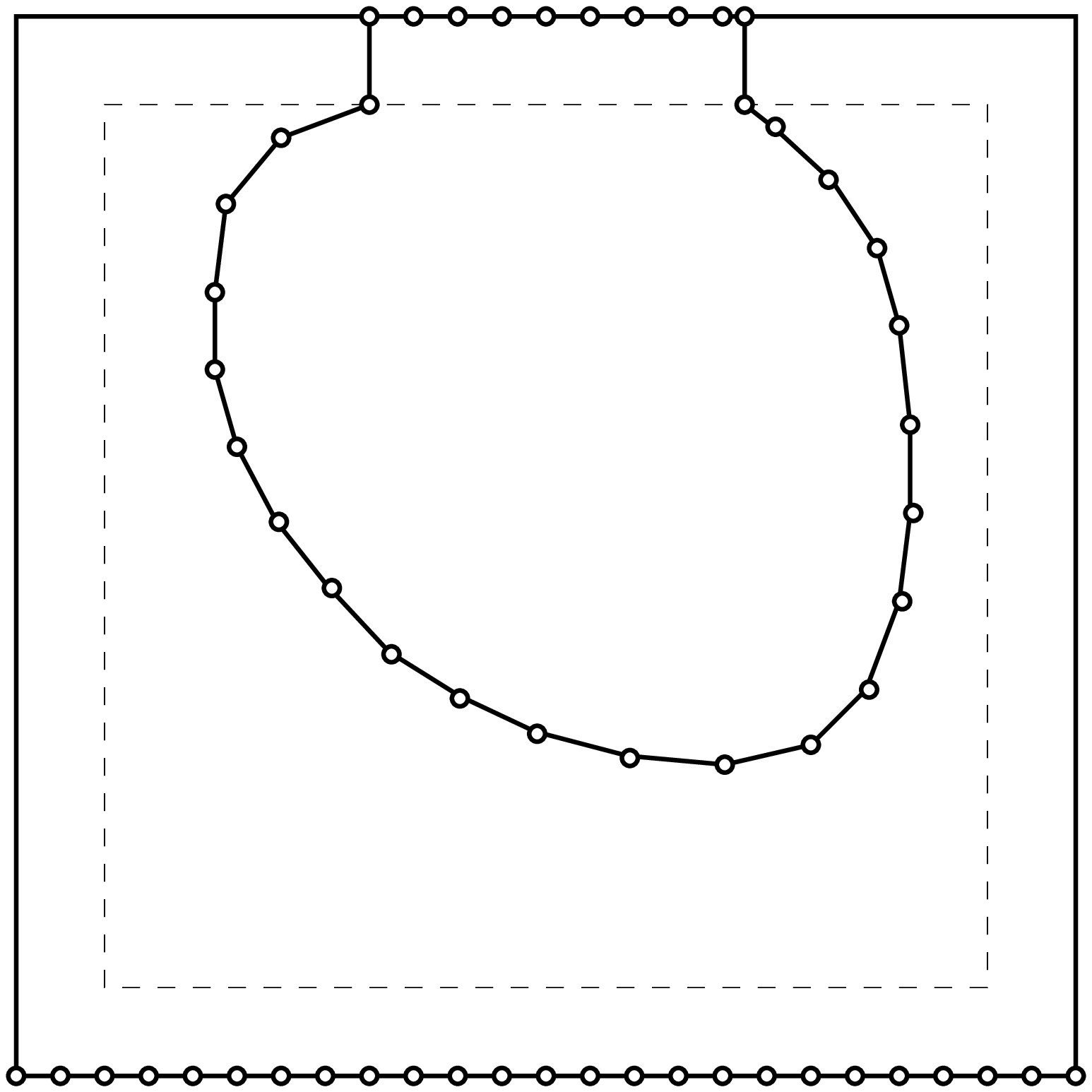,height=6cm}\hspace{1cm}
\psfig{file=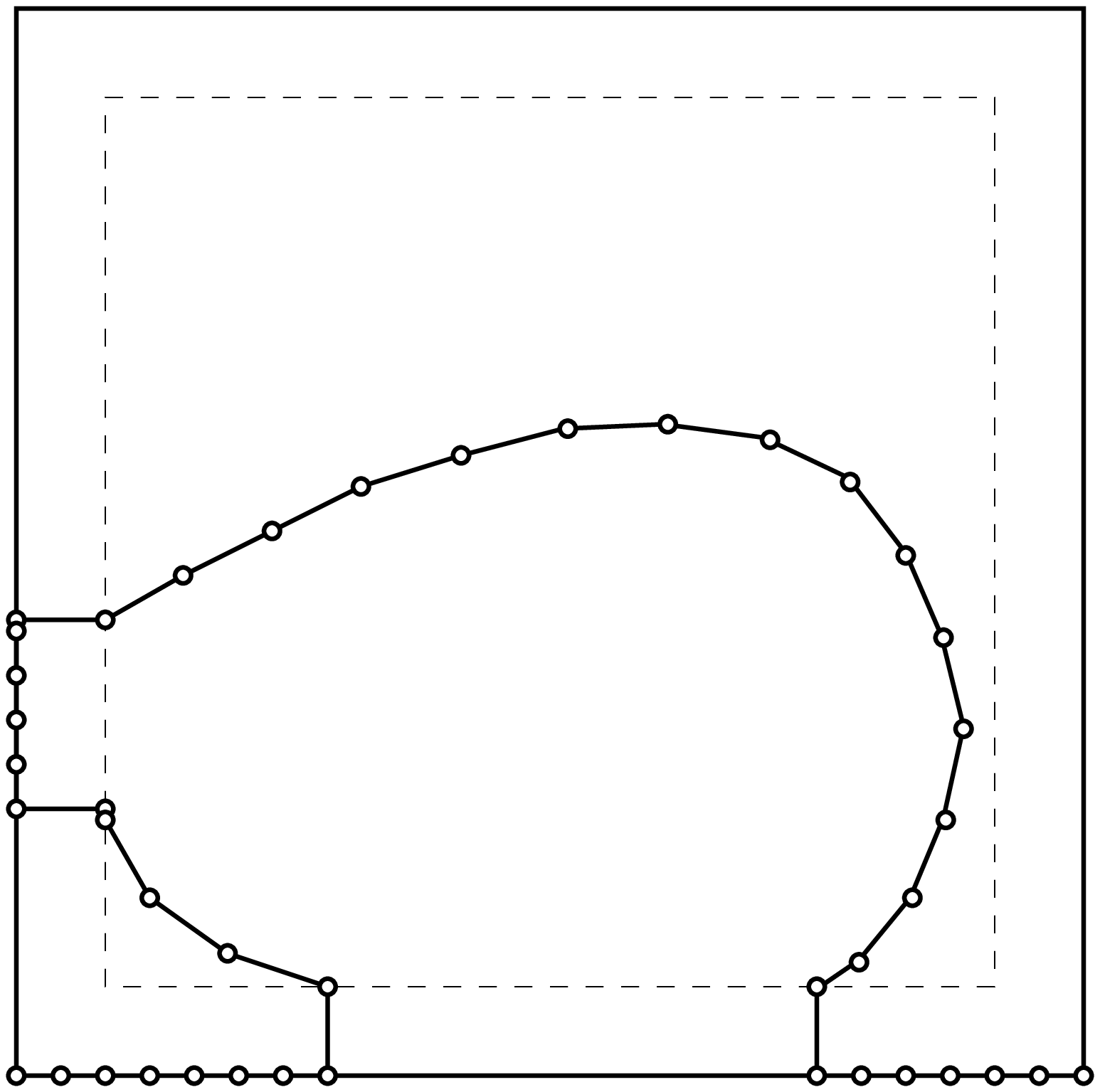,height=6cm}
}
\caption{The construction for $\bdf<0$. Left: $I=\eset$ (two polygonal lines:
one open and one closed. Right: $I\neq\eset$ (one open polygonal line).} 
\label{fig_bd_lb_cg_negh}
\end{figure} 
This is slightly more tricky. In this situation, one may be even more
pessimistic, since the duality is simply not defined when non-ferromagnetic
interactions are present! However, this turns out to be a false problem. Indeed,
we can use the following obvious identity to recover ferromagnetic interactions
(see footnote \ref{footnote}, p.~\pageref{footnote}),
\begin{equation*}
\Isbdf{+}{N}{\bdf} = \Isbdf{\pm}{N}{\abs\bdf}\,,
\end{equation*}
where $\pm$ correspond to the boundary condition $\overline{\gs}_i = 1$ if
$i(2)\geq 0$ and $\overline{\gs}_i = -1$ otherwise.

We then construct $\widehat\cP_N$ as in the first step and set
$I=\widehat\cP_N\cap \setof{x\in\rBoxTwoD}{x(2)=0}$. If $I=\eset$, then we
subdivide the set $\setof{x\in\rBoxTwoD}{x(2)=0}$ into segments of length
$\gd_N/2$ (except possibly for the last one, which might be shorter); this
defines a second (open) polygonal line $\widehat\cP_N'$ (with all its vertices
along the wall) (see Fig.~\ref{fig_bd_lb_cg_negh}). We then introduce a class
of {\em pair} of contours $(\gga,\gga')$, $\gga$ going through the vertices of
$\cP_N$ and defined as before, and $\gga'$ following the wall, going through
the vertices of $\cP_N'$ and staying inside small boxes along its edges,
similarly as for the other one ($\gga'$ is open). By construction $\gga$ and
$\gga'$ are disjoint. Duality then implies the following identity
\begin{align}
\Isbdf{\pm}{N}{\abs{\bdf}}(\{\gga,\gga'\}\subset\boldsymbol\gga(\,\cdot\,))
&= (\PFbdf{\pm}{N}{\abs{\bdf}})^{-1}\;
w(\gga)w(\gga')\;\sumtwo{\underline\gz:}{(\underline\gz,\gga,\gga')\text{
$\gL^*$-comp.}} w(\underline\gz) \nonumber\\
&= (1-e^{-\bigO(N)})\;\frac {\PFbdf{+}{N}{\abs{\bdf}}}
{\PFbdf{\pm}{N}{\abs{\bdf}}} \;\weightB{\abs{\bdf}}(\gga,\gga')\,.
\label{eq_dualityhneg}
\end{align}
The factor $(1-e^{-\bigO(N)})$ comes from the fact that we can apply duality
only to simply connected sets, and the exterior of $\gga$ is {\em not} simply
connected. We must therefore forbid families $\underline\gz$ for which duality
does not hold; since such families must contain at least one contour
surrounding $\gga$, we get the above correction.

We can now proceed as in the first case. The only additional work to do is to
analyze the ratio of partition functions in \eqref{eq_dualityhneg}, but this is
easy, since by duality
\begin{equation}\label{eq_ratio}
\frac {\PFbdf{+}{N}{\abs{\bdf}}} {\PFbdf{\pm}{N}{\abs{\bdf}}} =
\bigl( \bk{\gs_{t_{\rm l}}\gs_{t_{\rm r}}}^{\gb^*,\abs{\bdf}^*}_{\dBoxxrN}
\bigr)^{-1} \geq e^{\taubd(\gb,\abs{\bdf})\,(2N+1)}\,,
\end{equation}
where $t_{\rm l}=(-L-\tfrac12,-\tfrac12)$ and $t_{\rm r}=(L+\tfrac12,-\tfrac12)$
are the two dual sites at the lower left and lower right corners of $\dBoxxrN$,
and the last inequality follows from the upper bound (see
\cite{PfisterVelenik97} for example)
\begin{equation}\label{eq_upbd2ptfbd}
\bk{\gs_i\gs_j}^{\gb^*,\abs{\bdf}^*}_\dBoxxrN \leq
e^{-\taubd(\gb,\abs{\bdf})\norm{j-i}}\,,
\end{equation}
valid for any $i,j\in\thedwall$. We then see that the ratio of partition
function cancels the contribution from the sum over the open contour $\gga'$,
up to an error term $\exp\{\cO(N\gd_N)\}$.

\medskip
If $I\neq\eset$, the situation is simpler. Let's write $I=[x,y]$; then we
define a new polygonal line $\widehat\cP_N^\pm$: $\widehat\cP_N^\pm$ goes from
the lower right corner of $\rBoxTwoD$ to $a$ along the wall, then it follows
$\widehat\cP_N\setminus \setof{x\in\rBoxTwoD}{x(2)=0}$ up to $b$ and finally
goes from $b$ to the lower right corner of $\rBoxTwoD$ (see
Fig.~\ref{fig_bd_lb_cg_negh}). We subdivide as usual the part of
$\widehat\cP_N^\pm$ along the wall into segments of length $\gd_N/2$ and proceed
as in the first case, with $\widehat\cP_N^\pm$ replacing $\widehat\cP_N$, using
\eqref{eq_dualityhneg}. Summing over the open contour going through the
vertices of $\cP_N^\pm$ produces (up to the usual error term) a term
$\exp\{-\intSTbdf{\abs{\bdf}}(\widehat\cP_N^\pm)\,N\}$. Combining this with
\eqref{eq_ratio} and observing that
\begin{equation*}
\exp\{2\taubd(\gb,\abs{\bdf})\,N\}\;
\exp\{-\intSTbdf{\abs{\bdf}}(\widehat\cP_N^\pm)\, N\} =
\exp\{-\intSTbdf{\bdf}(\widehat\cP_N)\, N\}\,,
\end{equation*}
the conclusion follows as in the usual situation.

\bigskip\noindent
\underline {The upper bound.}
Let us now turn our attention to the proof of the upper bound. The basic
strategy is completely similar to that of the standard case, see
Subsection~\ref{dima_structure_ub}. The only serious modification concerns the
energy estimate, which should now associate the functional $\intSTbd$ to the
probability of skeletons. Again, the case $\bdf\geq 0$ is somewhat simpler than
the other, so we start with this one.

\medskip\noindent
{\it First case: $\bdf\geq 0$}

The basic problem we encounter when trying to make the energy estimate is the
same we met in the proof of the lower bound. Summing over an open contour
connecting two dual sites $i$ and $j$ might not yield a decay of order
$\exp\{-\tau_\gb(j-i)\}$ or $\exp\{-\taubd(\gb,\bdf)\norm{j-i}\}$ if $i$ and
$j$ are close enough to the wall but not on it (see \cite{PfisterVelenik98}).
However, the following bound, proven in \cite{PfisterVelenik97}, Lemma~5.1, is
sufficient to derive the energy estimate,
\begin{equation}\label{eq_upbddisj}
\sumtwo{\gl:\,i\ra j}{\gl\cap\cE(\thedwall)=\eset}\weightB{\bdf}(\gl) \leq
\exp\{-\tau_\gb(j-i)\}\,,
\end{equation}
for any $\bdf\geq 0$; $\cE(\thedwall)=\{e^*\subset\thedwall\}$. The definition
of skeletons will be done in such a way as to ensure that the additional
constraint $\gl\cap\cE(\thedwall)=\eset$ is automatically satisfied, see below.
We also need to extract the wall free energy when summing over contours joining
two dual sites belonging to $\thedwall$; this however is nothing else as
\eqref{eq_upbd2ptfbd}.

Let us now describe the construction of a skeleton $S=(u_1,\dots,u_n)$ of a
closed contour $\gga$. Remember that we have to define the skeletons in such a
way as to ensure that 1) the piece of the contour between two dual sites not
both on the wall must be edge-disjoint from the wall, and 2) the Hausdorff
distance between the contour $\gga$ and the polygonal line $\Pol(S)$ is smaller
than the cutoff parameter $s(N)$.

\begin{figure}[t]
\centerline{
\psfig{file=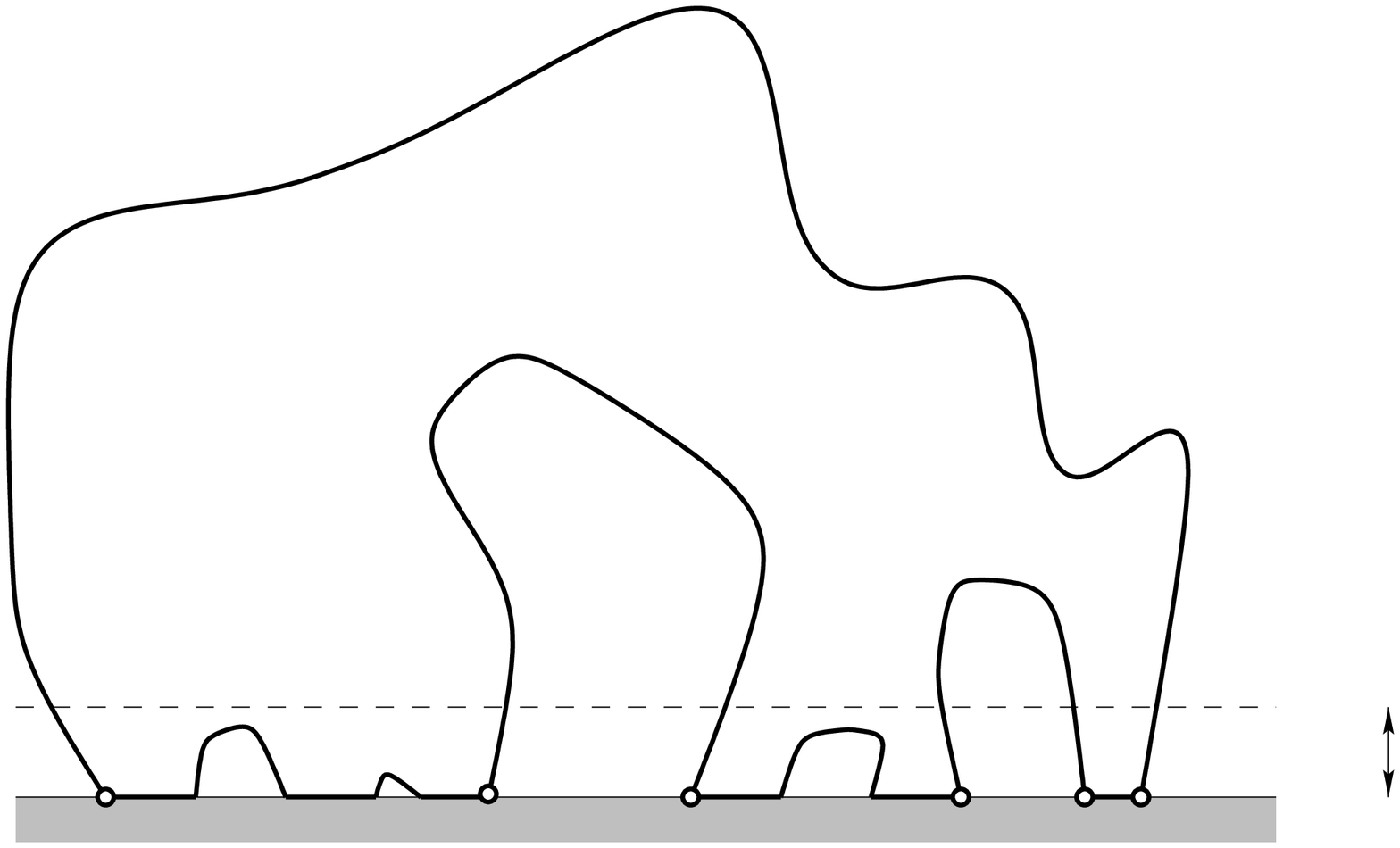,height=4.5cm}\hspace{3.7mm}
\psfig{file=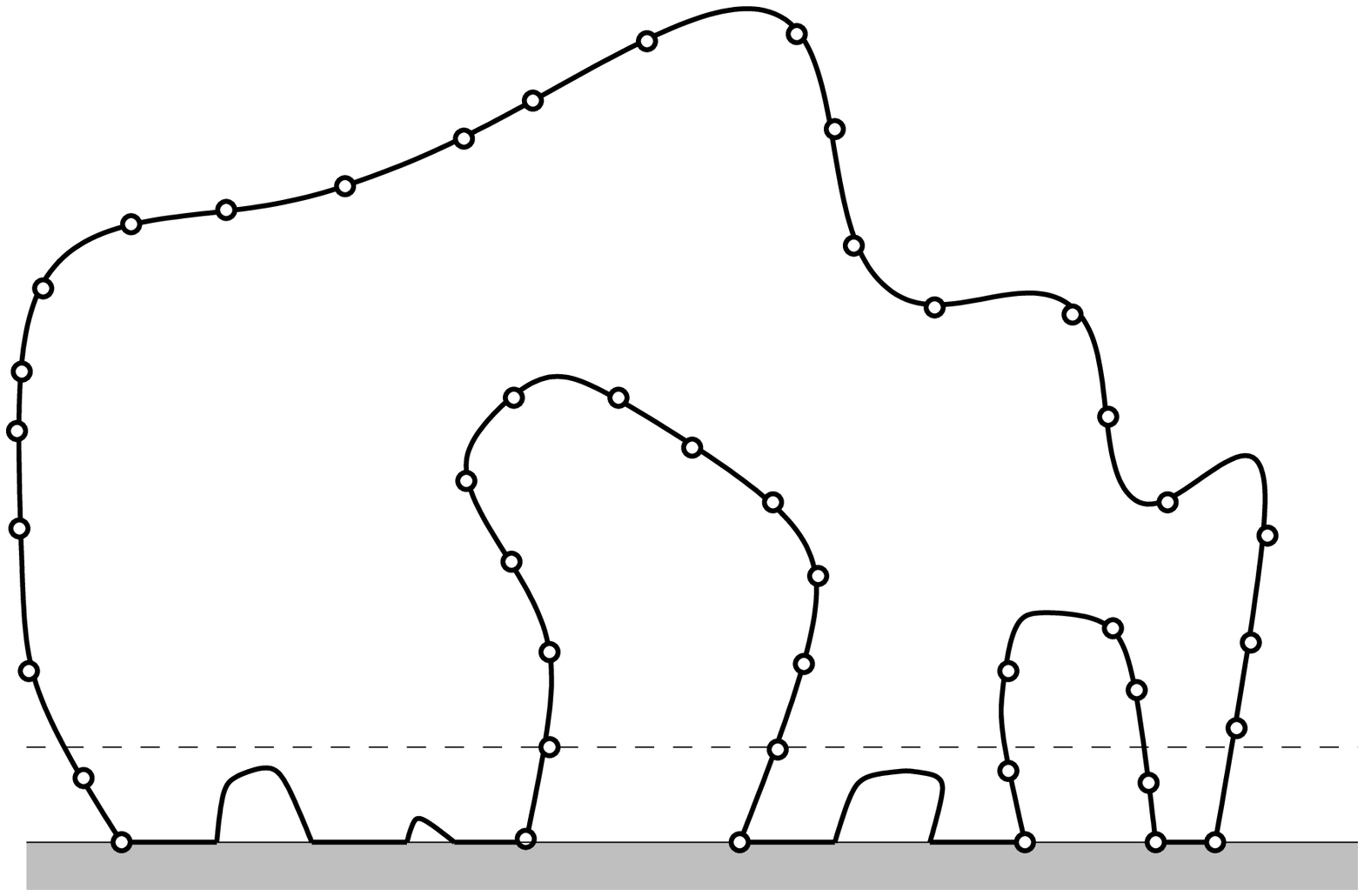,height=4.5cm}
}
\figtext{ 
\writefig        0.10   0.88    {\footnotesize $s(N)$}
\writefig       -6.77   0.83    {\footnotesize $v_1$}
\writefig       -4.70   0.83    {\footnotesize $v_2$}
\writefig       -1.22   0.83    {\footnotesize $v_{2m}$}
}
\caption{Left: A contour touching the wall and the family $(v_1,\dots,v_{2m})$.
Right: An $s$-skeleton for the contour.} 
\label{fig_bd_ub_skel}
\end{figure} 
For contours $\gga$ which do not touch the wall, the definition of skeletons is
the same as in Part~\ref{part_strongWulff}. Suppose
$\gga\cap\cE(\thedwall)\neq\eset$. Let us define
$(v_1,\dots,v_{2m})$ as the {\em minimal} family of dual sites satisfying the
following properties:
\begin{enumerate}
\item $v_k\in\thedwall\cap\gga$ for $k=1,\dots,2m$ and $v_k(1)< v_{k'}(1)$ if
$k< k'$;
\item $(v_1,\dots,v_m)$ split $\gga$ into pieces $\gga_1:v_1\ra
v_2,\dots,\gga_{2m}:v_{2m}\ra v_1$, such that
\begin{itemize}
\item $\gga_{2k}\cap\cE(\thedwall)=\eset$ for all $k=1,\dots,m$.
\item $\dH(\gga_{2k},\setof{x\in\bbR^2}{x(2)=-1/2})>s(N)$ for all $k=1,\dots,m$.
\item $\dH(\gga_{2k+1},\setof{x\in\bbR^2}{x(2)=-1/2})\leq s(N)$ for all $k=1,\dots,m$.
\end{itemize}
\end{enumerate}
We then say that $S=(u_1,\dots,u_n)$ is an $s$-skeleton of $\gga$ if
\begin{itemize}
\item All vertices of $S$ belong to $\gga$.
\item $v_1,\dots,v_{2m}$ are vertices of $S$.
\item The only vertices of $S$ along $\gga_{2k+1}$ are $v_{2k+1}$ and
$v_{2k+2}$, for all $k=1,\dots,m$.
\item The distance between any successive pair of vertices $u_l,u_{l+1}$ of $S$
along $\gga_{2k}$ satisfies $s(N)/2\leq\normsup{u_l-u_{l+1}} \leq 2s(N)$, for
all $k=1,\dots,m$. \item $\dH(\gga,\Pol(S))\leq s(N)$.
\end{itemize}
This definition has the nice property that either $u_l$ and $u_{l+1}$ both
belong to $\thedwall$, or the part of $\gga$ between these two sites is
edge-disjoint from $\thedwall$ (see Fig.~\ref{fig_bd_ub_skel}). This allows us
to use the estimates \eqref{eq_upbd2ptfbd} and \eqref{eq_upbddisj}. This yields
the following extension of \eqref{3.4.energy} \cite{PfisterVelenik97}
\begin{equation}\label{eq_energybd}
\dIsbd(\frS) \leq \exp\{-\intSTbd(\frS)\}\,.
\end{equation}
The analogue of the energy estimate \eqref{3.4.energyge} then follows easily,
since $\taubd(\gb,\bdf)\geq 0$ when $\bdf\geq 0$ and therefore it is still
possible to control the number of vertices of $\frS$ in terms of
$\intSTbd(\frS)$. This gives
\begin{equation}\label{eq_energygebd}
\dIsbd(\frS \geq r) \leq \exp\Bigl\{-r(1-\frac{C\log N}{s(N)})\Bigr\}\,.
\end{equation}
Using this and the estimates in the phase of small contours, which still hold in
the presence of a boundary field, the upper bound follows easily.

\medskip\noindent
{\it Second case: $\bdf< 0$}

As for the lower bound, we have to deal with the fact that, for $\bdf<0$, the
duality is not defined. The solution is the same as there: We just change
boundary conditions, i.e. we look at the measure $\Isbdf{\pm}{N}{\abs\bdf}$,
which was defined when we dealt with the lower bound.

Once we have done this, the main difference is that the family of
low-temperature contours of any configurations compatible with these boundary
conditions contains exactly one open contour, with endpoints $t_{\rm l} =
(-N-\tfrac12, -\tfrac12)$ and $t_{\rm r} = (N+\tfrac12, \tfrac12)$. It is
straightforward to generalize the notion of skeleton introduced in the
preceding case to the present situation. What we get by this procedure is a
family of skeletons $\frS^\pm=(S_0,S_1,\dots,S_n)$ containing exactly one
skeleton, $S_0$, with $\Pol(S_0)$ open with endpoints $t_{\rm l}$ and $t_{\rm
r}$.

\begin{figure}[t]
\centerline{
\psfig{file=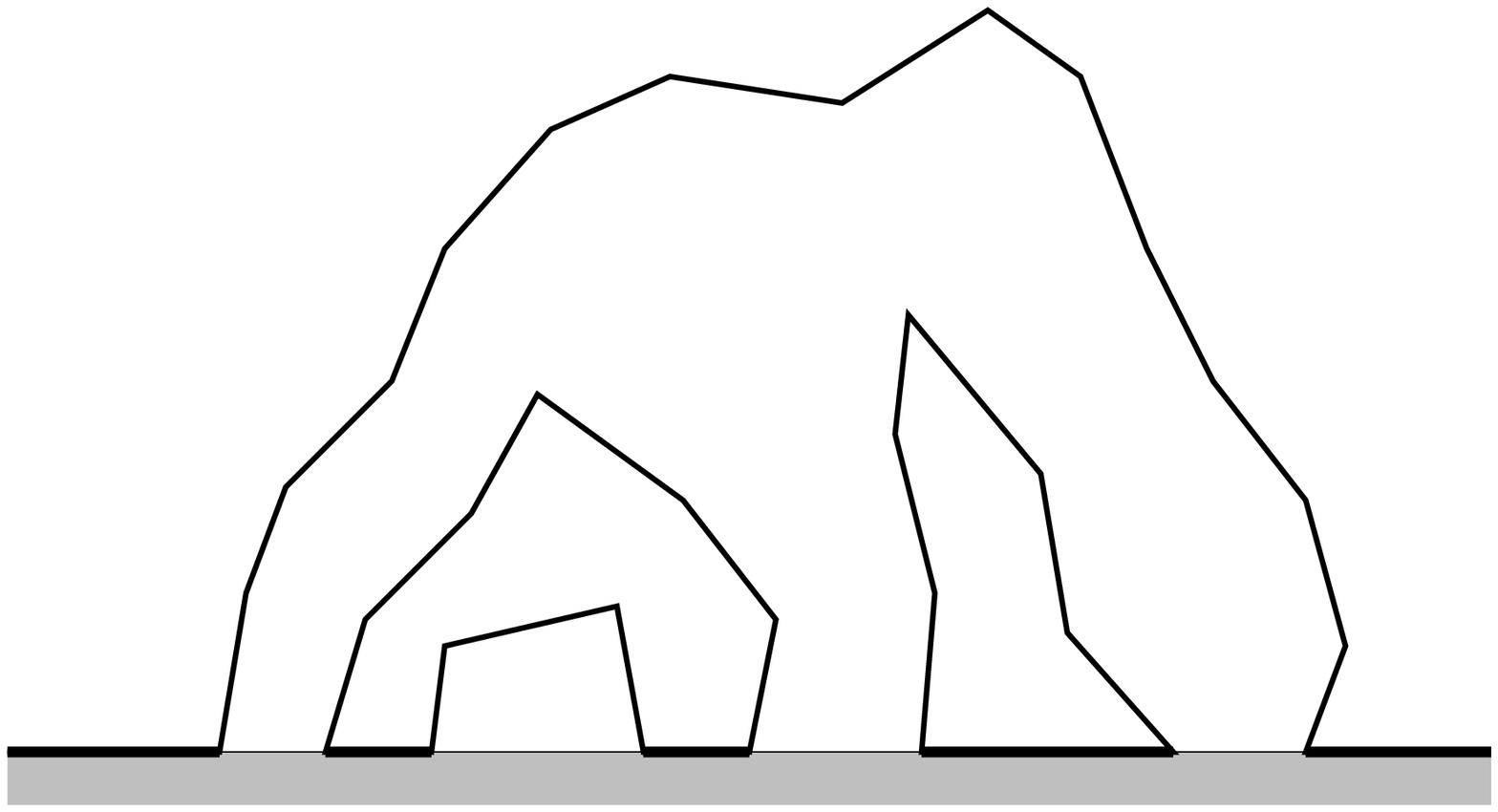,height=3.7cm}\hspace{1cm}
\psfig{file=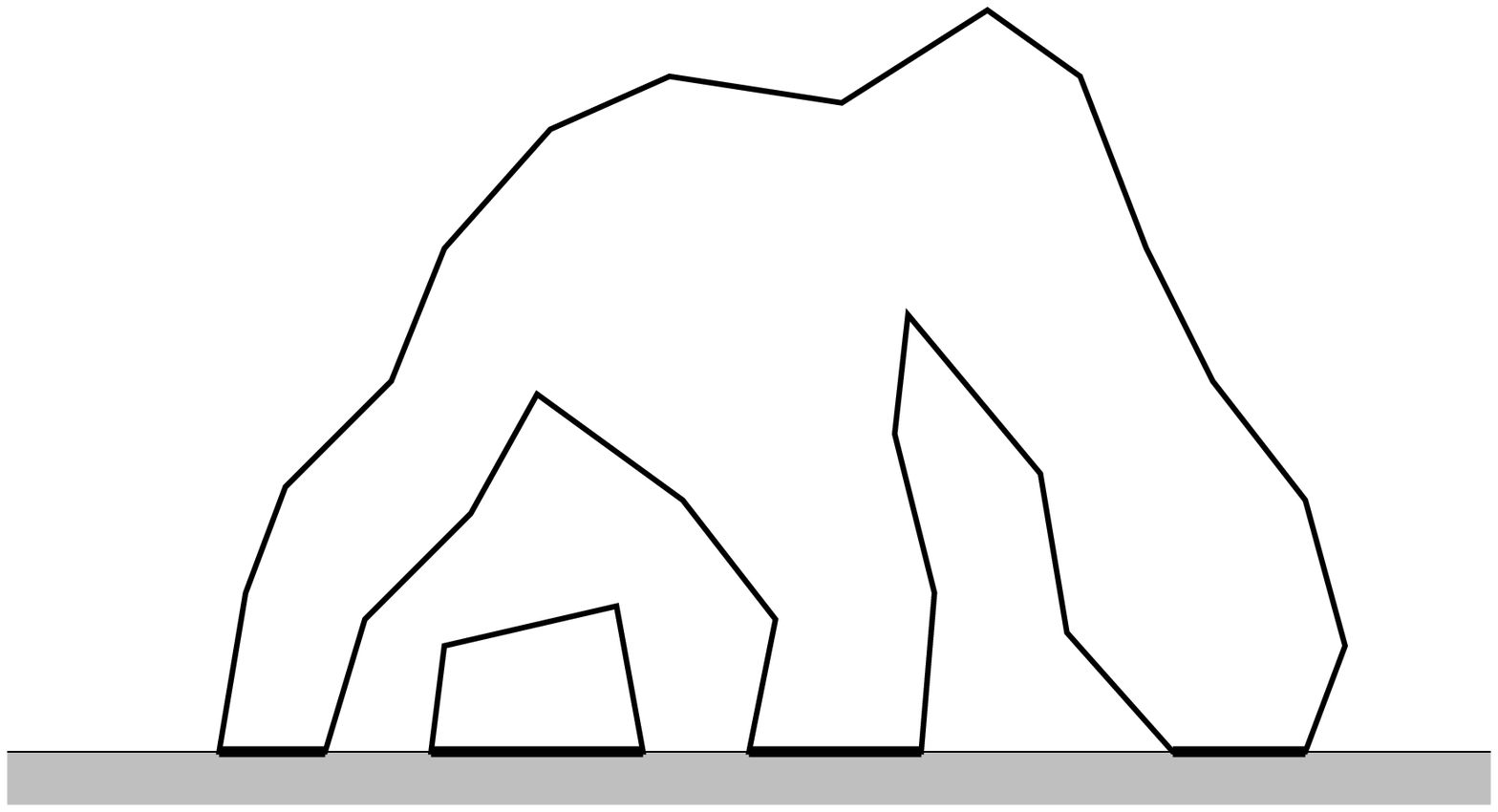,height=3.7cm}
}
\caption{Left: The family of polygonal lines associated to $\frS^\pm$. Right:
The family of {\em closed} polygonal lines associated to $\frS$.} 
\label{fig_bd_ub_Spm}
\end{figure} 
Since we want to compare the corresponding families of polygonal lines with the
solution of the variational problem, i.e. with the boundary of a convex body in
$\rBoxTwoD$, it is convenient to introduce another family $\frS$ of skeletons whose
associated polygonal lines are {\em closed}; $\frS$ possesses the same set of
vertices (except for $t_{\rm l}$ and $t_{\rm r}$, but with a different set of
edges, which is such that its associated family of polygonal lines satisfies
\begin{equation*}
\Pol(\frS) = \Pol(\frS^\pm) {\scriptstyle\triangle}
\setof{x\in\bbR^2}{-N/2-\tfrac12\leq x(1) \leq N/2+\tfrac12,\, x(2)=-\tfrac12}
\end{equation*}
where $\scriptstyle\triangle$ denotes symmetric difference (see
Fig.~\ref{fig_bd_ub_Spm}).

One then has the following relation
\begin{equation*}
\intSTbd(\frS) = \intSTbdf{\abs\bdf}(\frS^\pm) - (2N+1)\; \taubd(\gb,\abs\bdf)\,.
\end{equation*}
In particular, the following version of \eqref{eq_energybd} holds \cite{PfisterVelenik97}
\begin{align}
\Isbdf{\pm}{N}{\abs\bdf}(\frS^\pm) &\leq K_1\;
\exp\{-\intSTbd(\frS)\} &\bdf&> -\hw(\gb) \label{eq_energybd2}\\
\Isbdf{\pm}{N}{\abs\bdf}(\frS^\pm) &\leq K_2N^{3/2}\;
\exp\{-\intSTbd(\frS)\} &\bdf&\leq -\hw(\gb) \label{eq_energybd3}
\end{align}
The energy estimate \eqref{eq_energygebd} is slightly more delicate now, since
the wall free energy is negative. It turns out however that in the partial
wetting regime, $\bdf>-\hw(\gb)$, it is easy to reduce ourselves to a situation
similar to the case $\taubd(\gb,\bdf)\geq 0$. The case $\bdf\leq -\hw(\gb)$,
i.e. complete wetting, is more subtle, but happens not to give too much
problems as long as we consider volume-order large deviations (or, in fact,
deviations close enough to volume order).

Let us first consider the case of partial wetting; this regime is characterized
by $\abs{\taubd(\gb,\bdf)}<\tau^*_\gb$. Let us write $\intSTbd(\frS) =
T^++T^-$, where $T^+$ ($T^-$) is the positive (negative) part of the
functional. Then, since $T^+\geq (\tau^*_\gb/\taubd(\gb,\bdf))\,T^-$ and the
number of vertices along the wall is at most two-third of the total number
$\#(\frS)$, we have
\begin{equation*}
\#(\frS) \leq \frac{K}{s(N)(\tau^*_\gb+\taubd(\gb,\bdf))}\;\intSTbd(\frS)\,,
\end{equation*}
for some absolute constant $K$. This allows to prove that
\begin{equation}\label{eq_energygepw}
\Isbd{-}{N}(\frS \geq r) \leq \exp\Bigl\{-r(1-\frac{C\log N}{s(N)})\Bigr\}\,.
\end{equation}
When $\bdf\leq-\hw(\gb)$, one cannot establish so good an upper bound. The
best we can do is to use the fact that $T^-\geq (2N+1)\,\taubd(\gb,\bdf)$,
which turns out to be enough to prove the following, weaker, version of the
energy estimate
\begin{equation}\label{eq_energygecw}
\Isbd{-}{N}(\frS \geq r) \leq \exp\Bigl\{-r(1-\frac{C\log N}{s(N)})+C'\frac{N\log
N}{s(N)}\Bigr\}\,.
\end{equation}
The reason why such an estimate is still sufficient to get the desired result
is that the relevant values of $r$ are also of order $N$, so that the first
term can always be made to dominate the second one.

\bigskip
Once we have \eqref{eq_energygepw} and \eqref{eq_energygecw}, the proof is
concluded as usual, after observing that the estimate in the phase of small
contours still applies in the presence of the boundary field $\abs\bdf$.
\subsection{Ising model in $D\geq 3$}

The proof of Theorem~\ref{thm_sessile3D} is based on the $\bbL_1$-Theory 
introduced in Part~\ref{part_weakWulff}. 
We simply explain how the main ingredients of the proof 
should be modified
and refer to \cite{BodineauIoffeVelenik99} for details.\\

The arguments of geometric measure Theory can be extended easily to this
new setting. In particular, it is straightforward to check that the 
functional $\intSTbd$ is lower semi-continuous and that 
the approximation Theorems~\ref{thm approx} and \ref{theo ABCP} hold.

The main problem is to define proper mesoscopic phase labels 
for the measures with a boundary magnetic field. 
If $\eta \geq 0$, then the mesoscopic phase labels introduced in  
Part~\ref{part_weakWulff} satisfy the Assumptions A and B, as well as 
Conditions C1-C3 under the measure $\Isbd{+}{N}$.
Instead if $\eta < 0$, some problems occur because the FK measure
looses its ferromagnetic properties and the random coloring measures
are more complicated to deal with. 
Nevertheless, it is still possible to define mesoscopic phase labels
and to derive estimates as in Section~\ref{Coarse graining and mesoscopic phase
labels}.
 
Other difficulties have to be overcomed in order to implement the general 
philosophy of the $\bbL_1$-Theory.
In the case of a negative boundary magnetic field, the interface induced
by the field prevents us from applying directly the techniques developed 
to prove the exponential tightness Theorem \ref{prop 1}. 
Therefore an alternative approach similar to the one described in
Subsection \ref{ssec_2D} is required.
The analysis of the surface tension needs also some care.
We recall that the computation of surface tension is based on a localization 
procedure along the boundary of functions of bounded variation.
For a given test function either locally its boundary is in the bulk
and we recover the usual surface tension term or it
intersects the wall and arguments similar to those used in the bulk 
enable us to derive the wall free energy.
In this way the complexity of the problem is reduced because
the difficult analysis of the fluctuations of the microscopic interface 
between the wall and the bulk is replaced by soft $\bbL_1$ estimates.\\
\section{Open problems}
\setcounter{equation}{0}
As in the previous parts, there are still a lot of open problems. Most of those
presented before have natural analogues in the present situation. In the
following, we restrict ourselves to problems intrinsically related to the topics
discussed in this part.
\subsubsection{2D nearest-neighbors Ising model}
The fact that one is still unable to analyze non-pertur\-ba\-tively the
fluctuations of the phase separation line is only strengthened when we would
like to study boundary effects. Indeed, a general analysis of typical open paths
with endpoints at general positions with respect to the wall has not been done
even at low temperature. Problems related to this are the following:
\begin{enumerate}
\item Give a non-perturbative proof that the probability measure of a suitably
rescaled version of an open contour with endpoints on the wall converges weakly
to the measure of Brownian excursion when $\bdf\leq-\hw(\gb)$ (as was sketched in
the low-temperature case for $\bdf=-1$ in \cite{Dobrushin93}). This would
provide a way of analyzing the typical fluctuations of magnetization in the
complete wetting regime, and would complete the heuristic picture of the
wetting transition in the Grand-Canonical Ensemble.
\item Establish Ornstein-Zernike behavior for the boundary 2-point function
without having recourse to explicit computations. Even weaker lower bounds,
like those given in \cite{Al}, have not been proved in such a
constrained geometry.
\end{enumerate}
Another open problem is to investigate the full range of moderate deviations.
This may require an understanding of point 1. above.
\subsubsection{Higher dimensional nearest-neighbors Ising models}
If fluctuations of phase separation lines are not yet understood, the
situation is only much worse when considering their higher dimensional
counterparts; in fact, even perturbative results are not always available. Here
is a far from exhaustive list of related open problems.
\begin{enumerate}
\item Give a microscopic description of the behavior of phase boundaries in
the partial and complete wetting regimes in the Grand-Canonical Ensemble to put
some flesh on the heuristics given above.
\item Decide whether $\hw(\gb)=1$ or not. The corresponding results for the
SOS model \cite{Chalker} suggest that $\hw(\gb)<1$ in any dimension; numerical 
investigations confirm this in dimension~3 \cite{BinderLandau}.
\end{enumerate}
In fact, even much simpler problems related to behavior of higher dimensional
interfaces are still open: proof of the existence of a roughening transition in
$d=3$, proof of the unstability of the $(1,1,1)$ interface, ...

In some simpler models of the SOS type some (but not all!) of these problems can
be solved, but this does not seem to help in solving the original ones.
\subsubsection{The wall}
Another type of problems concerns properties of the wall. In particular, it
might be interesting to answer the following questions.
\begin{enumerate}
\item What happens if the interaction with the wall is more complicated (say,
non-nearest neighbor).
\item What happens if the boundary field is not homogeneous (for example, is a
``random'' configuration of $\bdf_1$ and $\bdf_2$ macroscopically equivalent to
some well-chosen homogeneous boundary field $\bdf=\overline\bdf$?).
\end{enumerate}

\part{Appendix}
\label{appendix}
\setcounter{section}{0}

\section{Appendix A : Proof of Theorem~\ref{thm Compactness}}
\label{appendix A}

Assumption~A  controls the number of zero $u_k$-blocks, whereas Assumption~B is
used to control the geometry of the mesoscopic phase labels. The dependence of $k_0$ on 
$\gd$ could be described as follows: we choose $k_0$ so large that
\begin{equation}
\label{knull}
\gr_k ~\leq ~\frac1{C(d)}\gd\qquad\text{for every}\ k\geq k_0,
\end{equation}
where $C(d)$ is a large enough fixed constant.  Three terms on the left hand side on 
\eqref{tightbound} correspond to three different exponential estimates:
\vskip 0.2cm
\noindent
\subsection{Estimate on the volume of zero $u_k$-blocks.}
The domination by Bernoulli measure \eqref{A} implies that
\begin{equation}
\label{volume}
\bbP_N 
\left( \#\{x\in\sTor{n-k}:u_k (x)=0\}\geq \gd\left(\frac{N}{2^k}\right)^d
\right)~\leq~c_2 \; \text{exp}\left\{ -\gd\left(\frac{N}{2^k}\right)^d\log
\frac{\gd}{\gr_k}\right\} .
\end{equation}

Each realization of the phase label function $u_k$ splits $\uTor$ into the disjoint union of
three mesoscopic regions:
$$
\uTor~=~\{x:u_k (x)=1\}\vee\{x:u_k (x)=-1\}\vee\{x:u_k (x)=0\} ~\df~
{\bf A}_+\vee {\bf A}_-\vee {\bf A}_0 .
$$
By the choice of the scale $k_0$ in \eqref{knull} the estimate \eqref{volume} is
 non-trivial for every $k\geq k_0$, and, in view of the target claim \eqref{tightbound}, we
can restrict attention only to such realizations of $u_k$ for which
\begin{equation}
\label{A_0}
\left| {\bf A}_0\right|~=~\int_{\uTor} 1_{\{u_k(x)=0\}}\text{d}x ~< ~\gd .
\end{equation}
This has the following important implication: if $u_k\in \cV\left( K_a ,2\gd\right)^{\text{c}}$,
 the area of the
boundary of any regular set $A$ such that ${\bf A}_+\subseteq A\subseteq\uTor\setminus {\bf A}_-$ is 
bounded below as
\begin{equation}
\label{partialA}
\left|\partial A\right|~\geq ~a .
\end{equation}
Using the Assumption~B of the Theorem we are going to construct such sets $A$ on the
finite $k_0$ scale; $A\in\cF_{n-k_0}$, and in such a fashion that all the boundary $k_0$-blocks
of $A$ will necessarily have zero $u_{k_0}$-labels. This reduction enables a uniform 
treatment of all coarser scales $k\geq k_0$. 

So let $k\geq k_0$, and assume that \eqref{A_0} holds.  
We denote by $A_-$ (resp. $A_+$) the set of all boxes $\sBox{n-k_0}$
in ${\bf A}_-$ (resp ${\bf A}_+$).
We say that $x\in\sTor{n-k_0}$ is $-*$~connected to $A_-$; 
$x\stackrel{-*}{\longleftrightarrow}A_-$, if there exists a $*$-connected 
chain of ``$-$'' $u_{k_0}$ blocks leading from $\sBox{n-k_0}(x)$ (and including it) to
$A_-$. Define now the complement $A^{\text{c}}$ as follows:
$$
A^{\text{c}}~=~A_-\bigcup_{x\stackrel{-*}{\longleftrightarrow}A_-}\sBox{n-k_0}(x) .
$$
By the virtue of the Assumption~B, ${\bf A}_+\subseteq A$. Moreover, 
by construction all
the $k_0$-blocks of $A$ attached to the boundary $\partial A^c$ have zero $u_{k_0}$-labels.
With a slight abuse of notation we proceed to denote this collection of boundary $k_0$-blocks
 as $\partial A$. By \eqref{partialA} the number of $k_0$-blocks in $\partial A$ is bounded
below by 
\begin{equation}
\label{kpartialA}
\#_{k_0}\left(\partial A\right)~\geq ~\frac{c(d)a}{2^{(d-1)k_0}}N^{d-1} .
\end{equation}
Since, however, the total number of $k_0$-blocks in the corresponding decomposition of $\uTor$
equals to $N^d /2^{dk_0}$ the estimate \eqref{kpartialA} alone is not sufficient for 
giving the desirable upper bound on the probability 
$\Joint_N \left( u_k\in\cV(K_a ,2\gd)^{\text{c}}\right)$.  The required 
entropy cancelation stems
from the fact that small connected contours of $\partial A$ 
cannot surround too much volume.

Let us decompose $A$ to the disjoint union of its maximal connected components:
$$
A~=~\bigvee_{i=1}^l A_i\qquad\text{respectively}\qquad 
\partial A~=~\bigvee_{i=1}^l \partial A_i .
$$
We shall quantify contours $\partial A_i$ according to
 the size (or the number of $k_0$-blocks ) in
 $A_i$.  Namely, the contour $\partial A_i$ is called small, if 
\begin{equation}
\label{Kd}
\#_{k_0}\left( A_i\right)~\leq ~K(d)\log N\qquad\text{or}
\qquad \left| A_i\right|~
\leq ~K(d)\frac{2^{dk_0}}{N^d}\log N,
\end{equation}
where $K(d)$ is a sufficiently large constant. Otherwise,
 the contour $\partial A_i$ is called large.

We claim that under \eqref{A_0} the following inclusion is valid:
\begin{equation}
\label{inclusion}
\left\{u_k\in\cV(K_a ,2\gd)^{\text{c}}\right\}~
\subseteq~\left\{\sum_{\partial A_i-\text{small}}|A_i | >\gd\right\}\bigcup
\left\{\sum_{\partial A_i-\text{large}} |\partial A_i | > a
\right\} .
\end{equation}
Indeed, if the total volume inside small contours is less
 than $\gd$, then repainting all the small
components $A_i$ into ``$-1$'' and all the large components $A_j$ 
into ``$+1$''  we produce a $\{\pm 1\}$-valued function 
 which is at most at the $\bbL_1$-distance $2\gd$ from $u_k$ and which,
 thereby, cannot belong to $K_a$.
\vskip 0.2cm
\noindent
\subsection{Peierls estimate on the size of large contours.}
\begin{equation}
\label{peierls}
\begin{split}
\Joint_N \left( \sum_{\partial A_i-\text{large}}| \partial A_i | >a\right)
~&=~ \Joint_N \left( \sum_{\partial A_i-\text{large}}\#_{k_0}(\partial A_i)> 
\frac{c(d)a}{2^{(d-1)k_0}}N^{d-1}
\right) \\
&\ \ \leq ~\text{exp}\left\{ -c_3 (d)\frac{a}{2^{(d-1)k_0}}N^{d-1}\right\} .
\end{split}
\end{equation}
This immediately follows from  Assumption~A, once the constant $K(d)$ in \eqref{Kd} has been
properly chosen.
\vskip 0.2cm
\noindent
\subsection{Estimate in the phase of small contours.}
The volume of small components
$A_i$ is related to the total number of $k_0$-blocks in these components as
$$
 \sum_{\partial A_i-\text{small}}|A_i |  
~=~\left(\frac{N}{2^{k_0}}\right)^{-d}
\sum_{\partial A_i-\text{small}}\#_{k_0}(A_i  ) .
$$
On the other hand, for every $l\in [1,...,n-k_0 ]$;
\begin{equation*}
\begin{split}
\sum_{\partial A_i-\text{small}}\#_{k_0}(A_i  )~&=~
\sum_{x\in\sTor{n-k_0}}\sum_{\partial A_i-\text{small}} 1_{\{x\in A_i\}}\\
&=~\sum_{t\in [0,...,2^l)^d}
\sum_{x\in\sTor{n-k_0-l}}\sum_{\partial A_i-\text{small}} 1_{\{\theta_{t\gD_0}x\in A_i\}} ,
\end{split}
\end{equation*}
where 
 $\gD_0\df2^{k_0 -n}$ is the step size on the embedded torus $\sTor{n-k_0}$, and  
$\theta_{\bullet}$ is the shift on this torus. Consequently,
\begin{equation}
\label{lsplit}
\Joint_N \left(\sum_{\partial A_i-\text{small}}|A_i |  >\gd\right)~\leq~\max_{t\in [0,...,2^l)^d}
\Joint_N \left(\sum_{x\in\sTor{n-k_0-l}}
\sum_{\partial A_i-\text{small}} 1_{\{\theta_{t\gD_0}x\in A_i\}} >
\gd\left(\frac{N}{2^{k_0 +l}}
\right)^d\right) .
\end{equation}
If, however, $2^l > K(d)\log N$, then no two distinct points on the torus $\sTor{n-k_0 -l}$ (or
 any shift of it) can belong to the same small component $A_i$.  This, in view of the domination
by the independent Bernoulli site percolation (Assumption~A), suggests an application of the 
B-K inequality. Since, by the choice of the scale $k_0$ in \eqref{knull};
$$
\epsilon_{k_0}~\df ~\bbP_{\text{perc}}^{\gr_{k_0}}\left(\exists~\text{a closed surface of zero
$u_{k_0}$-blocks around $x$}\right)~< ~\gd ,
$$
for every $x\in\sTor{n-k_0}$, we readily obtain that the right hand side of \eqref{lsplit} is bounded
above by
$$
c_4 (d)\text{exp}\left\{ -\gd\left(\frac{N}{2^{k_0 +l}}\right)^d
\log\left(\frac{\gd}{\epsilon_{k_0}}\right)\right\} .
$$
The proof of Theorem \ref{thm Compactness} is concluded.
\qed
\section{Appendix B : Proof of the three-point lower bound Lemma~\ref{triple}} 
\label{appendix B}
The proof of Lemma~\ref{triple} is based on the following positive
 stiffness property of the surface tension \cite{AkutsuAkutsu86}:
\begin{equation}
\label{B_stiff}
\min_{\theta\in[0,2\pi]}\left\{
\frac{{\rm d}^2}{{\rm d}\theta^2}
\st \left(\vec{n}(\theta )\right)~+~\st \left(\vec{n}(\theta )\right)
\right\}\ =\ \min_{\theta\in[0,2\pi]}R_{\gb}\left(
\vec{n}(\theta )\right)\ >\ 0.
\end{equation}
 where the unit normal $\vec{n}(\theta )$ is defined via 
$\vec{n}(\theta ) = (\cos\theta ,\sin \theta )$, and $R_{\gb}\left(
\vec{n}\right)$ is the radius of curvature of $\partial \cK $ at the point 
 supporting the tangent line orthogonal to $\vec{n}$. An integral 
version of \eqref{B_stiff} is the strong triangle inequality \cite{I1},
\cite{
Velenik97}: For any $u,v\in\bbR^2$:  
\begin{equation}
\label{B_strong}
\st\left( u\right)+
\st\left( v\right)-
\st\left( u+v\right)~\geq ~ c_1(\gb )\left(
\| u\|_2+
\| v\|_2 -\| u +v\|_2\right) .
\end{equation}
The latter inequality is used to control the fluctuations of the 
microscopic phase boundaries (in their random line representation of
 Section \ref{dima_skeletons}).

Let now an $(s,\gep )$-compatible triple of points $(u,w,v)$ be given.
Fix $K =K(\gb )$ large enough and define the ``oval'' neighborhood
 ${\bf N}_K (u,w )$ of $\{u,v\}$ as:
$$
{\bf N}_K (u,w ) ~\df ~\left\{ z\in\bbR^2 :\ 
\st\left(z- u\right)+
\st\left( w-z\right)-
\st\left( w-u\right)\leq K\log s\right\} .
$$
The oval neighborhood ${\bf N}_K (w,v )$ is defined exactly in the same 
fashion. Relations \eqref{3.4.qasub} and \eqref{3.4.OZ} readily imply that
that the main contribution to 
$\langle\gs_u\gs_w\rangle_{f}^{\gb^*}$ (respectively to 
$\langle\gs_w\gs_v\rangle_{f}^{\gb^*}$ ) comes from the 
paths $\gl_1$ (respectively $\gl_2$ ) which stay in 
 ${\bf N}_K (u,w )$ (respectively  ${\bf N}_K (w,v )$). More precisely,
\begin{equation}\label{eq_concentration}
\sumtwo{\gl_1 :u\to w}{\gl\in{\bf N}_K (u,w )}
q^{\gb^*}\left( \gl_1\right)\geq
\langle\gs_u\gs_w\rangle_{f}^{\gb^*}\left( 1+{\rm \small{o}}(1)\right) ,
\end{equation}
uniformly in all $(s,\gep )$-compatible triples. Any such path
 $\gl_1 =\left(\gl_1 (0), ..., \gl_1 (n_1)\right)$ could
 be decomposed as follows:
 Define 
$$
n_w~=~\max\left\{ k: \gl_k \in {\bf N}_K (u,w )\setminus {\bf N}_K (w,v )
\right\} ,
$$
and set $\gl_1^u = \left(\gl_1 (0) ,...,\gl_1 (n_w )\right)$, 
$\gl_1^w = \left(\gl_1 (n_w +1) ,...,\gl_1 (n_1)\right)$; 
$\gl_1 =\gl_1^u\vee\gl_1^w$.
 The decomposition  $\gl_2 =\gl_2^u\vee\gl_2^w$ is defined in a completely
symmetric way. Notice that, by the construction,
 the paths $\gl_1^u$ and $\gl_2^v$ are disjoint and compatible, and,
 by \eqref{B_strong};
$$
\max\left\{ \| \gl_1 (n^w )-w\|_2 , \| \gl_2 (n^w )-w\|_2\right\}~\leq
~c_2 (\gep )\log s .
$$
The claim of the lemma follows now from \eqref{3.4.qasuper} and 
\eqref{3.4.expbouns}.\qed

\end{document}